
\hsize=13.50cm    
\vsize=18cm       
\parindent=12pt   \parskip=0pt      
\pageno=1 


\hoffset=15mm    
\voffset=1cm    
 

\ifnum\mag=\magstep1
\hoffset=-2mm   
\voffset=.8cm   
\fi


\pretolerance=500 \tolerance=1000  \brokenpenalty=5000

\catcode`\@=11

\font\eightrm=cmr8         \font\eighti=cmmi8
\font\eightsy=cmsy8        \font\eightbf=cmbx8
\font\eighttt=cmtt8        \font\eightit=cmti8
\font\eightsl=cmsl8        \font\sixrm=cmr6
\font\sixi=cmmi6           \font\sixsy=cmsy6
\font\sixbf=cmbx6


\font\tengoth=eufm10       \font\tenbboard=msbm10
\font\eightgoth=eufm10 at 8pt      \font\eightbboard=msbm10 at 8 pt
\font\sevengoth=eufm7      \font\sevenbboard=msbm7
\font\sixgoth=eufm7 at 6 pt        \font\fivegoth=eufm5

 \font\tencyr=wncyr10       
\font\eightcyr=wncyr10 at 8 pt      
\font\sevencyr=wncyr10 at 7 pt      
\font\sixcyr=wncyr10 at 6 pt


\skewchar\eighti='177 \skewchar\sixi='177
\skewchar\eightsy='60 \skewchar\sixsy='60


\newfam\gothfam           \newfam\bboardfam
\newfam\cyrfam

\def\tenpoint{%
  \textfont0=\tenrm \scriptfont0=\sevenrm \scriptscriptfont0=\fiverm
  \def\rm{\fam\z@\tenrm}%
  \textfont1=\teni  \scriptfont1=\seveni  \scriptscriptfont1=\fivei
  \def\oldstyle{\fam\@ne\teni}\let\old=\oldstyle
  \textfont2=\tensy \scriptfont2=\sevensy \scriptscriptfont2=\fivesy
  \textfont\gothfam=\tengoth \scriptfont\gothfam=\sevengoth
  \scriptscriptfont\gothfam=\fivegoth
  \def\goth{\fam\gothfam\tengoth}%
  \textfont\bboardfam=\tenbboard \scriptfont\bboardfam=\sevenbboard
  \scriptscriptfont\bboardfam=\sevenbboard
  \def\bb{\fam\bboardfam\tenbboard}%
 \textfont\cyrfam=\tencyr \scriptfont\cyrfam=\sevencyr
  \scriptscriptfont\cyrfam=\sixcyr
  \def\cyr{\fam\cyrfam\tencyr}%
  \textfont\itfam=\tenit
  \def\it{\fam\itfam\tenit}%
  \textfont\slfam=\tensl
  \def\sl{\fam\slfam\tensl}%
  \textfont\bffam=\tenbf \scriptfont\bffam=\sevenbf
  \scriptscriptfont\bffam=\fivebf
  \def\bf{\fam\bffam\tenbf}%
  \textfont\ttfam=\tentt
  \def\tt{\fam\ttfam\tentt}%
  \abovedisplayskip=12pt plus 3pt minus 9pt
  \belowdisplayskip=\abovedisplayskip
  \abovedisplayshortskip=0pt plus 3pt
  \belowdisplayshortskip=4pt plus 3pt 
  \smallskipamount=3pt plus 1pt minus 1pt
  \medskipamount=6pt plus 2pt minus 2pt
  \bigskipamount=12pt plus 4pt minus 4pt
  \normalbaselineskip=12pt
  \setbox\strutbox=\hbox{\vrule height8.5pt depth3.5pt width0pt}%
  \let\bigf@nt=\tenrm       \let\smallf@nt=\sevenrm
  \normalbaselines\rm}

\def\eightpoint{%
  \textfont0=\eightrm \scriptfont0=\sixrm \scriptscriptfont0=\fiverm
  \def\rm{\fam\z@\eightrm}%
  \textfont1=\eighti  \scriptfont1=\sixi  \scriptscriptfont1=\fivei
  \def\oldstyle{\fam\@ne\eighti}\let\old=\oldstyle
  \textfont2=\eightsy \scriptfont2=\sixsy \scriptscriptfont2=\fivesy
  \textfont\gothfam=\eightgoth \scriptfont\gothfam=\sixgoth
  \scriptscriptfont\gothfam=\fivegoth
  \def\goth{\fam\gothfam\eightgoth}%
  \textfont\cyrfam=\eightcyr \scriptfont\cyrfam=\sixcyr
  \scriptscriptfont\cyrfam=\sixcyr
  \def\cyr{\fam\cyrfam\eightcyr}%
  \textfont\bboardfam=\eightbboard \scriptfont\bboardfam=\sevenbboard
  \scriptscriptfont\bboardfam=\sevenbboard
  \def\bb{\fam\bboardfam}%
  \textfont\itfam=\eightit
  \def\it{\fam\itfam\eightit}%
  \textfont\slfam=\eightsl
  \def\sl{\fam\slfam\eightsl}%
  \textfont\bffam=\eightbf \scriptfont\bffam=\sixbf
  \scriptscriptfont\bffam=\fivebf
  \def\bf{\fam\bffam\eightbf}%
  \textfont\ttfam=\eighttt
  \def\tt{\fam\ttfam\eighttt}%
  \abovedisplayskip=9pt plus 3pt minus 9pt
  \belowdisplayskip=\abovedisplayskip
  \abovedisplayshortskip=0pt plus 3pt
  \belowdisplayshortskip=3pt plus 3pt 
  \smallskipamount=2pt plus 1pt minus 1pt
  \medskipamount=4pt plus 2pt minus 1pt
  \bigskipamount=9pt plus 3pt minus 3pt
  \normalbaselineskip=9pt
  \setbox\strutbox=\hbox{\vrule height7pt depth2pt width0pt}%
  \let\bigf@nt=\eightrm     \let\smallf@nt=\sixrm
  \normalbaselines\rm}

\tenpoint


\def\pc#1{\bigf@nt#1\smallf@nt}         \def\pd#1 {{\pc#1} }


\catcode`\;=\active
\def;{\relax\ifhmode\ifdim\lastskip>\z@\unskip\fi
\kern\fontdimen2 \font\kern -1.2 \fontdimen3 \string;}

\catcode`\:=\active
\def:{\relax\ifhmode\ifdim\lastskip>\z@\unskip\fi\penalty\@M\ \fi\string:}

\catcode`\!=\active
\def!{\relax\ifhmode\ifdim\lastskip>\z@
\unskip\fi\kern\fontdimen2 \font \kern -1.1 \fontdimen3 \string!}

\catcode`\?=\active
\def?{\relax\ifhmode\ifdim\lastskip>\z@
\unskip\fi\kern\fontdimen2 \font \kern -1.1 \fontdimen3 \string?}

\def\^#1{\if#1i{\accent"5E\i}\else{\accent"5E #1}\fi}
\def\"#1{\if#1i{\accent"7F\i}\else{\accent"7F #1}\fi}

\frenchspacing


\newtoks\auteurcourant      \auteurcourant={\hfil}
\newtoks\titrecourant       \titrecourant={\hfil}

\newtoks\hautpagetitre      \hautpagetitre={\hfil}
\newtoks\baspagetitre       \baspagetitre={\hfil}

\newtoks\hautpagegauche     
\hautpagegauche={\eightpoint\rlap{\folio}\hfil\the\auteurcourant\hfil}
\newtoks\hautpagedroite     
\hautpagedroite={\eightpoint\hfil\the\titrecourant\hfil\llap{\folio}}

\newtoks\baspagegauche      \baspagegauche={\hfil} 
\newtoks\baspagedroite      \baspagedroite={\hfil}

\newif\ifpagetitre          \pagetitretrue  


\headline={\ifpagetitre\the\hautpagetitre
\else\ifodd\pageno\the\hautpagedroite\else\the\hautpagegauche\fi\fi}

\footline={\ifpagetitre\the\baspagetitre\else
\ifodd\pageno\the\baspagedroite\else\the\baspagegauche\fi\fi
\global\pagetitrefalse}


\def\raggedbottom{\topskip 10pt plus 36pt\r@ggedbottomtrue}



\def\pointir{\unskip . --- \ignorespaces}


\def\Bigbreak{\vskip-\lastskip\bigbreak}
\def\Medbreak{\vskip-\lastskip\medbreak}


\def\ctexte#1\endctexte{%
  \hbox{$\vcenter{\halign{\hfill##\hfill\crcr#1\crcr}}$}}


\long\def\ctitre#1\endctitre{%
    \ifdim\lastskip<24pt\vskip-\lastskip\bigbreak\bigbreak\fi
  		\vbox{\parindent=0pt\leftskip=0pt plus 1fill
          \rightskip=\leftskip
          \parfillskip=0pt\bf#1\par}
    \bigskip\nobreak}

\long\def\section#1\endsection{%
\vskip 0pt plus 3\normalbaselineskip
\penalty-250
\vskip 0pt plus -3\normalbaselineskip
\Bigbreak
\message{[section \string: #1]}{\bf#1\unskip}\pointir}

\long\def\sectiona#1\endsection{%
\vskip 0pt plus 3\normalbaselineskip
\penalty-250
\vskip 0pt plus -3\normalbaselineskip
\Bigbreak
\message{[sectiona \string: #1]}%
{\bf#1}\medskip\nobreak}

\long\def\subsection#1\endsubsection{%
\Medbreak
{\it#1\unskip}\pointir}

\long\def\subsectiona#1\endsubsection{%
\Medbreak
{\it#1}\par\nobreak}

\def\rem#1\endrem{%
\Medbreak
{\it#1\unskip} : }

\def\remp#1\endrem{%
\Medbreak
{\pc #1\unskip}\pointir}

\def\rema#1\endrem{%
\Medbreak
{\it #1}\par\nobreak}

\def\newparwithcolon#1\endnewparwithcolon{
\Medbreak
{#1\unskip} : }

\def\newparwithpointir#1\endnewparwithpointir{
\Medbreak
{#1\unskip}\pointir}

\def\newpara#1\endnewpar{
\Medbreak
{#1\unskip}\smallskip\nobreak}


\long\def\th#1 #2\enonce#3\endth{%
   \Medbreak
   {\pc#1} {#2\unskip}\pointir{\it #3}\medskip}

\long\def\tha#1 #2\enonce#3\endth{%
   \Medbreak
   {\pc#1} {#2\unskip}\par\nobreak{\it #3}\medskip}


\long\def\Th#1 #2 #3\enonce#4\endth{%
   \Medbreak
   #1 {\pc#2} {#3\unskip}\pointir{\it #4}\medskip}

\long\def\Tha#1 #2 #3\enonce#4\endth{%
   \Medbreak
   #1 {\pc#2} #3\par\nobreak{\it #4}\medskip}


\def\decale#1{\smallbreak\hskip 28pt\llap{#1}\kern 5pt}
\def\decaledecale#1{\smallbreak\hskip 34pt\llap{#1}\kern 5pt}
\def\puce{\smallbreak\hskip 6pt{$\scriptstyle\bullet$}\kern 5pt}



\def\displaylinesno#1{\displ@y\halign{
\hbox to\displaywidth{$\@lign\hfil\displaystyle##\hfil$}&
\llap{$##$}\crcr#1\crcr}}


\def\ldisplaylinesno#1{\displ@y\halign{ 
\hbox to\displaywidth{$\@lign\hfil\displaystyle##\hfil$}&
\kern-\displaywidth\rlap{$##$}\tabskip\displaywidth\crcr#1\crcr}}


\def\eqalign#1{\null\,\vcenter{\openup\jot\m@th\ialign{
\strut\hfil$\displaystyle{##}$&$\displaystyle{{}##}$\hfil
&&\quad\strut\hfil$\displaystyle{##}$&$\displaystyle{{}##}$\hfil
\crcr#1\crcr}}\,}


\def\system#1{\left\{\null\,\vcenter{\openup1\jot\m@th
\ialign{\strut$##$&\hfil$##$&$##$\hfil&&
        \enskip$##$\enskip&\hfil$##$&$##$\hfil\crcr#1\crcr}}\right.}


\let\@ldmessage=\message

\def\message#1{{\def\pc{\string\pc\space}%
                \def\'{\string'}\def\`{\string`}%
                \def\^{\string^}\def\"{\string"}%
                \@ldmessage{#1}}}



\def\up#1{\raise 1ex\hbox{\smallf@nt#1}}


\def\qed{\raise -2pt\hbox{\vrule\vbox to 10pt{\hrule width 4pt
                 \vfill\hrule}\vrule}}

\def\cqfd{\unskip\penalty 500\quad\qed\medbreak}

\def\virg{\raise .4ex\hbox{,}}   


\def\build#1_#2^#3{\mathrel{
\mathop{\kern 0pt#1}\limits_{#2}^{#3}}}


\def\boxit#1#2{%
\setbox1=\hbox{\kern#1{#2}\kern#1}%
\dimen1=\ht1 \advance\dimen1 by #1 \dimen2=\dp1 \advance\dimen2 by #1 
\setbox1=\hbox{\vrule height\dimen1 depth\dimen2\box1\vrule}%
\setbox1=\vbox{\hrule\box1\hrule}%
\advance\dimen1 by .6pt \ht1=\dimen1 
\advance\dimen2 by .6pt \dp1=\dimen2  \box1\relax}


\catcode`\@=12

\showboxbreadth=-1  \showboxdepth=-1



 \input amssym.def 
\input amssym.tex

\magnification=\magstep1
\hsize=17,5truecm   
\vsize=25.5truecm  
\hoffset=-0.9truecm  
\voffset=-0.8truecm
\topskip=1truecm
\footline={\tenrm\hfil\folio\hfil}
\raggedbottom
\abovedisplayskip=3mm 
\belowdisplayskip=3mm
\abovedisplayshortskip=0mm
\belowdisplayshortskip=2mm
\normalbaselineskip=12pt  
\normalbaselines
\def\noblackbox{\overfullrule=0pt}
\noblackbox

\def\et{{\acute et}}
\def\red{{r\acute ed}}

\def\g{\frak{g}}
\def\h{\frak{h}}
\def\C{{\bf C}}
\def\Q{{\bf Q}}
\def\Z{{\bf Z}}
\def\Br{{\rm Br}}
\def\br{{\rm Br}}
\def\Div{{\rm Div}}
\def\pic{{\rm Pic}}
\def\Pic{{\rm Pic}}

\def\G{{\bf G}}
\def\g{ \frak{g}  }
\def\Ker{{\rm Ker}}
\def\Coker{{\rm Coker}}
\def\ker{{\rm Ker}}
\def\Hom{{\rm Hom}}
\def\Ext{{\rm Ext}}
\def\X{{\cyr X}}
\def\k{{\overline k}}

\def\Grille{\setbox13=\vbox to 5\unitlength{\hrule width 109mm\vfill} 
\setbox13=\vbox to 65mm{\offinterlineskip\leaders\copy13\vfill\kern-1pt\hrule}
\setbox14=\hbox to 5\unitlength{\vrule height 65mm\hfill}
\setbox14=\hbox to 109mm{\leaders\copy14\hfill\kern-2mm\vrule height 65mm}
\ht14=0pt\dp14=0pt\wd14=0pt \setbox13=\vbox to
0pt{\vss\box13\offinterlineskip\box14} \wd13=0pt\box13}



\def\diagram#1{\def\normalbaselines{\baselineskip=0pt\lineskip=5pt}
\matrix{#1}}

\def\hfl#1#2#3{\smash{\mathop{\hbox to#3{\rightarrowfill}}\limits
^{\scriptstyle#1}_{\scriptstyle#2}}}

\def\gfl#1#2#3{\smash{\mathop{\hbox to#3{\leftarrowfill}}\limits
^{\scriptstyle#1}_{\scriptstyle#2}}}


{\bf R\'esolutions flasques des groupes lin\'eaires connexes} 


\bigskip

J.-L. Colliot-Th\'el\`ene

\medskip

{\bf Summary} {\it A connected reductive group $G$ over a
 field $k$  may be written as a quotient $H/S$, where the $k$-group
$H$ is an extension of a quasitrivial torus by a simply connected
semisimple group, and $S$ is a flasque $k$-torus, central in $H$.
The flasque torus $S$ is well-defined up to multiplication by a quasitrivial torus.
Such presentations $G=H/S$ lead to a simplified approach 
of the Galois cohomology of $G$ and of related objects,
such as the Brauer group of a smooth compactification of $G$.
 When $k$ is  a number field, one also recovers
 known formulas, in terms of $S$, for
the quotient of the group $G(k)$ of rational points by
$R$-equivalence, and for the abelian groups which measure
the lack of
 weak approximation and the failure of the
Hasse principle for principal homogeneous spaces.}

\bigskip

Sur un corps de nombres, une s\'erie de travaux
importants a \'etabli pour les groupes semi-simples
simplement connexes les propri\'et\'es suivantes.
Ils satisfont l'approximation faible,
leurs espaces  homog\`enes principaux satisfont le
principe de Hasse, leur nombre de Tamagawa est \'egal \`a 1.
Dans presque tous les cas, on a aussi \'etabli que la
$R$-\'equivalence sur les $k$-points d'un tel groupe
est triviale.

Ces propri\'et\'es ne valent  en g\'en\'eral pas pour 
un groupe r\'eductif connexe.
Pour mesurer le d\'efaut de ces propri\'et\'es, 
on utilise la cohomologie galoisienne. Ce genre de travail a \'et\'e fait
en particulier par Sansuc [Sa81],  Kottwitz [Ko84, Ko86], Borovoi 
[Bo96, Bo98], Gille [Gi97, Gi01]. Il implique l'utilisation de pr\'esentations convenables
des groupes lin\'eaires connexes.
On trouve dans la
litt\'erature deux types de pr\'esentations pour ces groupes, les {\it rev\^etements sp\'eciaux} et les $z$-{\it extensions}.

Soit $k$ un corps de caract\'eristique nulle.
Un rev\^etement sp\'ecial d'un $k$-groupe r\'eductif connexe $G$
est la donn\'ee d'une isog\'enie
$$ 1 \to \mu \to H \to G \to 1,$$
o\`u $H$ est le produit d'un $k$-groupe semi-simple connexe
par un $k$-tore quasi-trivial. Tout $k$-groupe
semi-simple connexe admet \'evidemment un tel rev\^etement.
On montre (voir [Sa81], Lemme 1.10) que pour tout $k$-groupe
r\'eductif connexe  $G$, il existe un entier $n>0$ et un $k$-tore quasi-trivial $P$
tel que le produit $G^n \times_{k} P$ admette un rev\^etement quasi-trivial.
Ces rev\^etements sp\'eciaux ont \'et\'e utilis\'es par Ono, Sansuc, Gille.

Par ailleurs tout $k$-groupe r\'eductif connexe $G$ admet une
$z$-extension, c'est-\`a-dire qu'il existe une
suite exacte
$$ 1 \to P \to H \to G \to 1,$$
o\`u $P$ est un $k$-tore quasi-trivial central
dans le $k$-groupe r\'eductif $H$, le groupe
d\'eriv\'e $H^{ss} \subset H$ de $H$ \'etant un $k$-groupe
(semi-simple) simplement connexe.
Les $z$-extensions ont \'et\'e utilis\'ees
par
Langlands, Deligne, Milne-Shih [Mi-Sh82], Kottwitz [Ko 84, Ko86], Borovoi [Bo96, Bo98],
Borovoi-Kunyavski\u{\i} [BoKu04].

Le but du pr\'esent texte, dont les principaux r\'esultats
ont \'et\'e annonc\'es dans la note [CT04], est de d\'evelopper
un  troisi\`eme type de pr\'esentation, introduit dans [CT04, CT05], et de montrer comment 
cela permet d'\'etudier de fa\c con plus souple  
la cohomologie galoisienne des groupes lin\'eaires.

Le point principal est que l'analogue des r\'esolutions
flasques de tores d\'evelopp\'ees il y a  trente ans (Sansuc et
l'auteur [CT/Sa77], Voskresenski\u{\i} [Vo77, Vo98]) existe pour tout $k$-groupe r\'eductif
connexe,  le r\^ole des tores
quasi-triviaux \'etant ici tenu par les groupes  extensions d'un
tore quasi-trivial par un groupe semi-simple simplement connexe :

Pour tout $k$-groupe
r\'eductif connexe $G$, il existe des suites exactes de $k$-groupes alg\'ebriques
$$1 \to S \to H \to G \to 1,$$
o\`u le groupe $H$ est une   extension d'un $k$-tore {\it quasi-trivial}
par le $k$-groupe semi-simple simplement connexe
rev\^etement universel du groupe d\'eriv\'e de $G$,
et o\`u le $k$-tore $S$, central dans $H$, est {\it flasque} (son groupe des cocaract\`eres est un module galoisien
$H^1$-trivial). Une telle suite exacte est appel\'ee une r\'esolution flasque de $G$. A multiplication par un
$k$-tore quasitrivial pr\`es, le $k$-tore
flasque $S$ est d\'etermin\'e par $G$.

Il y a aussi lieu d'introduire les r\'esolutions coflasques de $G$,
qui sont des suites exactes
$$1 \to P \to H \to G \to 1,$$
o\`u $P$ est un $k$-tore quasi-trivial  central
dans le $k$-groupe r\'eductif $H$, le groupe $H$
\'etant une extension d'un $k$-tore coflasque par
un $k$-groupe semi-simple simplement connexe.
Une extension de ce type est une $z$-extension, mais d'un type tr\`es particulier.

\medskip

Au paragraphe  0  on rappelle quelques propri\'et\'es de base des groupes alg\'ebriques lin\'eaires connexes. 
La d\'efinition, l'existence et les propri\'et\'es g\'en\'erales des r\'esolutions flasques et coflasques  font l'objet des paragraphes 1 \`a 4. 
 Au paragraphe 5, on voit comment ces r\'esolutions sont li\'ees
 aux torseurs universels sur les compactifications lisses des groupes lin\'eaires.
 Au paragraphe 6, on associe \`a tout $k$-groupe lin\'eaire connexe $G$
 un module galoisien discret de type fini, not\'e $\pi_{1}(G)$,  
 le groupe fondamental alg\'ebrique de $G$, et on en \'etablit les propri\'et\'es
 principales.
Notre d\'efinition de $\pi_{1}(G)$ passe par les r\'esolutions flasques. Elle  est ind\'ependante de la
d\'efinition originale de Borovoi ([Bo98],  \S 1)
et de celle de Merkur'ev ([Me98], \S 10).  
Au paragraphe 7, on  donne deux formules pour  le groupe de Brauer d'une compactification lisse
d'un $k$-groupe lin\'eaire connexe $G$, l'une en termes du tore flasque apparaissant dans une r\'esolution flasque de $G$, l'autre en termes du groupe fondamental alg\'ebrique de $G$, comme d\'efini au paragraphe  6
(avec le groupe fondamental alg\'ebrique de Borovoi, cette formule est d\'ej\`a \'etablie dans [BoKu00]).
Aux paragraphes 8 et 9, on montre comment une r\'esolution flasque 
d'un $k$-groupe lin\'eaire connexe $G$ donne de l'information d'une part sur
le groupe $G(k)/R$, quotient du groupe des points rationnels de $G$ par la $R$-\'equivalence ([CTSa77]),
d'autre part sur l'ensemble $H^1(k,G)$ classifiant les $G$-espaces  homog\`enes principaux.
On montre en particulier comment sur un corps local ou global  les r\'esolutions flasques permettent de d\'eduire de fa\c con rapide et fonctorielle la plupart des r\'esultats connus sur la cohomologie galoisienne des groupes lin\'eaires connexes  (Sansuc, Kottwitz, Borovoi, Gille), \`a partir des  r\'esultats fondamentaux sur les groupes semi-simples simplement connexes (Kneser, Harder, Chernousov).

Dans les travaux d'autres auteurs, en particulier ceux de Borovoi, un r\^ole important est jou\'e par
l'hypercohomologie de certains complexes de $k$-tores, de  longueur 2, li\'es au choix d'un tore maximal dans le groupe consid\'er\'e, choix qui est aussi fait par ces auteurs lors de la d\'efinition du groupe fondamental alg\'ebrique. L'approche par les r\'esolutions flasques propos\'ee ici \'evite l'utilisation de ces complexes, et (sauf au paragraphe  9) elle \'evite  le choix d'un tore maximal.
Cette approche est n\'eanmoins \'equivalente \`a celle des travaux cit\'es. C'est ce que l'on \'etablit
dans l'appendice A, o\`u l'on montre que les complexes de $k$-tores de Borovoi sont quasi-isomorphes 
\`a des complexes d\'eduits des r\'esolutions flasques.
On   montre ainsi que le groupe fondamental alg\'ebrique d'un $k$-groupe alg\'ebrique lin\'eaire ici d\'efini
est isomorphe \`a celui d\'efini par Borovoi. On fait aussi le lien avec les groupes
de cohomologie ab\'elianis\'ee de Borovoi [Bo96, Bo98]. La donn\'ee d'une compactification lisse
d'un groupe m\`ene  \`a la construction d'un autre complexe de $k$-tores, dont on montre dans 
l'appendice B qu'il est quasi-\'equivalent aux pr\'ec\'edents. 
Cela permet de faire le lien avec un travail de Borovoi et van Hamel [BovH06].

\medskip

Si la notion de r\'esolution flasque des groupes r\'eductifs connexes introduite dans [CT04, CT05] et d\'evelopp\'ee ici est nouvelle,
et s'il est \'etonnant qu'elle n'ait pas \'et\'e d\'egag\'ee plus
t\^ot, on peut dire qu'elle a \'et\'e pr\'epar\'ee par les travaux
ant\'erieurs de Voskresenski\u{\i} [Vo77, Vo98], Sansuc et l'auteur [CTSa77, CTSa87a, CTSa87b], Sansuc [Sa81],  
Borovoi [Bo96, Bo98], Gille [Gi01], Borovoi et Kunyavski\u{\i} [BoKu00, BoKu04].

\medskip
 
La th\'eorie d\'evelopp\'ee aux paragraphes 1 \`a 5 est originale. 
Les paragraphes 6 \`a 9 sont une pr\'esentation originale de concepts et
r\'esultats d\'ej\`a discut\'es par d'autres auteurs, au moins en caract\'eristique nulle ([Ko84, Ko86,
Bo96, Bo98, BoKu00, BoKu04]).
On s'est astreint \`a une pr\'esentation autonome de ces r\'esultats.
 Les exceptions sont : l'utilisation d'un th\'eor\`eme   de
Borovoi et Kunyavski\u{\i} [BoKu04], ind\'ependant du reste de leurs travaux (th\'eor\`eme {\bf 0.9} ci-apr\`es), 
le recours \`a certains arguments de [Bo98] dans les d\'emonstrations du th\'eor\`eme 8.4 (ii), du th\'eor\`eme 9.1 (iii)
et du th\'eor\`eme 9.4, et l'utilisation d'un th\'eor\`eme de P. Gille (appendice de [BoKu04])
dans les d\'emonstrations  du th\'eor\`eme 8.4 (i) et du th\'eor\`eme 9.3. 
Par ailleurs, jusqu'au milieu du paragraphe 8, on s'est astreint \`a d\'evelopper la th\'eorie
sur des corps de caract\'eristique arbitraire.

\bigskip

{\bf \S 0. Notations et  rappels}

\medskip

Soient $k$ un corps, $\k$ une cl\^oture {\it s\'eparable} de $k$ 
et $\g$ le groupe de Galois de $\k$ sur $k$. Un module galoisien
est un groupe ab\'elien muni de la topologie discr\`ete, \'equip\'e d'une
action continue de $\g$.
Etant donn\'es 
$X$ une $k$-vari\'et\'e alg\'ebrique  et $K/k$ une extension de corps,
on note $X_K=X\times_kK$. 
On note ${\overline X}=X\times_k\k$.
On note $K[X]=H^0(X_K,{\cal O}_{X_K})$. On note
 $\k[X]=H^0({\overline X},{\cal O}_{{\overline X}})$.
 Si $X$ et $Y$ sont deux $k$-vari\'et\'es, leur produit $X \times_kY$
 est g\'en\'eralement not\'e  $X \times Y$.
 Si  la $k$-vari\'et\'e  $X$ est int\`egre, on note $k(X)$ son corps des fonctions rationnelles.
Si $X$ est une $k$-vari\'et\'e g\'eom\'etriquement int\`egre, pour toute extension de corps $K/k$
on note $K(X)$ le corps des fonctions rationnelles de $X_K$.
On note $\G_{m,k}$, ou parfois simplement $\G_m$, le groupe multiplicatif sur $k$.
On note $\G_{a,k}$, ou parfois simplement $\G_a$, le groupe additif sur $k$.
On note $A^{\times}$ le groupe des unit\'es d'un anneau commutatif unitaire $A$.
Le groupe de Picard $\Pic(X)$ d'un sch\'ema $X$ est le groupe $H^1_{Zar}(X,{\cal O}_{X}^{\times})$.

\medskip

Sur un corps $k$ de caract\'eristique nulle, tout $k$-groupe lin\'eaire est lisse.
Sur un corps $k$ de caract\'eristique quelconque, un $k$-groupe r\'eductif est par
d\'efinition lisse.

Etant donn\'e un $k$-groupe alg\'ebrique lin\'eaire  $G$, on note ici
$G^*=\Hom_{\k-gp}({\overline G},\G_{m,\k})$ son groupe des caract\`eres.
C'est un module galoisien discret de type fini. Si $G$ est lisse et connexe,
le groupe ab\'elien sous-jacent est sans torsion.

Un $k$-groupe de type multiplicatif $M$ est un $k$-groupe qui sur $\k$  admet un plongement 
dans un produit de groupes multiplicatifs $\G_{m,\k} $. Un tel groupe est d\'etermin\'e par son groupe des caract\`eres $M^*=\Hom_{\k-gp}(\G_{m,\k}, {\overline M})$, qui est un groupe ab\'elien de type fini.
Le $k$-groupe $M$ est dit d\'eploy\'e si $\g$ agit trivialement sur $M^*$,
ce qui revient \`a dire que le groupe $M$ se plonge, sur $k$, dans un produit de $k$-tores $\G_{m,k}$.

Un $k$-tore $T$ est un $k$-groupe qui sur $\k$  est isomorphe \`a
un produit de groupes multiplicatifs $\G_{m,\k} $.   
On note $T_*=\Hom_{\k-gp}(\G_{m,\k}, {\overline T})$ son groupe des cocaract\`eres. On a l'isomorphisme de modules galoisiens $T_*=\Hom_{\Z}(T^*,\Z)$.

Pour $K/k$ une extension  finie s\'eparable de corps , on note
$R_{K/k}\G_{m}$ le $k$-tore descendu \`a la Weil du $K$-tore $\G_{m,K}$.

\medskip

On fera un usage constant et souvent tacite des propri\'et\'es suivantes.

\medskip

{\bf 0.1} 
Pour $X$ une $k$-vari\'et\'e lisse g\'eom\'etriquement int\`egre,
le module galoisien  $\k[X]^{\times}/\k^{\times}$ est un groupe ab\'elien libre
de type fini, et il est additif par rapport au produit de telles
vari\'et\'es
(Rosenlicht, cf. [CTSa77], Lemme 10 p.~188).

\medskip

{\bf 0.2} {\it (Lemme de Rosenlicht)}
Si $G$ est un $k$-groupe
lin\'eaire lisse 
connexe, le module galoisien $\k[G]^{\times}/\k^{\times}$ s'identifie au module
galoisien $G^*$ des caract\`eres de ${\overline G}$. 
Le groupe $k[G]^{\times}/k^{\times}$ s'identifie au groupe des caract\`eres de $G$ d\'efinis sur $k$.

\medskip

{\bf 0.3}  Soit $G$ un $k$-groupe r\'eductif connexe. On note $G^{ss} \subset G$ son groupe d\'eriv\'e.
C'est un $k$-groupe semi-simple distingu\'e dans $G$. On note $G^{tor}$ le $k$-groupe 
quotient de $G$ par $G^{ss}$. C'est un $k$-tore, c'est le plus grand quotient torique de $G$.
On a   la suite exacte 
exacte canonique de $k$-groupes r\'eductifs
$$1 \to G^{ss} \to G \to G^{tor} \to 1.$$
Le groupe des caract\`eres $(G^{tor})^*$ de $G^{tor}$ s'identifie au groupe des caract\`eres $G^*$ de $G$.
Si $G$ s'ins\`ere dans une suite exacte de $k$-groupes alg\'ebriques
$$ 1 \to G_{1} \to G \to T \to 1,$$
avec $G_{1}$ semi-simple et $T$ un $k$-tore, alors le sous-groupe $G_{1}  \subset  G$
co\"{\i}ncide avec $G^{ss} \subset G$, et ceci induit un isomorphisme $ G^{tor} \buildrel \simeq \over \rightarrow T$.

On note $G^{sc}$ le rev\^etement simplement connexe de $G^{ss}$. On a une suite exacte
canonique
$$ 1 \to \mu \to G^{sc} \to G^{ss} \to 1,$$
o\`u $\mu$ est un $k$-groupe de type multiplicatif fini, central dans $G^{sc}$, qu'on appelle
 le groupe fondamental de $G^{ss}$.

 \bigskip

{\bf 0.4} Supposons $k$ de caract\'eristique nulle.
Soit $G$ un $k$-groupe lin\'eaire  connexe. On note $G^{u} \subset G$ le radical unipotent de $G$.
C'est un $k$-groupe unipotent distingu\'e dans $G$, c'est le plus grand tel sous-groupe.
Le groupe $G^{u}$ admet une suite de composition dont les quotients successifs sont $k$-isomorphes au groupe additif $\G_{a,k}$. Sur toute $k$-vari\'et\'e affine, tout $G^{u}$-torseur (pour la topologie \'etale) est trivial. La $k$-vari\'et\'e sous-jacente au groupe $G^{u}$ est $k$-isomorphe \`a un espace affine sur $k$. Ceci implique  $\k^{\times} =\k[G^{u}]^{\times}$ et $\pic (G^{u} \times_{k}\k)=0$.

On a la suite exacte canonique
$$ 1 \to G^{u} \to G \to G^{\red} \to 1,$$
o\`u $G^{\red}$ est le plus grand $k$-groupe r\'eductif quotient  de $G$.
On note $G^{ss} \subset G^{\red}$ le groupe d\'eriv\'e de $G^{\red}$, et l'on note $G^{ssu}$ l'image r\'eciproque de $G^{ss}$ dans $G$. Le groupe $G^{ssu}$ est le noyau de l'homomorphisme
naturel de $G$ vers son plus grand quotient torique $G^{tor}$. 
 Le groupe des caract\`eres de $G^{tor}$ est le groupe des caract\`eres de $G$.
Le groupe 
$G^{ssu}$ est une extension de $G^{ss}$ par $G^{u}$. C'est un groupe connexe, simplement connexe.
On a ainsi les suites exactes canoniques de $k$-groupes connexes
$$1 \to G^{ssu} \to G \to G^{tor} \to 1$$
et
$$1 \to G^{u} \to G^{ssu} \to G^{ss} \to 1.$$
Le groupe $G^{ss}$ est simplement connexe (apr\`es passage \`a $\k$)
si et seulement si $G^{ssu}$ l'est.

\medskip

{\bf 0.5} Soit $1 \to G'' \to G \to G' \to 1$ une suite exacte de $k$-groupes alg\'ebriques lin\'eaires connexes, suppos\'es r\'eductifs si ${\rm car}(k)>0$. Cette suite induit une  suite exacte naturelle de modules galoisiens de type fini sur $\Z$
$$0 \to {G'}^* \to {G}^* \to {G''}^* \to \pic({{\overline G}'}) \to \pic({\overline G}) \to \pic({{\overline G}''}) \to 0. $$
L'exactitude jusqu'au terme $\pic({{\overline G}})$ fait l'objet du corollaire 6.11 de [Sa81], p. 43.
La surjectivit\'e de  $\pic({\overline G}) \to \pic({{\overline G}''})$
 est \'etablie dans la remarque 6.11.3  de [Sa81], p. 43.
 On a aussi la suite exacte de groupes ab\'eliens
 $$0 \to  ({G'}^*)^{\g} \to ({G}^*)^{\g}   \to ({G''}^*)^{\g} \to \Pic(G') \to \Pic(G) \to \Pic(G''),$$
 o\`u le groupe $({G}^*)^{\g}$ est le groupe des caract\`eres de $G$ d\'efinis sur $k$.
 Les groupes de Picard $\Pic(G)$ et $\Pic({\overline G})$ sont des groupes finis
 ([Sa81], Lemme 6.9 p. 41).

\medskip

{\bf 0.6} Soit $G$ un $k$-groupe semi-simple connexe. Le module galoisien $\pic({\overline G})$ 
 est isomorphe au groupe des caract\`eres  ${\mu}^*$ du groupe fondamental $\mu$ de ${\overline G}$
 ([Sa81], Lemme 6.9 (iii)). Il est nul si et seulement si $G$ est simplement connexe.

\medskip

{\bf 0.7}
Soit
$G$ un $k$-groupe alg\'ebrique extension d'un $k$-groupe de type
multiplicatif $K$ par un $k$-groupe de type multiplicatif $H$.
Si $G$ est commutatif ou $K$ connexe, alors $G$ est de type multiplicatif.
En particulier si $K$ et $H$ sont des $k$-tores, alors $G$
est un $k$-tore. Voir  [SGA3 II],
Expos\'e IX, Prop. 8.2 et Expos\'e XVII, Prop. 7.1.1

\medskip

{\bf 0.8}   Soit $M$ un $\g$-module $\Z$-libre de type fini.
Le module galoisien
$M$ est dit de permutation s'il poss\`ede une base
sur ${\Z}$  respect\'ee par
$\frak{g}$. Le module galoisien
$M$ est dit flasque, respectivement  coflasque,  
si pour tout sous-groupe
ouvert $\frak{h} \subset \frak{g}$, le groupe de cohomologie $H^1(\frak{h}, {\rm Hom}_{\Z}(M,\Z))$
respectivement  $H^1(\frak{h}, M)$,  est  nul.
Tout module de permutation est flasque et coflasque.

Un $k$-tore est dit quasi-trivial, resp. flasque, resp. coflasque, si son
groupe de caract\`eres est de permutation, resp. flasque, resp. coflasque.
On renvoie \`a [CTSa77], [CTSa87b] et [Vo98] pour les propri\'et\'es de base des
modules et tores flasques.

\medskip

 {\bf Th\'eor\`eme 0.9} (Borovoi-Kunyavski\u{\i} [BoKu04], Thm. 3.2)  {\it Soient $k$ un corps et
$G$ un $k$-groupe lin\'eaire connexe, suppos\'e r\'eductif si ${\rm car}(k)>0$. Soit $X$ une $k$-compactification lisse de $G$.
Le r\'eseau $\pic({\overline X})$ est  un module galoisien flasque.}

{\it Esquisse de d\'emonstration} La d\'emonstration proc\`ede par r\'eduction au cas des groupes r\'eductifs
quasi-d\'eploy\'es, par passage au corps des fonctions de la  vari\'et\'e des sous-groupes de Borel de $G$. Le cas o\`u $G$ est quasi-d\'eploy\'e, i.e.  poss\`ede un sous-groupe de Borel d\'efini sur $k$, \'etait connu ([CTSa77]) : il se ram\`ene au cas des tores ([Vo77], [CTSa77]). Dans [BoKu04], les auteurs
supposent le corps de base de caract\'eristique nulle. Leur d\'emonstration
s'\'etend aux corps de caract\'eristique quelconque, car les faits suivants valent sur un corps
$k$ quelconque.

(i) Tout $k$-tore alg\'ebrique admet une  $k$-compactification lisse (voir [CTHaSk05]).

(ii) Pour tout $k$-groupe r\'eductif connexe quasi-d\'eploy\'e $G$, on dispose de la d\'ecomposition de Bruhat, qui assure que $G$ contient un ouvert isomorphe au produit d'un $k$-tore et d'un espace affine
([SGA3 III],  Expos\'e XXII, \S 5.9).

(iii) Pour tout $k$-groupe r\'eductif connexe, la  $k$-vari\'et\'e des sous-groupes de Borel est une
$k$-vari\'et\'e (projective) lisse g\'eom\'etriquement connexe ([SGA3 III],  Expos\'e XXII, \S 5.9).

(iv) Si deux $k$-vari\'et\'es projectives, lisses, g\'eom\'etriquement int\`egres $X$ et $Y$ sont $k$-birationnellement \'equivalentes, alors il existe des $\g$-modules de permutation (de type fini)
$P_{1}^*$ et $P_{2}^*$ et un isomorphisme de $\g$-modules galoisiens discrets
${\rm Pic}({\overline X}) \oplus P_{1}^* \simeq {\rm Pic}({\overline Y})  \oplus P_{2}^*$ ([CTSa87b], Prop. 2.A.1 p. 461). \qed

\medskip

Sauf mention du contraire, la cohomologie utilis\'ee est la cohomologie \'etale,
qui lorsque l'on se restreint au cas du spectre d'un corps est la cohomologie galoisienne de ce corps.
Pour ${\cal F}$ un faisceau pour la topologie \'etale sur un sch\'ema $X$,  
 les groupes de cohomologie \'etale $H^n_{\et}(X,F)$ seront, sauf risque
d'ambigu\"{\i}t\'e, not\'es $H^n(X,{\cal F})$.

Pour tout sch\'ema $X$, on a $\Pic(X) = H^1_{Zar}(X,\G_{m}) \buildrel \simeq \over \rightarrow  H^1(X,\G_{m})$. Le groupe de Brauer de $X$ est le groupe $\Br(X)=H^2(X,\G_{m})$.

\bigskip

{\bf \S 1. Vari\'et\'es $\G_m$-quasi-triviales, vari\'et\'es $\G_m$-coflasques 
et vari\'et\'es finies-fac\-torielles}

\bigskip

{\bf D\'efinition 1.1} {\sl  Soit $X$ une $k$-vari\'et\'e
g\'eom\'etriquement int\`egre. On dit que 
$X$ est une $k$-vari\'et\'e $\G_m$-quasi-triviale
si les deux propri\'et\'es suivantes sont
satisfaites :

(i) Le module galoisien $\k[X]^{\times}/\k^{\times}$ est un module
de permutation.

(ii) On a ${\rm Pic}({\overline X})=0$.}

\bigskip

{\bf Proposition 1.2} {\sl  Soit $X$ une $k$-vari\'et\'e g\'eom\'etriquement int\`egre
$\G_m$-quasi-triviale.
Si $X$ est lisse, tout ouvert de $X$ est une $k$-vari\'et\'e
$\G_m$-quasi-triviale.} 

\medskip

{\it D\'emonstration} 
Pour $U$ ouvert non vide d'une $k$-vari\'et\'e lisse g\'eom\'etriquement
connexe $X$, de compl\'ementaire $F$,  on a la suite exacte de modules galoisiens
$$0 \to  \k[X]^{\times}/\k^{\times} \to \k[U]^{\times}/\k^{\times} \to
\Div_{{\overline F}}({\overline X})
\to \Pic({\overline X}) \to \Pic({\overline U}) \to 0.$$
Ici $\Div_{{\overline F}}({\overline X})$ est le groupe des diviseurs \`a support
dans ${\overline F}$. L'hypoth\`ese de lissit\'e assure que les diviseurs de Cartier co\"{\i}ncident avec les diviseurs de Weil,  le module galoisien  $\Div_{{\overline F}}({\overline X})$ est donc le module de permutation
de base les points de codimension 1 de ${\overline X}$ non dans ${\overline U}$.
Sous les hypoth\`eses de la proposition, on voit que $\Pic({\overline U})=0$
et que le module galoisien $\k[U]^{\times}/\k^{\times}$ est une extension
d'un module de permutation par un module de permutation.
Toute telle extension est $\g$-scind\'ee (lemme de Shapiro). 
Ainsi $\k[U]^{\times}/\k^{\times}$ est un module de permutation.
\cqfd

\bigskip

{\bf Proposition 1.3} {\it Soit
$X$ une $k$-vari\'et\'e g\'eom\'etriquement int\`egre $\G_m$-quasi-triviale.
Soit $T$ un $k$-tore. Si $T$ est flasque, ou si $\k^{\times} =  \k[X]^{\times} $,
la fl\`eche naturelle 
$H^1(k,T) \to H^1(X,T)$ induite par le morphisme structural
est bijective.
}

\medskip
{\it D\'emonstration}  
La suite spectrale de Leray pour le morphisme structural
$ p : X \to {\rm Spec}(k)$ et le faisceau \'etale d\'efini
par le $k$-tore $T$ donne la suite exacte courte
$$ 0 \to H^1(k,H^0({\overline X},{\overline T})) \to H^1(X,T) \to H^1({\overline X},{\overline T}).$$
Comme $T$ est un $k$-tore, l'hypoth\`ese ${\rm Pic}({\overline X})=0$
implique $H^1({\overline X},{\overline T})=0$. Par ailleurs, on a l'isomorphisme
de modules galoisiens $$H^0({\overline X},{\overline T}) \buildrel \simeq \over \rightarrow {\rm Hom}_{\bf Z}(T^*,
H^0({\overline X},{\bf G}_{m,\k}))={\rm Hom}_{\bf Z}(T^*,\k[X]^{\times} ).$$

Consid\'erons la suite exacte
$$ 1 \to \k^{\times} \to \k[X]^{\times} \to   P^* \to 0.$$
Comme $X$ est quasi-triviale, le module galoisien $P^*$ est un module de permutation.
D'apr\`es le lemme de Shapiro et le th\'eor\`eme 90 de Hilbert, on a $\Ext^1_{\g}(P^*,\k^{\times})=0$.
La suite exacte ci-dessus est donc scind\'ee comme suite de modules galoisiens.
Ceci donne naissance \`a la suite exacte scind\'ee de modules galoisiens
$$ 0 \to {\rm Hom}_{\bf Z}(T^*,\k^{\times}) \to
{\rm Hom}_{\bf Z}(T^*,\k[X]^{\times}) \to {\rm Hom}_{\bf Z}(T^*,P^*) \to
0.$$ 
Si $T$ est un $k$-tore  flasque, on a $H^1(k,{\rm Hom}_{\bf Z}(T^*,P^*) )=0.$
Dans le second cas consid\'er\'e dans la proposition, on a $P^*=0$. Ainsi dans ces deux cas  
l'application naturelle
$H^1(k,T) \to H^1(k,H^0({\overline X},{\overline T}))$ est bijective.
\cqfd

\bigskip

{\bf Corollaire 1.4} {\it  Soit
$X$ une $k$-vari\'et\'e g\'eom\'etriquement int\`egre $\G_m$-quasi-triviale,
et soit $T$ un $k$-tore. Si $T$ est flasque, ou si $\k^{\times} =  \k[X]^{\times} $,
tout torseur $Y \to X$ sous
$T$  est, en tant que   
$T$-torseur, et en particulier en tant que $k$-vari\'et\'e,
$k$-isomorphe \`a $T \times_kX$.}\cqfd

\bigskip

{\bf D\'efinition 1.5} {\it Soit $X$ une $k$-vari\'et\'e g\'eom\'etriquement int\`egre.
On dit que $X$ est une vari\'et\'e $\G_{m}$-coflasque si les deux propri\'et\'es suivantes sont
satisfaites :

(i) Le module galoisien $\k[X]^{\times}/\k^{\times}$ est un module
coflasque.

(ii) On a ${\rm Pic}({\overline X})=0$.}

\bigskip

{\bf Proposition 1.6} {\it Soit   $X$ une $k$-vari\'et\'e g\'eom\'etriquement int\`egre
$\G_{m}$-coflasque.
Si $X$ est lisse, tout ouvert de $X$ est une $k$-vari\'et\'e $\G_{m}$-coflasque.}
\medskip

{\it D\'emonstration} 
Soit $U \subset X$ un ouvert  non vide. Comme \`a la proposition 1.2, on  a la suite exacte
de modules galoisiens :
$$0 \to  \k[X]^{\times}/\k^{\times} \to \k[U]^{\times}/\k^{\times} \to
\Div_{{\overline F}}({\overline X})
\to \Pic({\overline X}) \to \Pic({\overline U}) \to 0.$$
Sous les hypoth\`eses de la proposition, ceci donne l'\'egalit\'e $\Pic({\overline U})=0$
et une suite exacte courte 
$$0 \to  \k[X]^{\times}/\k^{\times} \to \k[U]^{\times}/\k^{\times} \to
\Div_{{\overline F}}({\overline X}) \to 0.$$
Le module $\Div_{{\overline F}}({\overline X}) $ est un module de permutation.
Toute extension d'un module de permutation par un module coflasque est scind\'ee
(lemme de Shapiro et d\'efinition des modules coflasques).
Ainsi $\k[U]^{\times}/\k^{\times} \simeq  \k[X]^{\times}/\k^{\times} \oplus \Div_{{\overline F}}({\overline X})$ est un module coflasque. \cqfd

\bigskip

{\bf Proposition 1.7} {\it Soit   $X$ une $k$-vari\'et\'e g\'eom\'etriquement int\`egre
$\G_{m}$-coflasque. 

(i) Pour tout $k$-tore quasi-trivial   $P$, on a $H^1(X,P)=0$. 

(ii) Tout  torseur $Y \to X$
sous un $k$-tore quasi-trivial $P$ est, en tant que   
$P$-torseur, et en particulier en tant que $k$-vari\'et\'e,
$k$-isomorphe \`a $P \times_kX$.

(iii) Pour toute extension finie s\'eparable de corps $K/k$, on a $\pic(X_{K})=0$.}

\medskip

{\it D\'emonstration}
La suite spectrale de Leray pour le morphisme structural
$ p : X \to {\rm Spec}(k)$ et le faisceau \'etale d\'efini
par le $k$-tore $P$ donne la suite exacte courte
$$ 0 \to H^1(k,H^0({\overline X},{\overline P})) \to H^1(X,P) \to H^1({\overline X},{\overline P}).$$
Comme $P$ est un $k$-tore, l'hypoth\`ese ${\rm Pic}({\overline X})=0$
implique $H^1({\overline X},{\overline P})=0$. Par ailleurs, on a l'isomorphisme
de modules galoisiens $$H^0({\overline X},{\overline P}) \buildrel \simeq \over \rightarrow {\rm Hom}_{\bf Z}(P^*,
H^0({\overline X},{\bf G}_{m,\k}))={\rm Hom}_{\bf Z}(P^*,\k[X]^{\times} ).$$

La suite exacte 
$$ 1 \to \k^{\times} \to \k[X]^{\times} \to   Q^* \to 0,$$
o\`u $Q^*$ est un module coflasque, est $\Z$-scind\'ee. Elle
donne donc naissance \`a la suite exacte de modules galoisiens
$$ 0 \to {\rm Hom}_{\bf Z}(P^*,\k^{\times}) \to
{\rm Hom}_{\bf Z}(P^*,\k[X]^{\times}) \to {\rm Hom}_{\bf Z}(P^*,Q^*) \to
0.$$ Comme $Q^*$ est  un module coflasque, on a $H^1(k,{\rm Hom}_{\bf Z}(P^*,Q^*) )=0.$
Par ailleurs le th\'eor\`eme 90 de Hilbert assure $H^1(k,{\rm Hom}_{\bf Z}(P^*,\k^{\times})=0$.
Ainsi 
$H^1(k,H^0({\overline X},{\overline P}))=0$, ce qui \'etablit (i). Les \'enonc\'es (ii) et (iii)
sont des reformulations de (i).
\cqfd

\bigskip

{\bf D\'efinition 1.8} {\it Soit $X$ une $k$-vari\'et\'e g\'eom\'etriquement int\`egre.
On dit que  $X$ est une vari\'et\'e finie-factorielle si pour toute extension finie
s\'eparable de corps $K/k$ on a $\pic(X_{K})=0$.}

\bigskip

{\bf Proposition 1.9} {\it Soit   $X$ une $k$-vari\'et\'e g\'eom\'etriquement int\`egre
finie-factorielle.
Si $X$ est lisse, tout ouvert de $X$ est une $k$-vari\'et\'e finie-factorielle.}

\medskip

{\it D\'emonstration} Pour tout ouvert non vide $U$ de $X$ et toute extension finie
de corps $K/k$, la fl\`eche de restriction $\pic(X_{K}) \to \pic(U_{K})$ 
est surjective, puisque la $K$-vari\'et\'e $X_{K}$ est lisse. \cqfd
\medskip

{\it Remarque 1.9.1} Pour $k$ de caract\'eristique nulle et $X/k$ g\'eom\'etriquement int\`egre et lisse,
on peut montrer que si $X$ est finie-factorielle, alors elle est  ``universellement factorielle'' :
pour tout corps $F$ contenant $k$, on a $\pic(X_{F})=0$.

\bigskip

{\bf Proposition 1.10} {\it Soit   $X$ une $k$-vari\'et\'e g\'eom\'etriquement int\`egre
finie-factorielle.  Pour tout $k$-tore quasi-trivial   $P$, on a $H^1(X,P)=0$. Tout  torseur $Y \to X$
sous $P$ est, en tant que   
$P$-torseur, et en particulier en tant que $k$-vari\'et\'e,
$k$-isomorphe \`a $P \times_kX$.}

\medskip

{\it D\'emonstration} Le $k$-tore $P$ est un produit de $k$-tores $R_{K/k}\G_{m}$,
avec $K/k$ extension finie s\'eparable, et
$H^1(X,R_{K/k}\G_{m}) \buildrel \simeq \over \leftarrow \pic(X_{K})=0$. \cqfd

\bigskip

{\bf Proposition 1.11} {\it Soit   $X$ une $k$-vari\'et\'e g\'eom\'etriquement int\`egre.

(i) Si $X$ est $\G_{m}$-coflasque, elle est finie-factorielle.

(ii)  Si $X$ est finie-factorielle et si pour toute
extension finie s\'eparable $K/k$ de corps la fl\`eche 
$\br(K) \to \br(K(X))$ est injective, alors $X$ est $\G_{m}$-coflasque.

(iii)  Si $X$ est finie-factorielle et poss\`ede un $k$-point,
alors  $X$ est $\G_{m}$-coflasque.}

\medskip

{\it D\'emonstration} L'\'enonc\'e (i) n'est qu'une reformulation de la proposition 1.7.
Supposons $X$ finie-factorielle. Par passage \`a la limite sur les extensions finies s\'eparables de $k$,
on voit que l'on a $ \pic({\overline X})=0$.
Soit $K\subset \k$ une extension finie s\'eparable de $k$, soit $\frak{h}={\rm Gal}(\k/K)$.
 La suite spectrale de Leray pour le morphisme structural
$ p : X_{K} \to {\rm Spec}(K)$ et le faisceau \'etale d\'efini
par le $K$-tore $\G_{m,K}$ donne une injection
$H^1(\frak{h},\k[X]^{\times}) \hookrightarrow \pic(X_{K}) $. On  a donc
$H^1(\frak{h},\k[X]^{\times}) =0$. La suite exacte naturelle
$$1 \to  \k^{\times} \to \k[X]^{\times} \to     \k[X]^{\times}/\k^{\times} \to 1$$
donne naissance \`a la suite exacte de groupes de cohomologie
$$ H^1(\frak{h}, \k[X]^{\times}) \to  H^1(\frak{h}, \k[X]^{\times}/\k^{\times})
\to H^2(\frak{h},  \k^{\times}) \to H^2(\frak{h}, \k[X]^{\times}).$$
Le module galoisien $\k[X]^{\times}/\k^{\times}$ est un groupe ab\'elien de type fini sans torsion.
Compte tenu de $H^1(\frak{h},\k[X]^{\times}) =0$, on voit que 
l'on a $H^1(\frak{h}, \k[X]^{\times}/\k^{\times})=0$  pour tout $\frak{h}={\rm Gal}(\k/K)$
si  et seulement si la fl\`eche $H^2(\frak{h},  \k^{\times}) \to  H^2(\frak{h}, \k[X]^{\times})$ est injective pour tout
tel $\frak{h}$. C'est le cas si l'une des conditions mentionn\'ees en (ii) et (iii) est satisfaite. \cqfd

\medskip

{\it Remarque 1.11.1} Soit $k$ un corps et $C/k$ une conique lisse sans $k$-point.
Soit $P\in C$ un point ferm\'e de degr\'e 2, s\'eparable sur $k$, et soit $X$ le compl\'ementaire de
$P$ dans $C$. Alors  la $k$-vari\'et\'e  $X$ est finie-factorielle, mais elle n'est pas $\G_{m}$-coflasque.

\bigskip

{\bf \S 2. Groupes  quasi-triviaux et groupes coflasques}

\medskip

{\bf D\'efinition 2.1} {\it Soit $H$ un $k$-groupe lin\'eaire  connexe, suppos\'e r\'eductif
si ${\rm car}(k)>0$.
On dit que $H$ est 
un groupe   quasi-trivial si la $k$-vari\'et\'e alg\'ebrique $H$ est une vari\'et\'e 
$\G_m$-quasi-triviale, c'est-\`a-dire que les deux conditions suivantes sont satisfaites :

(1) Le groupe  $\k[H]^{\times}/{\overline k}^{\times}$ est 
un $\g$-module de permutation.

(2) Le groupe de Picard de ${\overline H}$ est nul.}

\bigskip

{\bf Proposition 2.2} {\it Soit $H$ un $k$-groupe lin\'eaire  connexe, suppos\'e r\'eductif
si ${\rm car}(k)>0$. Le $k$-groupe $H$ est quasi-trivial si et seulement si
les propri\'et\'es suivantes sont satisfaites :

(1.a) Le $k$-tore  $H^{tor}$ est un tore quasi-trivial.

(1.b) Le groupe des caract\`eres $H^{*}$ est un module de permutation.

(2.a) Le groupe connexe $H^{ssu}$ est simplement connexe.

(2.b) Le groupe semi-simple connexe $H^{ss}$, groupe d\'eriv\'e de $H^{\red}=H/H^u$, est simplement connexe.}

\medskip

{\it D\'emonstration} D'apr\`es {\bf 0.3}, les conditions (1.a) et (1.b) sont \'equivalentes,
et les conditions (2.a) et (2.b) sont \'equivalentes. D'apr\`es {\bf 0.2}, 
on a $\k[H]^{\times}/{\overline k}^{\times}=H^*$. 
Les conditions (1), (1.a) et (1.b) sont donc \'equivalentes.

En appliquant aux suites exactes de $k$-groupes
$$1 \to H^{ssu} \to H \to H^{tor} \to 1$$ 
et 
$$1 \to H^{u} \to H^{ssu} \to H^{ss} \to 1$$ 
les \'enonc\'es {\bf 0.4} et {\bf 0.5}, on obtient $\pic({\overline H}^{ss}) \buildrel \simeq \over \rightarrow \pic({\overline H}^{ssu})$ et $\pic({\overline H} ) \buildrel \simeq \over \rightarrow \pic({\overline H}^{ssu}) $. D'apr\`es {\bf 0.6},
le groupe semi-simple connexe $H^{ss}$ est simplement connexe si et seulement si $\pic({\overline H}^{ss})=0$. \cqfd

\medskip

On verra plus loin (Prop. 6.5) une autre caract\'erisation des groupes
 quasi-triviaux, qui justifie peut-\^etre
encore plus la terminologie : un groupe lin\'eaire  connexe $G$
est quasi-trivial si et seulement si son groupe fondamental
alg\'ebrique $\pi_1(G)$ est un module de permutation.

\bigskip

{\bf Proposition 2.3}
{\it   Soit
$$ 1 \to G_1 \to G_2 \to G_3 \to 1 $$
une suite exacte de $k$-groupes lin\'eaires connexes, suppos\'es
r\'eductifs si ${\rm car}(k)>0$.
Si $G_1$ et $G_3$ sont  quasi-triviaux, alors il en est
de m\^eme de $G_2$.}
\medskip

{\it D\'emonstration} D'apr\`es {\bf 0.5}, on a la suite exacte de modules galoisiens :
$$ 0 \to G_{3}^* \to  G_{2}^* \to  G_{1}^* \to
{\rm Pic}({{\overline G}_3}) \to {\rm Pic}({{\overline G}_2}) \to {\rm Pic}({{\overline G}_1}) \to 0.$$
Sous les hypoth\`eses de la proposition, on a ${\rm Pic}({{\overline G}_3})=0 $
et $ {\rm Pic}({{\overline G}_1})=0$, donc ${\rm Pic}({{\overline G}_2})=0$. Par ailleurs
le module galoisien $G_{2}^*$ 
est extension du module de permutation $G_{1}^*$ par le module de
permutation $G_{3}^* $, donc est un module de permutation.\cqfd

\bigskip

{\bf D\'efinition 2.4} {\it Soit $H$ un $k$-groupe lin\'eaire connexe, suppos\'e r\'eductif
si ${\rm car}(k)>0$.
On dit que $H$ est 
un groupe   coflasque si la $k$-vari\'et\'e alg\'ebrique $H$ est une vari\'et\'e 
$\G_m$-coflasque, c'est-\`a-dire que les deux conditions suivantes sont satisfaites :

(1) Le groupe  $\k[H]^{\times}/{\overline k}^{\times}$ est 
un $\g$-module coflasque.

(2) Le groupe de Picard de ${\overline H}$ est nul.}

\bigskip
\vfill\eject

{\bf Proposition 2.5}
{\it Soit $H$ un $k$-groupe lin\'eaire connexe, suppos\'e r\'eductif
si ${\rm car}(k)>0$. Le $k$-groupe $H$ est coflasque si et seulement si
les propri\'et\'es suivantes sont satisfaites :

(1.a) Le $k$-tore  $H^{tor}$ est un tore coflasque.

(1.b) Le groupe des caract\`eres $H^{*}$ est un module coflasque.

(2.a) Le groupe connexe $H^{ssu}$ est simplement connexe.

(2.b) Le groupe semi-simple connexe $H^{ss}$, groupe d\'eriv\'e de $H^{\red}=H/H^u$, est simplement connexe.}

{\it D\'emonstration} D'apr\`es {\bf 0.3}, les conditions (1.a) et (1.b) sont \'equivalentes. D'apr\`es 
{\bf 0.4},
les conditions (2.a) et (2.b) sont \'equivalentes. D'apr\`es {\bf 0.2}, 
on a $\k[H]^{\times}/{\overline k}^{\times}=H^*$. 
Les conditions (1), (1.a) et (1.b) sont donc \'equivalentes.

En appliquant aux suites exactes de $k$-groupes
$$1 \to H^{ssu} \to H \to H^{tor} \to 1$$ 
et 
$$1 \to H^{u} \to H^{ssu} \to H^{ss} \to 1$$ 
les \'enonc\'es {\bf 0.4} et {\bf 0.5}, on obtient $\pic({\overline H}^{ss}) \buildrel \simeq \over \rightarrow \pic({\overline H}^{ssu})$ et $\pic({\overline H} ) \buildrel \simeq \over \rightarrow \pic({\overline H}^{ssu}) $. D'apr\`es {\bf 0.6},
le groupe semi-simple connexe $H^{ss}$ est simplement connexe si et seulement si $\pic({\overline H}^{ss})=0$. \cqfd
\bigskip

{\bf Proposition 2.6}
{\it Soit $H$ un $k$-groupe lin\'eaire connexe, suppos\'e r\'eductif
si ${\rm car}(k)>0$. Le $k$-groupe $H$ est coflasque 
si et seulement si la $k$-vari\'et\'e $H$ est finie-factorielle, c'est-\`a-dire si pour toute extension finie s\'eparable $K/k$ on a  $\pic(H_{K})=0$.}

\medskip

{\it D\'emonstration} C'est une application directe de la proposition 1.11. \qed

\bigskip

{\bf \S 3.  R\'esolutions flasques d'un $k$-groupe lin\'eaire  connexe}

\medskip

La proposition suivante g\'en\'eralise 
un r\'esultat bien connu pour les 
$k$-tores alg\'ebriques ([CTS77], [Vo77], [Vo98]).

\medskip
{\bf Proposition-D\'efinition 3.1} 
{\it   Soit $G$ un $k$-groupe lin\'eaire connexe, suppos\'e r\'eductif si ${\rm car}(k)>0$.
Il existe un $k$-tore
flasque $S$, un $k$-groupe lin\'eaire connexe quasi-trivial $H$ 
 et une extension centrale de $k$-groupes
 $$ 1 \to S \to H \to G \to 1.$$ 
Une telle suite sera appel\'ee une {\bf r\'esolution flasque} du
$k$-groupe r\'eductif connexe $G$.
}

\medskip

{\it D\'emonstration}
Supposons d'abord  $G$  r\'eductif connexe. Soit $G^{ss}$ son groupe
d\'eriv\'e. Soit $Z$ la composante neutre du centre de $G$, qui est un $k$-tore
(c'est le radical de $G$). Ce $k$-tore est isog\`ene au $k$-tore $G^{tor}$, quotient de
$G$ par son groupe d\'eriv\'e.
Soit $G^{sc} \to G^{ss}$ 
le rev\^etement simplement connexe
de $G^{ss}$. 
Le morphisme  $G^{sc} \to G^{ss}$ est une isog\'enie centrale. L'image r\'eciproque du centre
de $G^{ss}$ est le centre de  $G^{sc}$.
 On dispose donc du $k$-homomorphisme surjectif 
  $G^{sc} \times Z \to G$.
Soit $(\alpha,\beta)$ dans le noyau $\mu$ de cet homomorphisme
(on sous-entend ici la donn\'ee d'une $k$-alg\`ebre commutative variable
dans laquelle on prend les points).
Son image $(\alpha',\beta) \in G^{ss}\times Z$ satisfait
$\alpha'.\beta=1 \in G$. Ainsi $\alpha'$ est dans le centre de $G^{ss}$,
et $\beta$ est d\'etermin\'e par $\alpha$.
Ceci implique que $\alpha$ est dans le centre $\nu$ de $G^{sc}$.
Ainsi $\mu \subset \nu \times Z$ est un $k$-groupe de type multiplicatif,
contenu dans le centre $ \nu \times Z$ de $G^{sc} \times Z$. En outre la projection sur le facteur
$G^{sc}$ induit un plongement $\mu \hookrightarrow \nu \subset G^{sc}$.
 Soit $Q \to Z$ un $k$-homomorphisme de $k$-tores,
avec $Q$ quasi-trivial. On obtient ainsi un $k$-homomorphisme surjectif
$ G^{sc}\times Q \to G$. Soit $M$ le noyau de cet homomorphisme. 
Soit $(g,q) \in M \subset  G^{sc}\times Q$. L'image de $(g,q)$
dans $G^{sc} \times Z$ par la fl\`eche \'evidente
est un couple $(g,p)$ qui est dans $\mu$. On a donc 
$g \in \nu$, et $(g,q) \in \nu \times Q$, qui est le centre
de $ G^{sc}\times Q$. Ainsi le $k$-groupe $M \subset \nu \times Z$ est 
un $k$-groupe de type multiplicatif central dans
 $ G^{sc}\times Q$.

On sait trouver
une suite exacte
$$ 1 \to M \to S \to P \to 1,$$
o\`u $S$ est un $k$-tore flasque et $P$ un $k$-tore quasi-trivial
([CTSa87a], (1.3.2)).
L'application dia\-gonale d\'efinit un plongement  du $k$-groupe de type multiplicatif $M$ 
dans le produit $(G^{sc} \times Q) \times S$. L'image est un sous-groupe central. Soit
$H$ le quotient du groupe r\'eductif $(G^{sc} \times Q) \times S$ par ce sous-groupe central. C'est un
$k$-groupe r\'eductif,  qui s'ins\`ere dans le diagramme commutatif de suites exactes de
$k$-groupes alg\'ebriques lin\'eaires suivant (o\`u l'on a remplac\'e la fl\`eche
donn\'ee $M \to S$ par son oppos\'ee) :

$$\diagram {&    & 1 && 1 && \cr
  &    & \downarrow{}{} && \downarrow{}{} && {}{}
&&    \cr 
1 &  \to & M    &\to & G^{sc} \times Q &
\to & G & \to & 1 \cr
 &    & \downarrow{}{} && \downarrow{}{} &&\downarrow{}{=} & 
 &       \cr
 1 & \to & S & \to & H & \to & G & \to &   1 \cr
&    & \downarrow{}{} && \downarrow{}{} && {}{} &&       \cr
  &  & P & = &P & &   & &  
\cr
&    & \downarrow{}{} && \downarrow{}{} &&  {}{}\cr
&    & 1 && 1. &&  
}
$$
Le quotient du groupe $H$ par le sous-groupe distingu\'e $G^{sc} \times 1 \subset G^{sc} \times  Q \subset H$
est un $k$-groupe  
extension du $k$-tore $P$ par le $k$-tore $Q$. D'apr\`es {\bf 0.7}, un tel groupe   est automatiquement un   $k$-tore. Comme toute extension d'un $k$-tore
quasi-trivial   par un $k$-tore quasi-trivial, dans la cat\'egorie des $k$-tores, est scind\'ee
(comme on voit au niveau des groupes
de caract\`eres, en utilisant le lemme de Shapiro), on conclut que le
quotient   de $H$ par $G^{sc}$ est   un $k$-tore quasi-trivial.
D'apr\`es~{\bf 0.3}, le groupe d\'eriv\'e $H^{ss}$ s'identifie \`a $G^{sc}$,
et le groupe $H^{tor}$ est un tore quasi-trivial : le groupe r\'eductif $H$ est un groupe quasi-trivial.

\medskip

Soient  maintenant $k$ un corps de caract\'eristique nulle et $G$ un $k$-groupe lin\'eaire connexe quelconque.
Soit  
$$1 \to S \to H_{1} \to G^{\red} \to 1$$
une r\'esolution flasque  
du $k$-groupe r\'eductif $G^{\red}$.
Soit $H$ le $k$-groupe produit fibr\'e de $H_{1}$ et de $G$ au-dessus de $G^{\red}$.
Le diagramme commutatif  de suites exactes de $k$-groupes lin\'eaires  
$$\diagram {&    &  && 1 && 1\cr
  &    & {}{} && \downarrow{}{} && \downarrow{}{}
&&    \cr 
 &   &    & & S  &
= &  S & &  \cr
 &    &  {}{} && \downarrow{}{} && \downarrow{}{} &&       \cr
 1 & \to & G^{u} & \to & H & \to &H_{1} &\to & 1 \cr
&    & \downarrow{=}{} && \downarrow{}{} && \downarrow{}{} &&       \cr
 1 &\to & G^{u}  &\to &G  &\to &  G^{\red} & \to &  1
\cr
&    &  {}{} && \downarrow{}{} && \downarrow{}{}\cr
&    &   && 1 & &1,  
}
$$
fournit  la r\'esolution flasque 
$$1 \to S \to H \to G \to 1$$
de $G$. \cqfd
\medskip

{\it Remarque 3.1.1} En suivant la d\'emonstration ci-dessus, on voit que l'on peut trouver 
une r\'esolution flasque $ 1 \to S \to H \to G \to 1$ telle que le $k$-tore $S$ et le $k$-tore
quasi-trivial $P=H^{tor}$ soient d\'eploy\'es par la plus petite extension finie galoisienne $K/k$
qui d\'eploie  \`a la fois le $k$-tore $G^{tor}$ (c'est \'equivalent \`a d\'eployer le $k$-tore $Z$
centre connexe de $G$) et le centre $\nu$ de $G^{sc}$.

\bigskip

{\bf Une construction alternative}

\medskip

Voici  une deuxi\`eme fa\c con
d'obtenir une r\'esolution flasque d'un $k$-groupe r\'eductif  connexe $G$.
Comme on l'a vu ci-dessus, il suffit de consid\'erer le cas
d'un groupe r\'eductif connexe.
On utilise l'existence d'une $z$-extension (cf. [Mi-Sh82], \S 3, Prop. 3.1  p.~297;
la litt\'erature sur ce sujet suppose souvent la caract\'eristique nulle.)
On dispose donc d'une suite exacte
$$ 1 \to P_1 \to G_1 \to G \to 1,$$
avec $P_1$ un $k$-tore quasi-trivial et 
$G_1$ un $k$-groupe r\'eductif dont le groupe d\'eriv\'e $G_1^{ss}$ est
(semi-simple) simplement connexe. Soit $T=G_{1}^{tor}$. Soit
$$1 \to S_2 \to P_2 \to T \to 1$$
une r\'esolution flasque du $k$-tore $T$ ([CTSa77],  [Vo77], [Vo98]),
c'est-\`a-dire une suite exacte de $k$-tores avec
 $S_2$  flasque
et
$P_2$  quasi-trivial.
Soit $H$ le groupe produit fibr\'e de $G_1$ et $P_2$ au-dessus de $T$.
On a la suite exacte
$$1 \to G_1^{ss} \to H \to P_2 \to 1.$$
Ainsi $H$ est un groupe r\'eductif quasi-trivial.
Par ailleurs le noyau de l'homomorphisme surjectif $H \to G$
compos\'e de $H \to G_1$ et de $G_1 \to G$
est un $k$-groupe alg\'ebrique extension du $k$-tore
$
P_1$ par le $k$-tore $S_2$. D'apr\`es {\bf 0.7}, c'est donc un $k$-tore $S$.
Comme il est  distingu\'e dans le $k$-groupe connexe $H$,
il est automatiquement central dans
$H$. Par ailleurs, ce $k$-tore $S$  est extension du
$k$-tore quasi-trivial $P_1$ par le $k$-tore flasque $S_2$.
Toute telle extension est automatiquement scind\'ee, comme on voit
sur les groupes de caract\`eres,
ainsi ce $k$-tore $S$ est $k$-isomorphe \`a $P_1 \times_k S_2$,
il est flasque.
\cqfd

\medskip

{\it Remarque 3.1.2} En caract\'eristique z\'ero, par une
construction toute diff\'erente, on donnera plus bas   (Th\'eor\`eme 5.4) une troisi\`eme  preuve de l'existence de r\'esolutions  flasques  pour les  groupes  lin\'eaires connexes.

\bigskip

{\bf Comparaison entre deux r\'esolutions flasques}

\medskip

{\bf Proposition 3.2 } { \it Soient $k$ un corps et $G$ un $k$-groupe lin\'eaire connexe, 
suppos\'e r\'eductif si ${\rm car}(k) > 0$.
Soient $1 \to S \to H \to G \to 1$
et $1 \to S_1 \to H_1 \to G \to 1$ deux r\'esolutions
flasques du $k$-groupe  $G$. 
Soient $P=H^{tor}$ et $P_{1}=H_{1}^{tor}$. 
Alors

(i) Il existe un isomorphisme de $k$-groupes
alg\'ebriques $S \times_k H_1 \simeq S_1 \times_k H$.

(ii) Les $k$-groupes   simplement connexes
$H^{ssu}$ et $H_{1}^{ssu}$ sont naturellement isomorphes, et les 
$k$-groupes semi-simples simplement connexes
$H^{ss}$  et $H_{1}^{ss}$ sont naturellement isomorphes.

(iii) 
Il existe  un isomorphisme de modules galoisiens
$S^* \oplus P_1^* \simeq S^*_1 \oplus P^*.$

(iv)  Il existe un isomorphisme
naturel entre ${\rm Coker} [S_{*} \to P_{*}]$ et
${\rm Coker} [S_{1*} \to P_{1*}]$.
}

\medskip

{\it D\'emonstration}
Soit $E$ le produit fibr\'e de $H$ et $H_1$ au-dessus
de $G$. On a donc le diagramme commutatif de suites exactes de
$k$-groupes 
$$\diagram {
&&1&& 1 &&1&&\cr
&&\uparrow  &&\uparrow && \uparrow && \cr
1 &\to &S& \to& H& \to &G& \to& 1 \cr
&&\uparrow {=}&&  \uparrow&& \uparrow&  \cr
1 &\to& S &\to& E &\to & H_{1}& \to &1 \cr
&&&&\uparrow &&\uparrow\cr
&&&&S_{1}&   \buildrel = \over \rightarrow&S_{1}&\cr
&&&&\uparrow &&\uparrow\cr
&&&&1&& 1 & .}$$

Comme le $k$-tore  $S_1$ est flasque et le groupe $H$ quasi-trivial,
d'apr\`es le corollaire 1.4  la projection
$E \to H$ admet une section $\sigma$. 
Quitte \`a multiplier cette section par un 
\'el\'ement de $S_{1}(k)$, on peut supposer
qu'elle envoie l'\'el\'ement neutre $e_H \in H(k)$ sur l'\'el\'ement neutre $e_E\in E(k)$. Cette section est
a priori un morphisme de vari\'et\'es alg\'ebriques.

Le morphisme $\theta: H \times H \to E$ 
d\'efini par $\theta(x,y)=\sigma(xy)\sigma(y)^{-1}\sigma(x)^{-1}$
envoie $H$ dans ${S_1}$, et il v\'erifie 
$\theta(x,e_H)=e_{S_1}$ et $\theta(e_H,y)=e_{S_1}$.
Un tel morphisme s'\'ecrit $\theta(x,y)=\alpha(x)\beta(y)$,
o\`u $\alpha$ et $\beta$ sont des $k$-morphismes de groupes
alg\'ebriques de $H$ dans $S_1$ (lemme de Rosenlicht, {\bf 0.2}).
 Par la formule ci-dessus, ces morphismes $\alpha$ et $\beta$ sont constants et envoient $H$ sur $e_{S_1}$. Ainsi $\sigma(xy)=\sigma(x)\sigma(y)$, et $\sigma$ est
un homomorphisme section de l'homomorphisme $E \to H$.
(Je remercie O. Gabber pour cet argument.)

Comme le noyau $S_1$ de cet homomorphisme est central,
on conclut que le $k$-groupe $E$ est $k$-isomorphe au
produit de groupes $S_1 \times_kH$. Le m\^eme argument vaut
pour $E \to H_1$, ce qui \'etablit l'\'enonc\'e (i) :
il existe un isomorphisme de $k$-groupes
alg\'ebriques $S \times_k H_1 \simeq S_1 \times_k H$.
Les projections $E \to H$ et $E \to H_{1}$ induisent des isomorphismes
$E^{ssu} \buildrel \simeq \over \rightarrow H^{ssu}$ et $E^{ssu} \buildrel \simeq \over \rightarrow H_{1}^{ssu}$.
L'\'enonc\'e (ii) en r\'esulte imm\'ediatement.

Toute extension d'un $k$-tore quasi-trivial par un $k$-tore flasque est
automatiquement scind\'ee. On voit alors que
le diagramme ci-dessus 
induit   un diagramme  de suites exactes scind\'ees
de $k$-tores
$$\diagram {
&& && 1 && &&\cr
&&   &&\uparrow &&   && \cr
   &  & &  & P&   && &   \cr
&&
&&  \uparrow&&  &  \cr
1 &\to& S &\to& E^{tor} &\to & P_{1}& \to &1 \cr
&&&&\uparrow && \cr
&&&&S_{1}&    & &\cr
&&&& \uparrow && \cr
&&&&1&&   & }$$
Les \'enonc\'es (iii) et (iv) en r\'esultent.
\cqfd

\medskip

{\it Remarque 3.2.1} Le diagramme de suites exactes scind\'ees ci-dessus induit
pour tout entier $n \geq 0$  un isomorphisme naturel
$$\Ker [H^n(k,S) \to H^n(k,P)] \simeq  \Ker [H^n(k,S_{1}) \to H^n(k,P_{1})].$$
Pour $n=1$, ceci donne en particulier un isomorphisme naturel
$$H^1(k,S) \simeq H^1(k,S_{1}).$$

\medskip

{\it Remarque 3.2.2}  Si l'on se donne trois r\'esolutions flasques $1 \to S_{i} \to H_{i} \to G \to 1$ ($i=1,2,3$) de $G$, on v\'erifie que le
compos\'e des isomorphismes naturels ${\rm Coker} [S_{1*} \to P_{1*}] \simeq
{\rm Coker} [S_{2*} \to P_{2*}]$ et ${\rm Coker} [S_{2*} \to P_{2*}] \simeq
{\rm Coker} [S_{3*} \to P_{3*}]$ construits au point (iv)  est l'isomorphisme naturel ${\rm Coker} [S_{1*} \to P_{1*}] \simeq
{\rm Coker} [S_{3*} \to P_{3*}]$.

\bigskip

{\bf  Groupe de Picard d'un groupe lin\'eaire et r\'esolution flasque}

\medskip

{\bf Proposition 3.3}
{\it Soit $G$ un $k$-groupe lin\'eaire connexe, suppos\'e r\'eductif si ${\rm car}(k)>0$.
Soit
$1 \to S \to H \to G \to 1$ une r\'esolution flasque de $G$.
Soit $P$ le $k$-tore quasi-trivial $H^{tor}$. Le noyau de l'application compos\'ee
  $S \to H \to H^{tor}=P$ est fini, et on  a la suite exacte naturelle
$$  (P^*)^{\g} \to (S^*)^{\g} \to \Pic(G) \to  0,$$
o\`u la fl\`eche $(S^*)^{\g} \to \Pic(G) $ associe \`a tout
caract\`ere  $\chi : S \to \G_{m,k}$ la classe du $\G_{m,k}$-torseur sur $G$
d\'eduit de la r\'esolution flasque en poussant le long de $\chi$.
}

\medskip

{\it D\'emonstration}
En utilisant {\bf 0.4} et {\bf 0.5}, on \'etablit $\Pic(H^{ssu})=0$, puis $\Pic(H)=0$. On voit aussi que l'on a $P^* \simeq H^*$. Appliquant  {\bf 0.5}  \`a la r\'esolution flasque, on obtient  la suite exacte annonc\'ee.  
Ceci vaut sur tout corps, en particulier pour $k={\overline k}$. On voit ainsi que
le conoyau de
$P^* \to S^*$  est le groupe $ \Pic({\overline G}) $, qui est fini (rappel {\bf 0.5}).
Par dualit\'e, la fl\`eche de $k$-tores $S \to P$ a un noyau fini.
\cqfd

\bigskip

{\bf Le diagramme fondamental attach\'e \`a une r\'esolution flasque d'un $k$-groupe lin\'eaire}

\medskip

Soit $G$ un  $k$-groupe  lin\'eaire  connexe, suppos\'e lisse si ${\rm car}(k)>0$.
 Soit $1  \to S \to H \to G  \to 1$
une r\'esolution flasque de $G$, 
avec $H$ extension du $k$-tore
quasi-trivial $P=H^{tor}$ par le $k$-groupe  
simplement connexe $ H^{ssu}$.
Le morphisme $H \to G$ induit un \'epimorphisme de $k$-tores
$H^{tor} \to G^{tor}$. Soit $M$ le $k$-groupe de type multiplicatif noyau de
$H^{tor} \to G^{tor}$. On a un homomorphisme  induit $S \to  M$.
Soit $\mu \subset S$ le
 $k$-groupe de type multiplicatif  
noyau de la fl\`eche compos\'ee $S \to H \to H^{tor}=P$,
qui est aussi le noyau de $S \to  M$.
D'apr\`es la  proposition 3.3, c'est un $k$-groupe de type multiplicatif  fini.
Le groupe $H^{ssu}$, resp. $G^{ssu}$ est le noyau de la fl\`eche
$H \to H^{tor}$, resp. de la fl\`eche $G \to G^{tor}$.
Le morphisme de $k$-groupes  $H \to G$ induit un
 morphisme de $k$-groupes $H^{ssu} \to G^{ssu}$.
On a donc le diagramme de complexes de $k$-groupes lin\'eaires suivant
 $$\diagram {&    & 1 && 1 && 1
\cr
  &    & \downarrow{}{} &&  \downarrow{}{} &&  \downarrow{}{}
&&    
\cr 
1 &\to &\mu &\to &H^{ssu} &\to & G^{ssu} &\to& 1
 \cr
 &    &  \downarrow{}{} &&  \downarrow{}{} &&  \downarrow{}{} &&     
   \cr
 1 & \to & S & \to & H & \to &G &\to  &  1 \cr
&    &  \downarrow{}{} &&  \downarrow{}{} &&  \downarrow{}{} &&   
    \cr
 1 &  \to & M    &\to & P &
\to & G^{tor} & \to & 1
\cr
&    &  \downarrow{}{} &&  \downarrow{}{} &&  \downarrow{}{}
\cr
&    & 1 && 1 && 1,
}
$$
o\`u 
 les deux  suites   horizontales inf\'erieures sont exactes.
Le $k$-groupe $G^{ssu}$ est une extension d'un
$k$-groupe semi-simple connexe par un $k$-groupe unipotent.
Il n'admet donc pas de $k$-morphisme non constant dans le
$k$-groupe de type multiplicatif quotient de $M$ par l'image de $S$. Ainsi
le diagramme ci-dessus est constitu\'e de suite exactes.
Le groupe
 $H^{ssu}$, qui est simplement connexe,
est le rev\^etement simplement connexe
de $G^{ssu}$, le $k$-groupe  de type multiplicatif fini $\mu$ s'identifie au groupe fondamental de $G^{ssu}$,
qui est  le noyau de $G^{sc} \to G^{ss}$.

Pour $G$ r\'eductif connexe, le
diagramme fondamental s'\'ecrit
$$\diagram {&    & 1 && 1 && 1
\cr
  &    & \downarrow{}{} &&  \downarrow{}{} &&  \downarrow{}{}
&&    
\cr 
1 &\to &\mu &\to &G^{sc } &\to & G^{ss } &\to& 1
 \cr
 &    &  \downarrow{}{} &&  \downarrow{}{} &&  \downarrow{}{} &&     
   \cr
 1 & \to & S & \to & H & \to &G &\to  &  1 \cr
&    &  \downarrow{}{} &&  \downarrow{}{} &&  \downarrow{}{} &&   
    \cr
 1 &  \to & M    &\to & P &
\to & G^{tor} & \to & 1
\cr
&    &  \downarrow{}{} &&  \downarrow{}{} &&  \downarrow{}{}
\cr
&    & 1 && 1 && 1.
}
$$

\bigskip

{\bf \S 4 R\'esolutions coflasques d'un $k$-groupe lin\'eaire connexe}

\medskip

L'\'enonc\'e ci-dessous g\'en\'eralise l'existence des r\'esolutions coflasques ``du deuxi\`eme type''
pour $G$ un $k$-tore alg\'ebrique ([CTSa87a],  (1.3.4)). 
Pour un $k$-groupe donn\'e $G$,
il distingue aussi une sous-classe particuli\`ere parmi toutes les $z$-extensions 
possibles pour ce groupe.
\medskip

{\bf Proposition 4.1}   
{\it   
Soit $G$ un $k$-groupe lin\'eaire connexe, suppos\'e r\'eductif si ${\rm car}(k)>0$.
Il existe un $k$-tore quasi-trivial  $P$,  un $k$-groupe lin\'eaire connexe (r\'eductif si ${\rm car}(k)>0$) extension d'un $k$-tore coflasque $Q$ par un groupe  (semi-simple si ${\rm car}(k)>0$) 
simplement connexe $H$, 
 et une extension centrale de $k$-groupes lin\'eaires connexes
 $$ 1 \to P \to H \to G \to 1.$$ }
 
\medskip

{\it D\'emonstration} Soit 
$$1 \to S \to H_{1} \to G \to 1$$
une r\'esolution flasque de $G$ (pour la d\'emonstration qui suit, on pourrait aussi partir d'une $z$-extension). On sait  (voir [CTSa87a],  Lemma 0.6 et Prop.  1.3)
qu'il existe une suite exacte de $k$-tores
$$ 1 \to S \to P \to Q \to 1,$$
o\`u $P$ est un $k$-tore quasi-trivial et $Q$ est un $k$-tore coflasque.

Consid\'erons alors le $k$-groupe alg\'ebrique lin\'eaire $H$ qui
est le quotient de $P \times_{k} H_{1}$ par l'image diagonale de $S$, qui est centrale et donc distingu\'ee dans 
$P \times_{k} H_{1}$.
Le $k$-groupe $H$ s'ins\`ere dans le diagramme commutatif de suites exactes de
$k$-groupes alg\'ebriques lin\'eaires suivant (o\`u l'on a remplac\'e la fl\`eche donn\'ee
$S \to P$ par son oppos\'ee) :

 $$\diagram {&    & 1 && 1 && \cr
  &    & \downarrow{}{} && \downarrow{}{} && {}{}
&&    \cr 
1 &  \to & S    &\to & H_{1} &
\to &G & \to & 1 \cr
 &    & \downarrow{}{} && \downarrow{}{} &&\downarrow{=}{} &&       \cr
 1 & \to & P& \to & H & \to &G &\to 
   &  1 \cr
&    & \downarrow{}{} && \downarrow{}{} && {}{} &&       \cr
  &  & Q & = &Q & &   & &  
\cr
&    & \downarrow{}{} && \downarrow{}{} &&  {}{}\cr
&    & 1 && 1 && .  
}
$$

On a $H^{ssu}=H_{1}^{ssu}$, donc $H_{1}^{ssu}$ est un groupe simplement connexe.
Le tore $H^{tor}$ est extension du $k$-tore coflasque $Q$ par le $k$-tore quasi-trivial $H_{1}^{tor}$.
Toute   extension d'un $k$-tore coflasque par un $k$-tore quasi-trivial est scind\'ee,
donc $H^{tor}$ est un $k$-tore coflasque.\cqfd

\bigskip

{\bf Proposition 4.2} { \it Soient $k$ un corps et $G$ un $k$-groupe lin\'eaire connexe, 
suppos\'e r\'eductif si ${\rm car}(k) > 0$.
Soient $1\to P \to H  \to G \to 1$ et $1 \to P_{1} \to H_{1} \to G \to 1$
deux r\'esolutions coflasques du $k$-groupe  $G$.
Soient $Q=H^{tor}$ et $Q_{1}=H_{1}^{tor}$ les $k$-tores coflasques associ\'es.

(i) Il existe un $k$-isomorphisme de $k$-groupes
alg\'ebriques $P \times_k H_1 \simeq P_1 \times_k H$.

(ii) Les $k$-groupes simplement connexes $H^{ssu}$ et $H_{1}^{ssu}$ sont naturellement isomorphes, et les 
$k$-groupes 
semi-simples simplement connexes
$H^{ss}$ et $H_{1}^{ss}$ sont naturellement isomorphes.

(iii) Il existe un $k$-isomorphisme de $k$-tores $P \times_{k} Q_{1} \simeq P_{1} \times_{k}Q$.

(iv)  Il existe un isomorphisme
naturel entre ${\rm Coker} [P_{*} \to Q_{*}]$ et
${\rm Coker} [P_{1*} \to Q_{1*}]$.
}

\medskip

{\it D\'emonstration}
Soit $E$ le produit fibr\'e de $H$ et $H_1$ au-dessus
de $G$. On a donc le diagramme commutatif de suites exactes de
$k$-groupes 
$$\diagram {
&&1&& 1 &&1&&\cr
&&\uparrow  &&\uparrow && \uparrow && \cr
1 &\to &P& \to& H& \to &G& \to& 1 \cr
&&\uparrow {=}&&  \uparrow&& \uparrow&  \cr
1 &\to& P &\to& E &\to & H_{1}& \to &1 \cr
&&&&\uparrow &&\uparrow\cr
&&&&P_{1}&  \buildrel \simeq \over \rightarrow &P_{1}&\cr
&&&&\uparrow &&\uparrow\cr
&&&&1&& 1 &.}$$

Comme le $k$-tore  $P_1$ est quasi-trivial et le $k$-groupe alg\'ebrique $H$ coflasque,
le $P_{1}$-torseur $E \to H$ est trivial (proposition 1.10 et proposition 2.6).
Le m\^eme argument que dans la d\'emonstration de la proposition 3.2
assure qu'alors le $k$-groupe  $E$ est $k$-isomorphe, comme $k$-groupe,
au $k$-groupe $P_{1} \times_{k} H$. On voit de m\^eme qu'il est 
$k$-isomorphe au $k$-groupe $P \times_{k} H_{1}$. Ceci \'etablit (i).
Les \'enonc\'es (ii) et (iii)  sont des cons\'equences imm\'ediates.
L'argument pour (iv) est identique \`a celui donn\'e dans la d\'emonstration de 
la proposition 3.2. \cqfd

\bigskip

{\it Remarque 4.2.1} Si $1\to P \to H  \to G \to 1$ et $1 \to P_{1} \to H_{1} \to G \to 1$
sont deux $z$-extensions du m\^eme $k$-groupe r\'eductif connexe $G$, on peut
montrer que les $k$-tores $H_{1}^{tor}$ et $H_{2}^{tor}$ sont stablement
$k$-birationnels, i.e. $k$-birationnels apr\`es multiplication de chacun d'eux
par un espace affine convenable.

\bigskip

{\bf \S 5. Propri\'et\'es birationnelles et comparaison avec les torseurs universels}

\medskip

La lecture de ce paragraphe suppose une familiarit\'e avec la notion de torseur universel ([CTSa87b]).

\medskip

{\bf Proposition 5.1} {\it Soit $X$ une $k$-vari\'et\'e projective,
lisse et g\'eom\'etriquement int\`egre. Supposons que
le groupe  ${\rm Pic}({\overline X})$ est libre de type fini.
Soit ${\cal T} \to X$ un torseur universel sur $X$.
Tout ouvert $U$ de ${\cal T}$
est une vari\'et\'e $\G_m$-quasi-triviale.}
\medskip
{\it D\'emonstration} 
 Dans [CTSa87b], \S 2.1, on a \'etabli ce r\'esultat
pour $U={\cal T}$. D'apr\`es le lemme 1.2, ceci implique le r\'esultat
pour tout ouvert.\cqfd

\bigskip

{\bf Proposition 5.2} 
{\it Soient $k$ un corps et $G$ un $k$-groupe lin\'eaire connexe, suppos\'e r\'eductif
si ${\rm car}(k) >0$.
Soit $1 \to S \to H \to G \to 1$
une r\'esolution flasque de $G$. 
Soit $X$ une $k$-compactification lisse de $G$.
Soit $S_0$ le $k$-tore de groupe des caract\`eres
le module galoisien ${\rm Pic}({\overline X})$.
Soit ${\cal T} \to X$ un torseur universel sur $X$ 
de fibre triviale en l'\'el\'ement neutre $e_{G}$ de $G(k)$. Alors :

(i) Il existe un $k$-isomorphisme de $k$-vari\'et\'es
$S_0 \times_k H \simeq S \times_k {\cal T}_G$, o\`u
${\cal T}_G$ d\'esigne la restriction de ${\cal T}$
au-dessus de $G$. 

(ii) Il existe des modules de permutation 
$P^*$ et $P_0^*$ et un isomorphisme de modules
galoisiens
$S^* \oplus P^* \simeq S^*_0 \oplus P_0^*$.}
\medskip
{\it D\'emonstration} On consid\`ere la $k$-vari\'et\'e
$Y$ produit fibr\'e de $H$ et de ${\cal T}_G$
au-dessus de $G$. Cette vari\'et\'e poss\`ede un $k$-point.
Elle est munie d'une part d'une structure de torseur
sur $H$ sous $S_0$, d'autre part d'une structure
de torseur sur ${\cal T}_G$ sous $S$.
D'apr\`es le  th\'eor\`eme {\bf 0.9} (Borovoi
et Kunyavski\u{\i}),
le $k$-tore $S_0$ est flasque.  
D'apr\`es le corollaire 1.4, on voit
qu'on a des isomorphismes de $k$-vari\'et\'es
$Y \simeq S_0 \times_k H$ et $Y \simeq S \times_k {\cal T}_G$.
Ceci \'etablit le point (i).

Pour \'etablir le point (ii), on applique le
foncteur  $Z \mapsto \k[Z]^{\times}/\k^{\times}$
\`a $Z= S_0 \times_k H$ et
\`a $Z= S \times_k {\cal T}_G$. 
Pour ce dernier, on utilise le fait (Prop. 5.1) que
l'ouvert $Z= {\cal T}_G$ de ${\cal T}$ est $\G_m$-quasi-trivial,
donc  le module
$\k[{\cal T}_{G}]^{\times}/\k^{\times}$ est un module de permutation. \cqfd

\bigskip

{\bf Corollaire 5.3} 
{\it Soient $k$ un corps et $G$ un $k$-groupe lin\'eaire connexe, suppos\'e r\'eductif
si ${\rm car}(k) >0$.
  Si  $G$  est
stablement $k$-rationnel, i.e. $k$-birationnel \`a un espace affine sur son corps
de base apr\`es produit par un espace affine, et si $G$ admet une $k$-compactification lisse,
alors il existe une suite exacte
$$ 1 \to P \to H \to G \to 1,$$
o\`u  $H$ est un $k$-groupe quasi-trivial et  $P$ est un $k$-tore quasi-trivial central dans
$H$.}

\medskip

{\it  D\'emonstration} Soit
$1 \to S \to H \to G \to 1$ une r\'esolution flasque de $G$.
Soit $X$ une
$k$-compactification lisse de $G$. Si $X$
est stablement $k$-rationnelle, alors $\Pic({\overline X})$ est stablement
de permutation ([CTSa87b], Prop. 2.A.1).
D'apr\`es la proposition 5.2, on voit
que $S^*$ est aussi stablement de permutation.
Il  existe donc des $k$-tores quasi-triviaux $P_1$ et $P_2$
tels que $S \times_k P_1 \simeq P_2$.
En rempla\c cant la r\'esolution flasque donn\'ee par
$$1 \to P_2 \to H \times_kP_1 \to G \to 1,$$
on obtient une suite du type annonc\'e. \cqfd

\medskip

\bigskip

{\bf Th\'eor\`eme 5.4}
{\it Soient $k$ un corps et $G$ un $k$-groupe lin\'eaire connexe, suppos\'e r\'eductif
si ${\rm car}(k) >0$. Soit $X$ une $k$-compactification lisse
de $G$. Soit  ${\cal T} \to X$ un torseur universel sur $X$ de fibre
triviale en l'\'el\'ement neutre $e_{G} \in G(k) \subset X(k)$.
 
 (i) L'ouvert  ${\cal T}_G={\cal T}\times_XG \subset {\cal T}$ peut \^etre muni
d'une structure de $k$-groupe alg\'ebrique connexe $H$
tel que la projection ${\cal T}_G \to G$ soit un homomorphisme
surjectif de noyau le $k$-tore $S_0$ de groupe des caract\`eres
le groupe ${\rm Pic}({\overline X})$. 

(ii) La suite exacte $1 \to S_{0} \to H \to G \to 1$ est une r\'esolution flasque de $G$.}

\medskip

{\it D\'emonstration}
 L'\'enonc\'e (i)  est un cas particulier du th\'eor\`eme 5.6 ci-dessous.
 La vari\'et\'e ${\cal T}_G$ est une vari\'et\'e $\G_m$-quasi-triviale 
(Prop. 5.1). Ainsi  $H$ est un groupe lin\'eaire connexe quasi-trivial.
D'apr\`es le th\'eor\`eme {\bf 0.9} (Borovoi-Kunyavski\u{\i}), 
le $k$-tore
 $S_0$ est   flasque. Ceci \'etablit le point (ii). \cqfd

 \medskip
 
 {\it Remarque 5.4.1}
 En  caract\'eristique z\'ero, d'apr\`es
 Hironaka, il existe une $k$-compactification lisse de la $k$-vari\'et\'e lisse $G$. 
  En caract\'eristique  z\'ero, on obtient   ainsi
 une d\'emonstration alternative de l'existence de r\'esolutions flasques (Proposition 3.1).

\medskip

{\it Remarque 5.4.2} Dans [Gi97], III.4.b, Gille \'etudie
les torseurs universels sur les compactifications lisses
de groupes semi-simples, et il \'etablit un lien
avec le rev\^etement simplement connexe d'un tel groupe.

\bigskip

Notre but maintenant est d'\'etablir le th\'eor\`eme 5.6 ci-dessous, 
ce qui se fera en adaptant un argument de Serre ([Se59], Chapitre VII, \S 15).
Commen\c cons par des pr\'eliminaires.
Soient $k$ un corps,
 $S$ un $k$-groupe de type multiplicatif lisse et $X$ une $k$-vari\'et\'e.
Le groupe de cohomologie \'etale $H^1_{\acute et}(X,S)$ classifie les torseurs sur $X$ sous $S$,
\`a isomorphisme non unique pr\`es.
 Tout endomorphisme d'un torseur sur $X$ sous $S$ est un automorphisme (de
$S$-torseur), donn\'e par un \'el\'ement de $\Hom_{k}(X,S)={\rm Mor}_{k}(X,S)$.
On a une variante, utile pour la suite, o\`u $X$ est ($k$-)point\'e
par un $k$-point $e_{X} \in X(k)$, et chaque torseur est ($k$-)point\'e au-dessus 
du point $e_{X}$ de
$X$. Tout endomorphisme d'un tel torseur point\'e est un automorphisme,
donn\'e par un \'el\'ement de 
 ${\rm Mor}_{k}(X,S)$  envoyant
le point marqu\'e de $X$ sur l'\'el\'ement neutre $e_{S}$ de $S$.
Pour tout entier $i$ on introduit le groupe $H^{i}_{e_{X}}(X,S)$ d\'efini comme
le noyau de l'\'evaluation au point $e_{X} \in X(k)$ :
$$H^{i}_{e_{X}}(X,S) = \Ker \hskip1mm   {\rm ev}_{e_{X}} : H^{i}(X,S) \to
H^{i}(k,S).$$

La suite spectrale $H^p(k,H^q({\overline X},S)) \Longrightarrow H^n(X,S)$
donne naissance \`a la suite exacte
$$ 0 \to H^1(k, H^0({\overline X},S))
\to H^1(X,S) \to H^0(k,H^1({\overline X},S)) \to H^2(k, H^0({\overline X},S)) \to H^2(X,S).$$

Par comparaison avec le point marqu\'e de $X$, on en d\'eduit la suite
exacte
$$ 0 \to H^1(k, H^0_{e_{X}}({\overline X},S))
\to H^1_{e_{X}}(X,S) \to H^0(k,H^1({\overline X},S)) \to H^2(k, H^0_{e_{X}}({\overline X},S)) \to H^2_{e_{X}}(X,S).$$

On dit qu'un foncteur contravariant $F$ d'une cat\'egorie ${\cal C}$ de 
$k$-vari\'et\'es, stable par produits, dans
une cat\'egorie de modules est additif si
pour $X,Y$ dans ${\cal C}$, la fl\`eche naturelle \break
$F(X) \oplus F(Y) \to F(X\times_kY)$ d\'efinie
comme la somme des deux images inverses des projections
est un isomorphisme.  

Supposons ${\overline X}$ connexe et normale. On a alors
$H^0({\overline X},S)/H^0(\k,S)=\Hom_{\Z}(S^*, \k[X]^{\times}/\k^{\times})$,
et le foncteur envoyant $X$ sur $H^0({\overline X},S)/H^0(\k,S)$
est additif sur la cat\'egorie des $k$-vari\'et\'es
$X$ telles que ${\overline X}$ soit normale connexe.
La suite exacte ci-dessus s'\'ecrit alors 
$$ 0 \to H^1(k, \Hom_{\Z}(S^*, \k[X]^{\times}/\k^{\times}))
\to H^1_{e_{X}}(X,S) \to H^0(k,H^1({\overline X},S))
  \to H^2(k, \Hom_{\Z}(S^*,
\k[X]^{\times}/\k^{\times})).$$

Soient $X$ et $Y$ deux $k$-vari\'et\'es lisses,
g\'eom\'etriquement int\`egres, g\'eom\'etriquement rationnelles.
Alors le foncteur qui \`a $X$ associe le module galoisien 
$H^1({\overline X},S)$ est additif. 
On \'etablit ce r\'esultat par r\'eduction au cas
$S=\G_m$ et au cas $G=\mu_n$. Ce dernier cas se traite par
la suite de Kummer : on utilise le r\'esultat pour $H^1({\overline X},\G_m)$
(et dans ce cas on utilise le fait que les vari\'et\'es sont
g\'eom\'etriquement rationnelles, voir [CTSa77], Lemme 11 p. 188)
et le r\'esultat pour $H^0({\overline X},\G_m)/n$ (Rappel {\bf 0.1}).

\medskip

{\bf Lemme 5.5} {\it 
Soit $S$ un $k$-groupe de type multiplicatif.
Sur la cat\'egorie des
$k$-vari\'et\'es lisses point\'ees g\'eom\'etriquement connexes,
g\'eom\'etriquement rationnelles, le foncteur
qui \`a $X$ associe $H^1_{e_{X}}(X,S)$ est additif.
}

\medskip

En d'autres termes, \'etant donn\'ees deux telles $k$-vari\'et\'es
point\'ees $(X, e_{X})$ et $(Y,e_{Y})$, tout torseur ${\cal T}$ sur $X\times_kY$ sous $S$
est isomorphe \`a l'image, par la fl\`eche de multiplication
$S \times_k S \to S$, des torseurs $p_X^*\circ i_X^* ({\cal T})$
et $p_Y^*\circ i_Y^* ({\cal T})$, o\`u $p_X$ et $p_Y$ sont les projections
de $X \times_kY$ sur les deux facteurs, et o\`u $i_X : X \times {e_{Y}}  
\hookrightarrow
X
\times_kY$, resp. $i_Y : {e_{X}} \times_kY \hookrightarrow X
\times_kY$ sont les immersions naturelles.

\medskip
{\it D\'emonstration}
La somme des suites exactes ci-dessus pour $X$ et $Y$ s'envoie
naturellement dans la suite exacte pour $X\times_kY$.
L'additivit\'e des foncteurs $\Hom_{\Z}(S^*, \k[X]^{\times}/\k^{\times})$
et $H^1({\overline X},S)$ et le lemme des 5 donnent  le r\'esultat.
\cqfd

\bigskip

{\bf Th\'eor\`eme 5.6}
 {\it  Soient $k$ un corps, $G$ un $k$-groupe
alg\'ebrique  lin\'eaire lisse connexe et
$S$ un $k$-groupe de type multiplicatif lisse.  Soit $p : Y \to G$ un torseur sur
$G$ sous $S$, de fibre triviale en l'\'el\'ement neutre $e_{G} \in G(k)$.
Il existe alors une structure de $k$-groupe alg\'ebrique
lin\'eaire sur $Y$, telle que
de plus $p : Y \to G$ soit un   homomorphisme de $k$-groupes alg\'ebriques,
de noyau $S$. }

\medskip

{\it D\'emonstration}
La $k$-vari\'et\'e lisse $G$ est   g\'eom\'etriquement rationnelle.

Soit $p : H  \to G$ un torseur sur $G$ sous le $k$-groupe de type multiplicatif lisse
$S$. Supposons de plus $H$ \'equip\'e d'un $k$-point $e_{H}$ s'envoyant
sur   $e_{G} \in G(k)$.

Consid\'erons la multiplication 
$$ m : G \times G \to G.$$
D'apr\`es le lemme 5.5, tout \'el\'ement $\alpha \in H^1_{e_{G}\times e_{G} }(G \times G,S)$
s'\'ecrit $p_1^*(\alpha_{G \times e_{G}}) + p_2^*(\alpha_{e_{G} \times G})$. 
Ainsi pour tout \'el\'ement $\alpha \in H^1_{e_{G}}(G,S)$,
on a $m^*(\alpha)=p_1^*(\alpha)+p_2^*(\alpha)$.
On d\'ecrit ce fait en disant que  tout \'el\'ement de $H^1_{e_{G}}(G,S)$ est ``primitif".

Ceci \'etablit  l'existence
d'un diagramme commutatif de $k$-morphismes
de vari\'et\'es point\'ees
$$\diagram { \varphi &  :  & H                         &\times &H                        & \to      & H \cr
                                      &    & \downarrow{p}{}  &           &\downarrow{p}{}  &  & \downarrow{p}{} \cr
                               m &  :  & G                         &\times &G                        & \to      & G, }
$$
o\`u $m$ d\'esigne la multiplication de $G$ et $\varphi$ est compatible 
avec
la multiplication $S \times S\to S$,
en ce sens que l'on a 
$$\varphi(t_1h_1,t_2h_2)=t_1t_2\varphi(h_1,h_2).$$

On notera qu'un  diagramme comme ci-dessus
peut \^etre modifi\'e par un $k$-morphisme
$H \to H$ donn\'e par $h \to f(h).h$, o\`u
$f$ est
le compos\'e de la
projection $H \to G$ par un morphisme point\'e
$G \to S$, c'est-\`a-dire par un 
homomorphisme de $G$ dans $S$ ({\bf 0.2}, lemme de Rosenlicht).
Ceci simultan\'ement sur chaque facteur $H \to G$, les deux facteurs
\`a gauche et celui \`a droite.

Le morphisme $\varphi(e_{H},h) : H \to H$ est un morphisme de
$S$-torseurs point\'es sur $G$. On a donc
$\varphi(e_{H},h)=\chi_2(p(h)).h,$ o\`u $p : H \to G$
est la projection naturelle, et 
$\chi_2 : G \to S$ est un caract\`ere. 
De m\^eme, 
$\varphi(h,e_{H})=\chi_1(p(h)).h,$
avec $\chi_1 : G \to S$  un caract\`ere.

\medskip
Si l'on remplace $\varphi(h_1,h_2)$ par
$$(\chi_1(p(h_1)))^{-1}. (\chi_2(p(h_2)))^{-1}. \varphi(h_1,h_2),$$
on a maintenant un nouveau morphisme $\varphi = H \times H \to H$ 
satisfaisant
$\varphi(e_{H},h)=h=\varphi(h,e_{H})$. Autrement dit, $e_{H}$ est un \'el\'ement neutre pour
l'op\'eration $\varphi$.
Noter que ce morphisme donne encore lieu \`a un diagramme
commutatif comme ci-dessus, i.e. la projection $H \to G$
est compatible avec l'op\'eration $\varphi : H \times H \to H$,
avec l'action de $S \times S \to S$
et (par projection) avec la multiplication $G \times G \to G$.

\medskip

A-t-on associativit\'e ?
Soient $h_1,h_2,h_3 \in H$.  En utilisant le diagramme ci-dessus,
l'associativit\'e de la multiplication pour $G$ et
le lemme de Rosenlicht,
on voit qu'il existe
un $k$-morphisme de groupes $\sigma : H \times_kH\times_kH \to S$
tel que
$$\varphi(h_1,\varphi(h_2,h_3))= \sigma(h_1,h_2,h_3). 
\varphi(\varphi(h_1,h_2),h_3) \in H.$$

Comme $\varphi$ est compatible avec l'action de $S$,
on voit que $\sigma(t_1h_1,t_2h_2,t_3h_3)=\sigma(h_1,h_2,h_3).$
Ainsi $\sigma$ passe au quotient par l'action de 
$S \times_kS\times_kS$, et est induit par
un $k$-morphisme de groupes  $\tau : G \times_k G\times_k G \to S$.
Une autre application du lemme de Rosenlicht 
assure  que $\tau(g_1,g_2,g_3)= \chi_1(g_1).\chi_2(g_2).\chi_3(g_3),$
o\`u les $\chi_i : G \to S$ sont des caract\`eres.
 
En utilisant la propri\'et\'e
$\varphi(e_{H},h)=h=\varphi(h,e_{H})$, on voit que l'on a $\tau(e_{G},g_2,g_3)=1 \in S$.
Ceci \'etablit $\chi_2=1$ et $\chi_3=1$. Par sym\'etrie, on obtient de m\^eme 
$\chi_1=1$. Ceci \'etablit l'associativit\'e de $\varphi$.

\medskip

Il reste \`a \'etablir l'existence de l'inverse. 
On va adapter l'argument
de Serre ([Se59], bas de la page 183). On commence par montrer
que l'application $g \mapsto g^{-1}$ induit l'application
$i_G : x \to -x$ sur le groupe $H^1_{e_{G}}(G,S)$.
Ceci n'est pas formel, sur une courbe elliptique,
il y a des \'el\'ements du groupe de Picard de degr\'e non nul
invariants par $g\mapsto -g$. Mais ici, tout \'el\'ement
 de $H^1_{e_{G}}(G,S)$ est ``primitif". 
 Notons $\alpha \in H^1_{e_{G}}(G,S)$ la classe du torseur point\'e $p : H \to G$.
 Le compos\'e de
$G \to G\times G$ donn\'e par $u \mapsto (u,u^{-1})$
et de la multiplication $m : G \times G \to G$
est la fl\`eche constante \'egale \`a 1.
Du fait que $\alpha$ est primitif, on d\'eduit
imm\'ediatement $(-1)^*(\alpha)=-\alpha$.
Par ailleurs, de fa\c con triviale, 
la fl\`eche $ s \mapsto s^{-1}$ dans $S$ 
induit la fl\`eche $i_S : x \mapsto -x$ dans
$H^1_{e_{G}}(G,S)$. 
On a donc $i_G^*(\alpha)=i_S^*(\alpha)$
et on trouve ainsi un morphisme point\'e
$ \theta : H \to H$ au-dessus de l'application inverse $G \to G$,
tel que de plus $i(sh)=s^{-1}h$. Ce morphisme est un
isomorphisme de $k$-vari\'et\'es, car il est induit par la projection
$H \times_{G}G \to H$, la fl\`eche $G \to G$ \'etant l'isomorphisme
de $k$- vari\'et\'es donn\'e par $g \mapsto g^{-1}$.
Le morphisme $\lambda : x \mapsto \varphi(\theta(x),x)$ envoie
$H$ dans l'image r\'eciproque 
de $1_G$ par $H \to G$, c'est-\`a-dire dans $S$.
Ce morphisme satisfait $\lambda(sx)=\lambda(x)$.
Il d\'efinit donc un morphisme point\'e de
$G$ dans $S$, c'est-\`a-dire un caract\`ere $\chi : G \to S$.
Tout caract\`ere $\chi : G \to S$
d\'efinit un isomorphisme de $S$-torseurs
$H \to H$ via $x \to \chi(p(x)).x$.
Rempla\c cons alors le $k$-isomorphisme point\'e (de $k$-vari\'et\'es) $\theta$ par $\xi : H \to H$ d\'efini
par le $k$-isomorphisme point\'e (de $k$-vari\'et\'es)
$\xi(x)=\chi(p(x))^{-1}.\theta(x)$.
On a alors $\varphi(\xi(x),x)=1$. 
On a donc d\'efini un inverse \`a gauche
pour la loi de composition $\varphi$ sur $H$.
Comme chacun sait, c'est alors un inverse
\`a droite pour la loi $\varphi$, qui est
associative et poss\`ede un \'el\'ement neutre.
\cqfd
\medskip

{\it Remarques 5.6.1}

(1) Dans [Se59] (Chapitre VII, \S 15),
Serre \'etudie  les torseurs primitifs sur une vari\'et\'e ab\'elienne $A$,
sous un groupe lin\'eaire connexe (commutatif) $B$.
Serre utilise le fait que tout morphisme de $A$ dans $B$ est
constant. Cette propri\'et\'e  vaut encore si la vari\'et\'e ab\'elienne $A$ est
remplac\'ee par un groupe $G$ semi-simple connexe. Mais elle ne vaut plus
pour $G$  un groupe r\'eductif connexe quelconque.

(2)
 Le cas particulier du th\'eor\`eme 5.6  o\`u $k$ est
alg\'ebriquement clos et $S$ est un $k$-groupe ab\'elien fini  est le th\'eor\`eme 2 de   [Mi72].

(3) Le th\'eor\`eme 5.6 ne s'\'etend pas \`a des groupes $S$ quelconques.
M. Florence m'a donn\'e des contre-exemples, sur un corps $k$  alg\'ebriquement clos,
avec $S=PGL_{n}$ et $G=\G_{m}^2$.

\medskip

On d\'eduit du th\'eor\`eme le r\'esultat essentiellement \'equivalent suivant.

\medskip

{\bf Corollaire 5.7}
{\it Pour tout $k$-groupe lin\'eaire lisse connexe $G$ et tout $k$-groupe de type multiplicatif lisse $S$
 la fl\`eche naturelle  
$$\Ext_{k-gp}(G,S) \to \Ker[H^1(G,S) \to H^1(k,S)],$$
est un isomorphisme. }

{\it D\'emonstration} 
La fl\`eche dont il est ici question est la suivante. La donn\'ee d'une extension de $k$-groupes alg\'ebriques 
$$1 \to S \to G_{1} \to G \to 1$$ 
donne par oubli un torseur sur $G$ sous $S$, dont la restriction au-dessus de l'\'el\'ement neutre $e_{G} \in G(k)$  poss\`ede un $k$-point. Cette application induit un homomorphisme
de groupes ab\'eliens
$\Ext_{k-gp}(G,S) \to  \Ker[H^1_{\et}(G,S) \to H^1(k,S)]$ (la fl\`eche $H^1_{\et}(G,S)=H^1(G,S) \to H^1(k,S)$
\'etant d\'efinie par la restriction \`a $e_{G}$),
 dont le th\'eor\`eme 5.6 montre qu'il est surjectif.
L'argument donn\'e dans la d\'emonstration de la proposition 3.2 
montre que cet homomorphisme est injectif. \cqfd

\bigskip

{\bf \S 6. Le groupe fondamental alg\'ebrique des groupes alg\'ebriques lin\'eaires connexes via les r\'esolutions flasques}

\bigskip

{\bf Proposition-D\'efinition  6.1}  {\it Soit $G$ un $k$-groupe lin\'eaire
connexe, suppos\'e r\'eductif si ${\rm car}(k)>0$. Soit
$1 \to S \to H \to G \to 1$ une r\'esolution flasque de $G$.
Soit $P$ le $k$-tore quasi-trivial $H^{tor}$. On dispose donc
d'un morphisme de $k$-tores $S \to P$.
Le  module galoisien discret de type fini ${\rm Coker} [S_{*} \to P_{*}]$ ne d\'epend pas
de la r\'esolution flasque choisie. On le note
 $\pi_{1}(G)$ et on l'appelle le groupe fondamental alg\'ebrique de $G$.}

\medskip

{\it D\'emonstration} Le seul point \`a \'etablir est que  le
quotient ${\rm Coker} [S_{*} \to P_{*}]$ est canoniquement d\'efini.
C'est l'objet de la proposition 3.2 (iv) et de la remarque 3.2.2. \cqfd

\bigskip

{\bf Proposition 6.2} {\it 
La suite  
$$0 \to S_{*} \to P_{*} \to \pi_{1}(G) \to 0   $$
est exacte. C'est une r\'esolution coflasque du module galoisien $\pi_{1}(G)$.}

\medskip

{\it D\'emonstration} La fl\`eche $S \to P$   a un noyau fini (proposition 3.3),
donc la fl\`eche $S_{*} \to P_{*}$ est injective.
Ceci \'etablit l'exactitude de la suite.
Comme $S_{*}$ est un module coflasque et $P_{*}$ un module de permutation, ceci
est une r\'esolution coflasque du module galoisien de type fini $ \pi_{1}(G)$
([CTSa77], \S 1; CTSa87a], \S 0).\cqfd

\bigskip

{\it Remarques 6.2.1}

(1) On verra dans  l'appendice A que le groupe fondamental alg\'ebrique ici d\'efini co\"{\i}ncide avec 
le groupe d\'efini par Borovoi dans [Bo96] et [Bo98], groupe qui sera ici not\'e $\pi_{1}^{Bor}(G)$.

(2) Compte tenu des propri\'et\'es des r\'esolutions coflasques ([CTSa87a], Lemma 0.6),
le module galoisien
$S_{*}$ est d\'etermin\'e 
 \`a addition pr\`es d'un module de permutation  par le module galoisien $\pi_{1}(G)$.

(3) Borovoi et Kunyavski\u{\i} [BoKu04] partent de la donn\'ee du module galoisien
$\pi_1^{Bor}(G)$, et consid\`erent (op. cit., \S 1) une r\'esolution
coflasque du module galoisien $\pi_1^{Bor}(G)$ :
$$0 \to S_{*} \to P_{*} \to \pi_{1}^{Bor}(G) \to 0. $$
Ceci d\'efinit 
$S_*$ \`a addition pr\`es d'un module de permutation.
En caract\'eristique nulle, 
pour $G \subset X$ une $k$-compactification lisse,
ils montrent que le module galoisien $\Pic({\overline X})$ est flasque
(op. cit., \S 2; th\'eor\`eme {\bf 0.9} ci-dessus ). Ils montrent enfin (op. cit., \S 3, Thm. 3.21 p. 309)
que  $\Pic({\overline X})$ et $S^*$ sont isomorphes 
\`a addition 
pr\`es de modules de permutation.

\bigskip

\bigskip

On note $M_{\g}$ le groupe ab\'elien des coinvariants d'un 
 $\g$-module continu discret, c'est-\`a-dire le plus grand $\frak{g}$-quotient de $M$ sur lequel
 le groupe de Galois $\g$ agit trivialement.
 
 \medskip

{\bf Proposition 6.3} {\it Soit $G$ un $k$-groupe lin\'eaire connexe, suppos\'e r\'eductif si ${\rm car}(k) > 0$.
Les groupes ab\'eliens  finis $\Pic(G)$ et $(\pi_{1}(G)_{\g})_{tors}$ sont duaux l'un de l'autre.}

\medskip

{\it D\'emonstration}  
Soit
$1 \to S \to H \to G \to 1$ une r\'esolution flasque de $G$.
Soit $P$ le $k$-tore quasi-trivial $H^{tor}$.
On dispose  de la suite exacte de la proposition 6.2 :
$$ 0 \to S_{*} \to P_{*} \to \pi_{1}(G) \to 0.$$
Cette derni\`ere suite donne naissance de fa\c con \'evidente au diagramme commutatif  
$$\diagram{
&&0 
&&0&&0 && \cr
&&{\downarrow} &&{\downarrow}&&{\downarrow}&& \cr
&& \Hom_{\g}(\pi_{1}(G),\Z) & \to & \Hom_{\g}(\pi_{1}(G),\Q) & \to & \Hom_{\g}(\pi_{1}(G),\Q/\Z) & & \cr
&&{\downarrow}&&{\downarrow}&&{\downarrow}&& \cr
0& \to& \Hom_{\g}(P_{*},\Z) & \to & \Hom_{\g}(P_{*},\Q) & \to & \Hom_{\g}(P_{*},\Q/\Z)) & \to & 
H^1(\g,\Hom_{\Z} (P_{*},\Z))\cr
&&{\downarrow}&&{\downarrow}&&{\downarrow}&& \cr
0& \to & \Hom_{\g}(S_{*},\Z) & \to & \Hom_{\g}(S_{*},\Q) & \to & \Hom_{\g}(S_{*},\Q/\Z) & \to & H^1(\g,\Hom_{\Z} (S_{*},\Z))\cr
&&{\downarrow}&&{\downarrow}&&&& \cr
&&\pic(G)&& 0 &&&& \cr
&&{\downarrow}&&&&&& \cr
&&0 &&&&&&}
$$
o\`u le conoyau de $ \Hom_{\g}(P_{*},\Z) \to  \Hom_{\g}(S_{*},\Z)$, c'est-\`a-dire de
$(P^*)^{\g} \to (S^*)^{\g}$, a \'et\'e identifi\'e, gr\^ace \`a la proposition 3.3,
\`a $\Pic(G)$. Que $\Hom_{\g}(P_{*},\Q)  \to \Hom_{\g}(P_{*},\Q)$ soit surjectif r\'esulte simplement
du fait que  $\Ext^1_{\g}(M,\Q)=0$ pour tout $\g$-module continu discret. Comme $P_{*}$ et donc aussi
$P^*=\Hom_{\Z}(P_{*},\Z)$ sont des modules de permutation, on a $H^1(\g,\Hom_{\Z} (P_{*},\Z))=0$.

Le lemme du serpent donne alors
la suite exacte
$$ \Hom_{\g}(\pi_{1}(G),\Q)   \to   \Hom_{\g}(\pi_{1}(G),\Q/\Z) \to \Pic(G) \to 0,$$
soit encore
$$\Hom_{\Z}(\pi_{1}(G)_{\g}, \Q) \to \Hom_{\Z}(\pi_{1}(G)_{\g}, \Q/\Z)  \to \Pic(G) \to 0.$$
Il reste alors \`a appliquer le lemme g\'en\'eral suivant. Pour tout groupe ab\'elien  de type fini $A$,
on a une suite exacte naturelle
$$ \Hom_{\Z}(A,\Q) \to \Hom_{\Z}(A,\Q/\Z) \to \Hom_{\Z}(A_{tors},\Q/\Z) \to 0. \hskip3mm \qed$$

\bigskip

{\it Remarque 6.3.1}  
Soient $k$ un corps de caract\'eristique nulle 
et $G$ un $k$-groupe r\'eductif connexe. Comme le montre Borovoi ([Bo98], Prop. 1.10), le groupe  $\pi_{1}^{Bor}(G)$ admet aussi une interpr\'etation 
en termes du centre d'un groupe dual de Langlands (connexe) de $G$.
Plus pr\'ecis\'ement, 
 le groupe $\Hom_{\Z}(\pi_{1}^{Bor}(G),\C^*)$, muni de sa structure de $\g$-module
 via l'action sur $\pi_{1}(G)$, s'identifie au groupe des points complexes $Z({\hat G})(\C)$ du centre  $Z({\hat G})$
 d'un dual de Langlands connexe de $G$. On  a donc 
 $$(Z({\hat G})(\C))^{\g} \simeq \Hom_{\g}(\pi_{1}^{Bor}(G),\C^*)=\Hom((\pi_{1}^{Bor}(G))_{\g}, \C^*)$$
 d'o\`u il r\'esulte que le groupe ab\'elien fini $\pi_{0}((Z({\hat G})(\C))^{\g})$ des composantes connexes 
 de $(Z({\hat G})(\C))^{\g}$ utilis\'e
par Kottwitz dans [Ko86]
 s'identifie \`a $\Hom((\pi_{1}^{Bor}(G)_{\g})_{tors}, \C^*)$ soit au dual du groupe ab\'elien fini
 $(\pi_{1}^{Bor}(G)_{\g})_{tors}$.  La proposition 6.3  et l'identification
 $\pi_{1}(G)=\pi_{1}^{Bor}(G)$ faite dans l'appendice A montrent  alors
que le groupe $ \pi_{0}((Z({\hat G})(\C))^{\g})$  
n'est autre que le groupe
 $ \Pic(G)$. Ce dernier r\'esultat est d\'ej\`a \'etabli par Kottwitz 
 dans  [Ko84] (2.4.1).

 \bigskip
 
 Etant donn\'e un $k$-groupe fini de type multiplicatif $\mu$ et
$ {\mu^*}=\Hom_{k-gp}(\mu,\G_{m,k})$ son groupe des caract\`eres,
qui est un $\g$-module fini, on note $\mu(-1)$ le $\g$-module fini $\Hom_{\Z}({\mu^*}, \Q/\Z)$. On peut aussi d\'efinir
 ce groupe de la fa\c con suivante.
 Soit $N>0$ un entier annulant ${\mu^*}$.
On a $$\Hom_{\Z}({\mu^*}, \Q/\Z)=\Hom_{\Z}({\mu^*}, \Z/N)=\Hom_{k-gp}(\mu_{N},\mu),$$
o\`u la derni\`ere \'egalit\'e  est obtenue par dualit\'e. Ainsi $\mu(-1)=\Hom_{k-gp}(\mu_{N},\mu).$
  
  \medskip
  
  {\bf Proposition 6.4} {\it Soit $G$ un $k$-groupe lin\'eaire connexe, suppos\'e r\'eductif si
  ${\rm car}(k)>0$. Soit  $\mu$ le noyau de $G^{sc} \to G^{ss}$.
On a la suite exacte naturelle
de modules galoisiens
$$ 0 \to \mu(-1) \to \pi_{1}(G) \to G^{tor}_{*} \to 0.$$}
{\it D\'emonstration}
Soit
$1 \to S \to H \to G \to 1$ une r\'esolution flasque de $G$.
Soit $P$ le $k$-tore quasi-trivial $H^{tor}$.
Notons $T=G^{tor}$. Le diagramme fondamental du  paragraphe 3 donne  naissance \`a la suite exacte de $k$-groupes de type
multiplicatif 
$$ 1 \to \mu \to S \to P \to T \to 1.$$
La suite des groupes de caract\`eres associ\'ee est la suite exacte de modules galoisiens
$$ 0 \to T^* \to S^* \to P^* \to {\mu}^* \to 0.$$
On peut couper cette suite exacte en deux suites exactes courtes
$$ 0 \to T^* \to S^*  \to R^* \to 0$$
et 
$$0 \to R^*\to P^* \to {\mu}^* \to 0,$$
o\`u $R^*, S^*,T^*$ sont sans $\Z$-torsion, et ${\mu^*}$ est un groupe de torsion.

Appliquant \`a ces deux suites le foncteur $\Hom_{\Z}(\bullet,\Z)$, et combinant  les suites exactes courtes
obtenues, on obtient une suite exacte naturelle
$$ 0 \to \Ext^1_{\Z}({\mu^*},\Z) \to {\rm Coker }(S_{*} \to P_{*}) \to T_{*} \to 0,$$
soit encore
$$ 0 \to \Ext^1_{\Z}({\mu^*},\Z) \to \pi_{1} (G) \to T_{*} \to 0,$$
La consid\'eration de la suite exacte
$$0 \to \Z \to \Q \to \Q/\Z \to 0$$
donne un isomorphisme $\Hom_{\Z}({\mu^*}, \Q/\Z)=\Ext^1_{\Z}({\mu^*}, \Z)$.
\cqfd

\bigskip

Les deux cas extr\^emes de l'\'enonc\'e pr\'ec\'edent, cas qu'on peut bien s\^ur \'etablir directement,
 sont celui o\`u
$G=T$ est un $k$-tore, auquel cas $\pi_{1}(T)= T_{*}$ et celui o\`u $G$ est
un $k$-groupe semi-simple, auquel cas $\pi_{1}(G)=\mu(-1)$. On voit donc que pour $G$ un
 $k$-groupe semi-simple, le groupe fondamental alg\'ebrique $\pi_{1}(G)$ diff\`ere  du groupe fondamental  $\mu=\Ker [G^{sc} \to G]$ par une torsion galoisienne.

\medskip

{\it Remarque 6.4.1} Soit
$L/k$ la plus petite extension (automatiquement) galoisienne d\'eployant (au sens
des groupes de type multiplicatif)  \`a la fois
$G^{tor}$ et le centre de $G^{sc}$. 
On a vu (Remarque 3.1.1) qu'il existe une r\'esolution flasque $1 \to S \to H \to G \to 1$
telle que $S$ et $H^{tor}$ soient d\'eploy\'es par $L/k$. De la d\'efinition de $\pi_1(G)$  il r\'esulte
donc que le module galoisien $\pi_1(G)$ est d\'eploy\'e par $L$. Soit $K \subset L$ la plus petite sous-extension (automatiquement) galoisienne de $k$ d\'eployant le module galoisien $\pi_1(G)$. 
On n'a pas forc\'ement $K=L$, comme l'on voit d\'ej\`a en prenant  pour  $G$ 
un descendu \`a la Weil d'un  groupe $SL_{n}$.
Existe-t-il une r\'esolution flasque 
$1 \to S \to H \to G \to 1$ telle que $S$ et $H^{tor}$ soient d\'eploy\'es par $K/k$ ?
De la proposition 6.4 il r\'esulte que $K$ d\'eploie
\`a la fois $G^{tor}$ et le $k$-groupe de type multiplicatif $\mu$ noyau de $G^{sc} \to G^{ss}$.  L'extension $K/k$ est-elle la plus petite extension (galoisienne) satisfaisant ces propri\'et\'es ?

\bigskip

{\bf Proposition 6.5} {\it Soit $G$ un $k$-groupe   lin\'eaire connexe, suppos\'e r\'eductif si
${\rm car}(k) >0$.

(i) On a $\pi_{1}(G)=\pi_{1}(G^{\red})$.

(ii)  Si $G$ est r\'eductif, il  est semi-simple si et seulement si
 $\pi_{1}(G)$ est fini.

(iii)  Le groupe $G$ est   simplement connexe si et seulement
si  $\pi_{1}(G)=0$.
 
(iv) Le groupe   $\pi_{1}(G^{ssu})=\pi_{1}(G^{ss})$
s'identifie \`a la torsion de $\pi_{1}(G)$.

(v)  Le groupe $\pi_{1}(G)$  est sans torsion si et seulement si $G^{ssu}$ est
simplement connexe, si et seulement si $G^{ss}$ est simplement connexe.

(vi) Le $k$-groupe   $G$ est quasi-trivial si et seulement si $\pi_1(G)$ est un module de permutation.

(vii) Le $k$-groupe   $G$ est coflasque si et seulement si $\pi_1(G)$ est un module coflasque.
}

\medskip

{\it D\'emonstration}  
Rappelons que $G$ est quasi-trivial si et seulement si $\mu=1$ et $T=G^{tor}$
est un $k$-tore quasi-trivial. De m\^eme $G$ est coflasque si et seulement si $\mu=1$
et $T=G^{tor}$ est un $k$-tore coflasque.
L'\'enonc\'e est une cons\'equence imm\'ediate du diagramme fondamental (\S 3)
et de la proposition 6.4.
 \cqfd

\bigskip 

{\bf Proposition 6.6}  {\it Soient $G_{1},G_{2},G_{3}$ des $k$-groupes lin\'eaires connexes,
suppos\'es r\'eductifs si ${\rm car}(k)>0$.

(i) Tout   morphisme de $k$-groupes alg\'ebriques $\lambda : G_{1 } \to G_{2}$   
 induit un homomorphisme naturel
de modules galoisiens $\lambda_{* } : \pi_{1}(G_{1}) \to \pi_{1}(G_{2})$. 

(ii) Si $ \lambda : G_{1 } \to G_{2}$ et $\mu : G_{2 } \to G_{3}$ sont des
 morphismes de $k$-groupes alg\'ebriques, alors les homomorphismes
$(\mu \circ \lambda)_{*}$ et $ \mu_{*}\circ \lambda_{*} \in \Hom_{\g}(\pi_{1}(G_{1}),\pi_{1}(G_{3}))$ co\"{\i}ncident.
}

\medskip

{\it D\'emonstration} 
Soit  $1 \to S_{1} \to H_{1} \to G_{1} \to 1$ une r\'esolution flasque de $G_{1}$ et $1 \to S_{2} \to H_{2} \to G_{2} \to 1$ une r\'esolution flasque de $G_{2}$. Soit $E=H_{1}\times_{G_{2}}H_{2}$ le produit fibr\'e, o\`u l'homomorphisme $H_{1} \to G_{2}$ est le compos\'e de $H_{1} \to G_{1}$ et de
$\lambda : G_{1} \to G_{2}$. On a la suite exacte de $k$-groupes 
$$1 \to S_{2} \to E \to H_{1} \to 1.$$
Comme $H_{1}$ est un groupe quasi-trivial et $S_{2}$ est un $k$-tore flasque,
cette suite exacte est scind\'ee (voir la d\'emonstration de la proposition 3.2).
Il existe donc un diagramme commutatif de morphismes de $k$-groupes
$$\diagram{
1 &\to &S_{1} & \to & H_{1} & \to &G_{1} &\to &1 \cr
&&{\downarrow}   && {\downarrow} && {\downarrow}&&                                 \cr
1 &\to &S_{2} &\to& H_{2} &\to &G_{2} &\to& 1.\cr
}
$$
Deux choix diff\'erents de $H_{1} \to H_{2}$ diff\`erent par un $k$-morphisme  de
$H_{1}$ dans $S_{2}$ (par le lemme de Rosenlicht, c'est n\'ecessairement un $k$-morphisme de groupes alg\'ebriques).
Le diagramme ci-dessus  induit un diagramme commutatif de morphismes de
modules galoisiens
$$\diagram{(S_{1})_{*} & \to & (H_{1})_{*} \cr
{\downarrow} &  & {\downarrow} \cr
(S_{2})_{*} & \to & (H_{2})_{*}
}$$
Ceci d\'efinit un homomorphisme de modules galoisiens
$\pi_{1}(G_{1}) \to \pi_{1}(G_{2})$. Comme deux choix diff\'erents de $H_{1} \to H_{2}$
induisent des homomorphismes $(H_{1})_{*} \to (H_{2})_{*}$ qui diff\`erent  par l'image d'un
 homomorphisme de $(H_{1})_{*}$ dans $(S_{2})_{*}$, on voit que cet homomorphisme
est ind\'ependant du choix de $H_{1} \to H_{2}$. Que l'homomorphisme ne d\'epende pas
du choix des deux r\'esolutions a d\'ej\`a \'et\'e \'etabli (voir la d\'efinition 6.1). 
Ceci \'etablit la proposition.\cqfd

\bigskip

{\bf Proposition 6.7}  {\it Soit $G$ un $k$-groupe r\'eductif connexe. Soit $n>0$ un entier inversible dans
$k$. On a un isomorphisme naturel $H^1({\overline G},\mu_{n}) \simeq {\rm Hom}_{\Z}(\pi_{1}(G),\Z/n)$.}

\medskip

{\it D\'emonstration}  
Soit
$1 \to S \to H \to G \to 1$ une r\'esolution flasque de $G$.
Notons $p$ l'homomorphisme $H \to G$.
Soit $P$ le $k$-tore quasi-trivial $H^{tor}$.
 De la proposition 6.2 on tire  l'isomorphisme
$$ {\rm Hom}_{\Z}(\pi_{1}(G),\Z/n) \buildrel \simeq \over \rightarrow\Ker [P^*/n \to S^*/n], $$
soit encore 
$$ {\rm Hom}_{\Z}(\pi_{1}(G),\Z/n) \buildrel \simeq \over \rightarrow \Ker [H^*/n \to S^*/n], $$
compte tenu de l'isomorphisme $P^* \buildrel \simeq \over \rightarrow H^*$ (cons\'equence de {\bf 0.5}).

De la suite exacte
$$1 \to S \to H \to G \to 1$$
on tire la suite exacte de faisceaux  \'etales sur $G$ :
$$ 1 \to \G_{m,G} \to p_{*}\G_{m,H} \to S^* \to 0$$
([CTSa87b], Prop.  1.4.2 p. 383). Localement pour la topologie  \'etale sur $G$,
la fibration $H \to G$ s'identifie \`a un tore d\'eploy\'e au-dessus d'un sch\'ema r\'egulier. On a donc 
$R^1p_{*}(\G_{m,H})=0$.
Le lemme du serpent appliqu\'e \`a la multiplication par $n$
sur la suite exacte ci-dessus donne donc un isomorphisme
$$ R^1p_{*}\mu_{n} \buildrel \simeq \over \rightarrow  S^*/n.$$
On a par ailleurs de fa\c con \'evidente $\mu_{n,G} \buildrel \simeq \over \rightarrow p_{*}\mu_{n,H}$.
La suite des termes de bas degr\'e de la suite spectrale de Leray pour le morphisme
$p$ et le faisceau $\mu_{n,H}$, sur $\k$, s'\'ecrit donc
$$0 \to H^1({\overline G},\mu_{n}) \to H^1({\overline H},\mu_{n}) \to  S^*/n.$$
On a $\pic({\overline H})=0$ et $P^* \buildrel \simeq \over \rightarrow H^*$.
De la suite de Kummer sur ${\overline H}$
on d\'eduit donc un isomorphisme naturel $  H^*/n \buildrel \simeq \over \rightarrow H^1({\overline H},\mu_{n}) $.
On a donc 
$$ H^1({\overline G},\mu_{n}) \buildrel \simeq \over \rightarrow {\rm Ker} [H^*/n \to S^*/n]$$
ce qui \'etablit la proposition.
\cqfd

\medskip

{\it Remarque 6.7.1} Supposons $k$ de caract\'eristique nulle. On sait ([Mi72], Thm. 1)
 que le groupe fondamental de Grothendieck $\pi_{1}^{Groth}({\overline G})$ est un groupe profini {\it ab\'elien} (noter que le groupe fondamental ``alg\'ebrique'' qui intervient dans le titre de l'article de Miyanishi est le groupe fondamental de Grothendieck, ce n'est pas le groupe fondamental alg\'ebrique
 consid\'er\'e  ici).
La proposition ci-dessus implique que le module galoisien  $$\pi_{1}^{Groth}({\overline G})(-1)={\rm Hom}({\hat {\Z}}(1),\pi_{1}^{Groth}({\overline G}))$$
est naturellement le compl\'et\'e profini du module galoisien de type fini $\pi_{1}(G)$ d\'efini ici.

\bigskip

{\bf Proposition 6.8} {\it  Soit $k$ un corps de caract\'eristique nulle.
Si $1 \to G_{1} \to G_{2 } \to G_{3} \to 1$
est une suite exacte de $k$-groupes lin\'eaires connexes, on a une suite exacte courte induite
de modules galoisiens
$$ 0 \to \pi_{1}(G_{1}) \to  \pi_{1}(G_{2}) \to   \pi_{1}(G_{3})   \to 0.$$
}

{\it D\'emonstration}  
D'apr\`es {\bf 0.5}, on a la suite exacte naturelle
$$ 0 \to G_{3}^* \to G_{2}^* \to G_{1}^* \to \pic({\overline G}_{3}) \to \pic({\overline G}_{2}) \to
\pic({\overline G}_{1}) \to 0.$$
On a $  T_{i}^* \buildrel \simeq \over \rightarrow G_{i}^*  $, o\`u $T_{i}=G_{i}^{tor}$, et
on a 
$ \mu_{i}^*  \buildrel \simeq \over \rightarrow \pic({\overline G}_{i}^{ss}) \simeq 
\pic({\overline G}_{i})$, ceci de fa\c con fonctorielle.
Le groupe $\pic({\overline G}_{3})$ \'etant fini, on voit que la fl\`eche naturelle
$(T_{1})_{*} \to (T_{2})_{*}$ est injective.  Par ailleurs, la fl\`eche
$\Hom_{\Z}(\mu_{1}^*, \Q/\Z) \to \Hom_{\Z}(\mu_{2}^*, \Q/\Z)$, 
qui s'identifie \`a la fl\`eche naturelle $\mu_{1}(-1) \to \mu_{2}(-1)$, est injective.
On v\'erifie facilement qu'un diagramme
de r\'esolutions flasques
$$\diagram{
1 &\to &S_{1} & \to & H_{1} & \to &G_{1} &\to &1 \cr
&&{\downarrow}   && {\downarrow} && {\downarrow}&&                                 \cr
1 &\to &S_{2} &\to& H_{2} &\to &G_{2} &\to& 1\cr
}
$$
comme consid\'er\'e \`a la proposition 6.6 induit un diagramme commutatif
$$ \diagram{
 0 &\to &\mu_{1}(-1) &\to &{\rm Coker }[(S_{1})_{*} \to (P_{1})_{*}]& \to &{T_{1}}_{*}& \to &0 \cr
 &&   {\downarrow} && {\downarrow} && {\downarrow}  &&  \cr
 0& \to& \mu_{2}(-1) &\to& {\rm Coker }[(S_{2})_{*} \to (P_{2})_{*}] &\to &{T_{2}}_{*}& \to& 0
.}$$
Comme les fl\`eches verticales de gauche et de droite sont injectives, on conclut que la fl\`eche m\'ediane, c'est-\`a-dire $\pi_{1}(G_{1}) \to \pi_{1}(G_{2})$, est injective.

Soit $p$ l'homomorphisme: $G_{2} \to G_{3}$. Soit $n>0$ un entier. On v\'erifie que l'on a des isomorphismes naturels de faisceaux \'etales sur ${\overline G}_{3}$ :
$\Z/n  \buildrel \simeq \over \rightarrow p_{*} \Z/n$ et $H^1({\overline G}_{1},\Z/n)  \buildrel \simeq \over \leftarrow R^1p_{*} \Z/n$, o\`u les
termes de gauche sont ``constants'', i.e. proviennent de faisceaux sur $k$.
La suite exacte
$$0 \to H^1({\overline G}_{3},p_{*}\Z/n) \to  H^1({\overline G}_{2},\Z/n) \to  H^0({\overline G}_{3},R^1p_{*} \Z/n)$$
des termes de bas degr\'e de la  suite spectrale de Leray pour ${\overline p} $
s'\'ecrit donc
$$0 \to H^1({\overline G}_{3}, \Z/n) \to  H^1({\overline G}_{2},\Z/n) \to  H^1({\overline G}_{1},\Z/n).$$

Pour \'etablir qu'un complexe (born\'e)  $A^{\bullet}$ de groupes ab\'eliens de type fini est exact, il suffit
de voir que le complexe ${\rm Hom}_{\Z}(A^{\bullet},\Q/\Z)$  est exact, et pour cela il suffit de 
montrer que pour tout entier positif $n$, le complexe 
${\rm Hom}_{\Z}(A^{\bullet},\Z/n)$ est exact.
Combinant ceci avec la proposition 6.7, on voit que le complexe
$$\pi_{1}(G_{1}) \to  \pi_{1}(G_{2}) \to   \pi_{1}(G_{3})   \to 0$$
est exact. Comme on a \'etabli ci-dessus l'injectivit\'e de $\pi_{1}(G_{1}) \to  \pi_{1}(G_{2})$,
ceci ach\`eve la d\'emonstration. \cqfd
\medskip
{\it Remarque 6.8.1} Pour le groupe $\pi_{1}^{Bor}(G)$, la proposition ci-dessus est \'enonc\'ee 
 dans [Bo98] (Lemma 1.5).

\bigskip

{\bf \S 7.  Groupe de Brauer d'une compactification lisse}

\bigskip

 Etant donn\'es $k$ un corps et $X$ une $k$-vari\'et\'e alg\'ebrique, on note
 $${\rm Br}_{1}(X) = {\rm Ker} [{\rm Br}(X) \to {\rm Br}({\overline X})].$$

{\bf Th\'eor\`eme 7.1} {\it Soient $k$   et $G$ un
$k$-groupe lin\'eaire connexe, suppos\'e r\'eductif si $p={\rm car}(k)>0$.
Soit $X$ une $k$-compactification lisse de $G$ et soit
$1 \to S \to H \to G \to 1$
une r\'esolution flasque du $k$-groupe $G$.
On a alors
$${\rm Br}_{1}(X)/{\rm Br}(k) \buildrel \simeq \over \rightarrow H^1(k,S^*).$$
Le groupe $H^1(k,S^*)$ s'identifie  \`a un sous-groupe de ${\rm Br}(X)/{\rm Br}(k) $,
le quotient \'etant un groupe de torsion $p$-primaire.
 }

\medskip

{\it D\'emonstration} 
Le premier isomorphisme r\'esulte de la proposition 5.2 et de la formule ${\rm Br}_{1}(X)/{\rm
Br}(k)\buildrel \simeq \over \rightarrow H^1(k,{\rm Pic}({\overline X}))$, valable pour toute $k$-vari\'et\'e projective lisse
g\'eom\'etriquement int\`egre poss\'edant un $k$-point. 
La $\k$-vari\'et\'e ${\overline G}$ est $\k$-birationnelle \`a un espace affine.
La $\k$-vari\'et\'e projective, lisse, connexe ${\overline X}$ est donc 
$\k$-birationnelle \`a un espace projectif ${\bf P}^n_{\k}$. On sait que la partie premi\`ere
\`a la caract\'eristique du groupe de Brauer est un invariant birationnel des vari\'et\'es projectives, lisses, connexes ([Gr68])
et que le groupe de Brauer de ${\bf P}^n_{\k}$ est trivial ([Gr68]).
Ainsi ${\rm Br}({\overline X})$ est nul si  ${\rm car}(k)=0$ et il est de torsion 
$p$-primaire si  $p={\rm car}(k)>0$. Ceci
ach\`eve d'\'etablir l'\'enonc\'e.  \cqfd

\bigskip

Soit  $\g$ le groupe de Galois de ${\overline k}$ sur $k$.
 A tout $\g$-module continu $M$ et tout entier naturel $i$ on associe le groupe
$$\X^{i}_{\omega}(k,M) = {\rm Ker} [H^{i}(\g,M) \to  \prod_{\h} H^{i}(\h,M)],$$ o\`u $h$ parcourt les
sous-groupes ferm\'es procycliques de $\g$.
Avec la d\'efinition de Borovoi du groupe fondamental alg\'ebrique, l'\'enonc\'e suivant
a \'et\'e \'etabli par Borovoi et Kunyavski\u{\i} [BoKu00].

\medskip

{\bf Th\'eor\`eme 7.2} {\it  
Soient $k$ un corps et $G$ un
$k$-groupe lin\'eaire connexe, suppos\'e r\'eductif si $p={\rm car}(k)>0$.
Soit $X$ une $k$-compactification lisse de $G$.
On a 
$${\rm Br}_{1}(X)/{\rm Br}(k) \simeq {\cyr X}^1_{\omega}(k,\Hom_{\Z}(\pi_{1}(G),\Q/\Z )).$$
Le groupe ${\cyr X}^1_{\omega}(k,\Hom_{\Z}(\pi_{1}(G),\Q/\Z  ))$ s'identifie  \`a un sous-groupe de ${\rm Br}(X)/{\rm Br}(k) $,
le quotient \'etant un groupe de torsion $p$-primaire.
 }

\medskip

{\it D\'emonstration} Soit $1 \to S \to  H \to G \to 1$ une r\'esolution flasque de $G$. 
Soit $P$ le $k$-tore quasi-trivial $H^{tor}$. D'apr\`es la d\'efinition de $\pi_{1}(G)$, on a la suite exacte de modules galoisiens
$$ 0 \to S_{*} \to P_{*} \to \pi_{1}(G) \to 0.$$ 
Appliquant le foncteur $\Hom_{\Z}(\bullet,\Q/\Z)$, on en d\'eduit la suite exacte de modules galoisiens
$$0 \to \Hom_{\Z}(\pi_{1}(G), \Q/\Z  ) \to P^*\otimes_{\Z} (\Q/\Z )  \to S^*\otimes_{\Z} (\Q/\Z )  \to 0$$
(le groupe ab\'elien $\Q/\Z$ est injectif).
Ceci donne naissance \`a la suite exacte de groupes ab\'eliens
$$(P^*\otimes_{\Z}(\Q/\Z ))^\g   \to (S^*\otimes_{\Z}(\Q/\Z) )^\g \to H^1(\g, \Hom_{\Z}(\pi_{1}(G),\Q/\Z ))
\to H^1(\g,P^*\otimes_{\Z}(\Q/\Z ))  .$$
Tensorisant la suite exacte
$$0 \to \Z \to \Q \to \Q/\Z \to 0$$
par $P^*$ et par $S^*$, et prenant la suite de cohomologie associ\'ee, on obtient le diagramme 
commutatif de suites exactes
$$\diagram {(P^*\otimes_{\Z}\Q )^\g & \to &(P^*\otimes_{\Z}(\Q/\Z ) )^\g &\to &H^1(\g,P^*) &\to& 0\cr
\downarrow&& \downarrow&& \downarrow&&\cr
(S^*\otimes_{\Z}\Q )^\g &\to &(S^*\otimes_{\Z}(\Q/\Z ))^\g &\to &H^1(\g,S^*)& \to& 0.
}$$
La fl\`eche $P^* \to S^*$ a un conoyau fini (Proposition 3.3).
Ainsi la fl\`eche induite $P^*\otimes_{\Z}\Q \to S^*\otimes_{\Z}\Q$ est surjective, et il en est donc aussi ainsi de $(P^*\otimes_{\Z}\Q)^\g \to (S^*\otimes_{\Z}\Q)^\g$. Par ailleurs $H^1(\g,P^*)=0$.
Ceci implique que le conoyau de la fl\`eche $(P^*\otimes_{\Z}(\Q/\Z ))^\g   \to (S^*\otimes_{\Z}(\Q/\Z ))^\g$
s'identifie \`a $H^1(\g,S^*)$, on  a donc la suite exacte
$$0 \to H^1(\g,S^*)  \to H^1(\g, \Hom_{\Z}(\pi_{1}(G),\Q/\Z  ))
\to H^1(\g,P^*\otimes_{\Z}(\Q/\Z ))  .$$
On a des  suites exactes compatibles en rempla\c cant $\g$ par tout sous-groupe ferm\'e,
en particulier tout sous-groupe ferm\'e procyclique $\frak{h}$. On a donc une suite exacte induite
$$0 \to {\cyr X}^1_{\omega}(\g,S^*)  \to {\cyr X}^1_{\omega}(\g, \Hom_{\Z}(\pi_{1}(G),\Q/\Z  ))
\to {\cyr X}^1_{\omega}(\g,P^*\otimes_{\Z}(\Q/\Z ))  .$$
Comme $S^*$ est flasque, on a $H^1({\frak{h}},S^*)=0$ pour tout sous-groupe ferm\'e procyclique,
donc ${\cyr X}^1_{\omega}(\g,S^*) =H^1(\g,S^*)$. Par ailleurs, comme 
$P^*$ est un module de permutation, on a ${\cyr X}^1_{\omega}(\g,P^*\otimes_{\Z}(\Q/\Z ))=0$
(ceci est clair pour $P=\Z$, et se ram\`ene \`a ce cas par le lemme de Shapiro). On a donc bien
$$H^1(\g,S^*) \buildrel \simeq \over \rightarrow  {\cyr X}^1_{\omega}(\g, \Hom_{\Z}(\pi_{1}(G),\Q/\Z  )),$$
ce qui  compte tenu du th\'eor\`eme 7.1
 \'etablit le th\'eor\`eme.\cqfd

\bigskip

{\it Remarques 7.2.1}
 
(1) A la torsion $p$-primaire pr\`es, on calcule donc le groupe de Brauer d'une compactification
lisse de $G$, et donc du groupe de Brauer non ramifi\'e du corps des 
fonctions $k(G)$ 
 sans avoir \`a construire une compactification lisse explicite.
Il serait int\'eressant de donner une d\'emonstration directe du fait que
le groupe de Brauer non ramifi\'e de $k(G)$ 
est donn\'e par les formules
des th\'eor\`emes  7.1 et 7.2.

(2) Des cas particuliers des th\'eor\`emes 7.1 et 7.2 avaient \'et\'e \'etablis par Sansuc [Sa81] pour $k$ un corps de nombres.
Pour $G=T$ un $k$-tore, le th\'eor\`eme 7.2 se lit 
${\rm Br}(X)/{\rm Br}(k) \simeq {\cyr X}^2_{\omega}(k,T^*)$ (\`a la torsion $p$-primaire pr\`es). Cet \'enonc\'e avait \'et\'e \'etabli dans
[CTSa87a]. 
Pour  un $k$-groupe semi-simple $G$ de groupe fondamental $\mu = \Ker [
G^{sc} \to G ],$
 le th\'eor\`eme 7.2  se lit  
${\rm Br}(X)/{\rm Br}(k) \simeq {\cyr X}^1_{\omega}(k,\mu^*)$ (\`a la torsion $p$-primaire pr\`es).
Cet \'enonc\'e est  \'etabli dans
[CTKu98] sur un corps $k$ de caract\'eristique nulle quelconque, par un argument assez d\'etourn\'e (recours aux corps finis).

(3) Pour la comparaison entre  les r\'esultats ci-dessus et les r\'esultats  de 
Borovoi et Kunyavski\u{\i} [BoKu00, BoKu04], voir les remarques 6.2.1 et
l'appendice A.

\bigskip

{\bf \S 8. Cohomologie galoisienne des groupes lin\'eaires connexes }

\bigskip

Dans ce paragraphe et le suivant, \'etant donn\'e un $k$-groupe lin\'eaire lisse connexe,
nous cherchons \`a d\'eterminer le groupe $G(k)/R$ et l'ensemble $H^1(k,G)$.
Pour les notions de base sur la $R$-\'equivalence sur un groupe alg\'ebrique, le lecteur consultera [CTSa77] et [Gi97]).
En caract\'eristique z\'ero, la suite exacte
$$ 1 \to G^u \to G \to G^{\red} \to 1$$
induit une bijection $G(k)/R \buildrel \simeq \over \rightarrow G^{\red}(k)/R$ et une bijection
$H^1(k,G) \buildrel \simeq \over \rightarrow H^1(k,G^{\red})$. Le premier \'enonc\'e r\'esulte simplement du fait
que $G^u$ est un espace affine et que de la suite ci-dessus on tire un isomorphisme
de $k$-vari\'et\'es $G^u \times G^{\red} \simeq G$. Pour le second \'enonc\'e, voir
[Sa81], Lemme 1.13 p. 20.

On se limitera donc dans la suite \`a consid\'erer le cas o\`u $G$ est un $k$-groupe r\'eductif connexe.

\bigskip

Soit $G$ un $k$-groupe r\'eductif connexe.
Soit
$$1  \to S \to H  \to G  \to 1$$
une r\'esolution flasque de $G$, 
avec $H^{ss}$ simplement connexe et 
$H^{tor}$ quasi-trivial :
$$ 1 \to H^{ss} \to H \to H^{tor} \to 1.$$

Comme $S$ est central dans $H$,
cette suite exacte donne naissance \`a une suite exacte de groupes  
$$ 1 \to S(k) \to H(k) \to G(k) \to H^1(k,S) $$
qui s'ins\`ere dans une suite exacte
d'ensembles point\'es de cohomologie galoisienne
$$ 1 \to S(k) \to H(k) \to G(k) \to H^1(k,S) \to H^1(k,H) \to H^1(k,G) \to H^2(k,S) \hskip2mm (8.1)$$
On peut  pr\'eciser les fibres 
des applications  
(Serre, [SeCG], chap. I, \S 5). En particulier toute fibre de l'application 
$H^1(k,G) \to H^2(k,S)$ est soit vide soit un quotient d'un ensemble
$H^1(k,{}_{c}H)$ par une action de $H^1(k,S)$, le groupe ${}_{c}H$
\'etant un groupe obtenu par torsion de $H$, et dont on v\'erifie qu'il est
comme $H$ quasi-trivial.

On a en outre
la suite exacte de groupes
$$ 1 \to H^{ss}(k) \to H(k) \to H^{tor}(k)$$
qui s'ins\`ere dans la suite exacte d'ensembles point\'es
$$ 1 \to H^{ss}(k) \to H(k) \to H^{tor}(k) \to H^1(k,H^{ss}) \to H^1(k,H) \to H^1(k,H^{tor}) $$
qui compte tenu du th\'eor\`eme 90 de Hilbert donne la suite exacte d'ensembles point\'es
$$ 1 \to H^{ss}(k) \to H(k) \to H^{tor}(k) \to H^1(k,H^{ss}) \to H^1(k,H) \to 1 \hskip2mm (8.2)$$
\bigskip

{\bf Th\'eor\`eme 8.1} {\it Soit $k$ un corps.

(i) Une r\'esolution flasque $1 \to S \to H \to G \to 1$ 
d'un $k$-groupe r\'eductif connexe $G$ induit une  suite exacte de groupes 
$$ H(k)/R \to G(k)/R \to \Ker [H^1(k,S) \to H^1(k,H)] \to 1.$$

(ii) Le quotient de $G(k)/R$ par l'image de $H(k)$,
c'est-\`a-dire l'image de $G(k)$ dans $H^1(k,S)$,
est un quotient ab\'elien de $G(k)/R$ qui ne d\'epend pas
de la r\'esolution flasque choisie.

(iii) Si $k$ est un corps de type fini sur le corps premier,
ou si $k$ est un corps de type fini sur un corps alg\'ebriquement clos
de caract\'eristique z\'ero, alors ce quotient est fini.}

\medskip

{\it D\'emonstration}
Consid\'erons la suite exacte (8.1). Comme le $k$-tore 
$S$ est flasque, l'application $G(k) \to H^1(k,S)$
passe au quotient par la $R$-\'equivalence ([CTSa77], Prop. 12).
La $R$-\'equivalence \'etant fonctorielle, l'homomorphisme
$H(k) \to G(k)$ induit un homomorphisme
$H(k)/R \to G(k)/R$. Ceci donne la suite exacte du (i).

Soient $1 \to S \to H \to G \to 1$ et $1 \to S_1 \to H_1 \to G \to 1$
deux r\'esolutions flasques de $G$. Les projections du produit fibr\'e
$E=H\times_{G} H_1$ sur $H$ et sur $H_1$ admettent des sections (\S 3).
Ainsi le quotient  de $G(k)$ par l'image de $H(k)$ et le quotient de
$G(k)$ par l'image de $H_1(k)$ co\"{\i}ncident tous deux avec
le quotient de  $G(k)$ par l'image de $E(k)$. Ceci \'etablit le point  (ii).

La finitude \'enonc\'ee en (iii) r\'esulte de la finitude
de $H^1(k,S)$ pour $S$ un $k$-tore flasque et $k$ comme dans l'\'enonc\'e
 ([CTSa77];  [CTGiPa04, Thm. 3.4]). \qed

\medskip

{\it Remarque 8.1.1} Pour
$k$ de type fini sur le corps premier et $G/k$ r\'eductif, c'est un probl\`eme ouvert de savoir si 
le groupe $G(k)/R$ est fini. La proposition ci-dessus ram\`ene
le probl\`eme au cas o\`u $G$ est un groupe r\'eductif 
quasi-trivial. La suite exacte (8.2) ne semble pas permettre de r\'eduire le probl\`eme
au seul cas des groupes semi-simples simplement connexes.

Pour  $k$ corps de nombres ou $k$ corps local
et $G$ semi-simple simplement connexe ou adjoint,
le
quotient $G(F)/R$ peut \^etre non trivial pour
$F/k$ extension (transcendante) convenable. 
Pour le cas simplement connexe, voir 
Merkur'ev [Me93]). 
Pour le
cas  adjoint, voir Merkur'ev [Me96]. A noter que dans ces deux cas,
 on peut montrer que pour toute r\'esolution flasque $1 \to S \to H \to G \to 1$
et toute extension de corps $F/k$, on a $H^1(F,S)=0$ (voir [CTGiPa04], Cor. 4.11).

\bigskip

{\bf Proposition 8.2} {\it  Soient $k$ un corps, $X$ un  $k$-sch\'ema  et
$G$ un $k$-groupe lin\'eaire  connexe, suppos\'e r\'eductif
 si ${\rm car}($k$)>0$. Soit
 $$1 \to S \to H \to G \to 1$$
 une r\'esolution flasque de $G$, et soit $P=H^{tor}$.
 Une telle r\'esolution induit une application naturelle
 $$ H^1(X,G) \to \Ker [H^2(X,S) \to H^2(X,P)].$$
 Il existe un homomorphisme naturel
 $$ \Ker [H^2(X,S) \to H^2(X,P)] \to \Hom(\Pic(G), \Br(X)).$$
 L'application naturelle compos\'ee 
 $$H^1(X,G) \to   \Hom(\Pic(G), \Br(X))$$
 ne d\'epend pas du choix de la r\'esolution flasque.
}

\medskip

{\it D\'emonstration} 
On a le diagramme commutatif de suites exactes de $k$-groupes alg\'ebriques
$$\diagram {1& \to &S& \to& H& \to &G& \to& 1 \cr 
&   & \downarrow && \downarrow &&\downarrow &&\cr
1& \to &   M   & \to& P& \to &T& \to& 1, }
$$
o\`u $T=G^{tor}$ et $P=H^{tor}$.
Soit $X$ un $k$-sch\'ema. De ce diagramme on d\'eduit (en passant par $H^1(X,T)$) que la fl\`eche naturelle
$H^1(X,G) \to H^2(X,S)$ 
a son image dans le noyau de $H^2(X,S) \to H^2(X,P)$.

D'apr\`es la proposition 3.3,
la r\'esolution flasque   induit une suite exacte
$$   (P^*)^{\g} \to (S^*)^{\g} \to \Pic(G) \to 0,$$
o\`u la  fl\`eche $(S^*)^{\g} \to \Pic(G)$ 
admet la description simple suivante :
la donn\'ee d'un caract\`ere $S \to \G_{m,k}$ associe \`a l'extension
donn\'ee par la r\'esolution flasque une extension de $G $ par  $\G_{m,k},$ 
donc un \'el\'ement de $\Ext^1_{k-gp}(G, \G_{m,k})$, groupe qui s'envoie
(de fa\c con isomorphique, 
voir le corollaire 5.7) vers $\Pic(G)$ par l'application d'oubli \'evidente.

Le diagramme d'accouplements
$$\diagram  {H^2(X, S) & \times & (S^*)^{\g} & \to & \Br(X) \cr
              \downarrow &                 &   \uparrow  &   &      \downarrow{=}{}  
                \cr
H^2(X, P) & \times & (P^*)^{\g} & \to & \Br(X)}$$
 \'etant clairement commutatif,  on en d\'eduit un accouplement entre
 le noyau de $H^2(X,S) \to H^2(X,P)$ et le conoyau de
 $(P^*)^{\g} \to (S^*)^{\g}$, conoyau dont on vient de rappeler qu'il s'identifie \`a $\Pic(G)$.
On a donc une application induite
$$ H^1(X,G)  \times \Pic(G) \to \Br(X).$$
Cette application est induite
 par l'accouplement \'evident
$$H^1(X,G) \times \Ext^1_{k-gp}(G, \G_{m,k}) \to \Br(X).$$
Elle est donc ind\'ependante de la r\'esolution flasque de $G$.
\cqfd

\medskip

{\it Remarque} Sur une $k$-vari\'et\'e lisse $X$, l'existence d'une application naturelle
$$H^1(X,G) \to \Hom(\Pic(G),\Br(X))$$ est d\'ej\`a \'etablie par Sansuc ([Sa81], Prop. 6.10).
Il montre en effet que  la donn\'ee d'un torseur $Y$
sur une telle $k$-vari\'et\'e 
$X$ sous le $k$-groupe $G$ 
d\'efinit   un homomorphisme $\Pic(G) \to \Br(X)$.
Comme la suite est  ``naturelle'', deux $G$-torseurs isomorphes donnent naissance
\`a la m\^eme application $\Pic(G) \to \Br(X)$.
Ceci donne bien une application
$H^1(X,G) \to \Hom(\Pic(G),\Br(X)).$

\bigskip

{\bf Proposition 8.3} {\it Soient $k$ un corps 
et  $G$ un $k$-groupe r\'eductif connexe.
Soit
$$1  \to S \to H  \to G  \to 1$$
une r\'esolution flasque de $G$. Soit $P=H^{tor}$.
Une telle r\'esolution induit une application
$$ H^1(k,G) \to \Ker [H^2(k,S) \to H^2(k,P)]$$
dont toute fibre est vide ou est un quotient
d'un ensemble $H^1(k,{}_{c}H^{ss})$ pour
un $k$-groupe semi-simple simplement connexe ${}_{c}H^{ss}$
convenable. 

Par ailleurs une telle r\'esolution induit
un complexe d'ensembles point\'es
$$ H^1(k,H^{ss}) \to H^1(k,G) \to \Hom(\Pic(G),\Br(k)).$$
}

{\it D\'emonstration} 
L'\'enonc\'e  r\'esulte des suites (8.1) 
et (8.2) et de la proposition 8.2.
\cqfd

\bigskip

Le th\'eor\`eme suivant repose sur un grand nombre de r\'esultats ant\'erieurs.
En particulier, il utilise   et  \'etend  des r\'esultats de [Gi97],
[CTGiPa04] et [BoKu04].

\medskip

{\bf Th\'eor\`eme 8.4} {\it Soit $k$ un corps de caract\'eristique
nulle, de dimension cohomologique $\leq 2$, tel que sur toute extension
finie $K$ de $k$, indice et exposant des $K$-alg\`ebres simples centrales  
co\"{\i}ncident. On suppose en outre que la dimension cohomologique
de l'extension ab\'elienne maximale de $k$ est $\leq 1$
(cette hypoth\`ese n'est requise qu'en pr\'esence de facteurs de type $E_{8}$.)
Soit $1 \to S \to H \to G \to 1$ une r\'esolution flasque 
d'un $k$-groupe r\'eductif connexe $G$.

(i) Cette r\'esolution induit 
 un isomorphisme
$G(k)/R \buildrel \simeq \over \rightarrow H^1(k,S)$.

(ii)  Cette r\'esolution induit une bijection
$H^1(k,G) \buildrel \simeq \over \rightarrow \Ker [H^2(k,S) \to H^2(k,P)]$.}

\medskip

{\it D\'emonstration} 
Pour $k$ comme dans l'\'enonc\'e et $H$ un 
$k$-groupe   semi-simple
simplement connexe, on a $H^1(k,H)=1$.
On renvoie \`a
 [CTGiPa04], Thm. 1.2, pour l'historique de la d\'emonstration de ce r\'esultat.
Pour $H$ un $k$-groupe quasi-trivial, l'\'enonc\'e $H^1(k,H)=1$
r\'esulte alors de  la suite exacte (8.2).

Pour $H$ un $k$-groupe semi-simple simplement connexe et $k$ comme ci-dessus,
$H(k)/R=1$.
On trouvera dans [CTGiPa04], Thm. 4.5, la d\'emonstration de ce r\'esultat
ainsi qu'un historique.
C'est un th\'eor\`eme d\'elicat de Philippe Gille ([Gi01], appendice de [BoKu04])
qu'une suite exacte de $k$-groupes r\'eductifs
$$ 1 \to H_1 \to H \to T \to 1$$
avec $H_1$ semi-simple simplement connexe et $T$ un $k$-tore,
et $k$ comme dans l'\'enonc\'e ci-dessus, induit une suite exacte
$$ H_1(k)/R \to H(k)/R \to T(k)/R \to 1.$$
 Si le tore $T$ est quasi-trivial,
on a $T(k)/R=1$. Ainsi pour $H$ quasi-trivial sur $k$  comme ci-dessus, on a $H(k)/R=1$. 
L'\'enonc\'e (i) r\'esulte alors du th\'eor\`eme 8.1.

L'\'enonc\'e (ii) est une variante d'un th\'eor\`eme
de Borovoi et Kunyavski\u{\i} ([BoKu04],  th\'eor\`eme 6.7).
L'injectivit\'e de l'application
$H^1(k,G) \to \Ker [H^2(k,S) \to H^2(k,P)]$
 r\'esulte  de la proposition 8.3
et de la trivialit\'e des ensembles $H^1(k,H_1)$ pour tout 
$k$-groupe semi-simple simplement connexe $H_1$.
Indiquons les points-cl\'es utilis\'es pour \'etablir la surjectivit\'e. 
Pour $k$ comme dans le th\'eor\`eme et $G$ un $k$-groupe
semi-simple  connexe, le th\'eor\`eme 2.1
de [CTGiPa04] \'etablit que la fl\`eche
$H^1(k,G) \to H^2(k,\mu)$
d\'eduite de la suite naturelle
$1 \to \mu \to G^{sc} \to G \to 1$
est une bijection.
Une variante d'un argument de Borovoi  ([Bo93], Thm. 5.5)
permet alors d'en d\'eduire que pour $k$ comme ci-dessus
et $L=({\overline H},\kappa)$ un $k$-lien avec ${\overline H}$
un ${\overline k}$-groupe semi-simple, tout  \'el\'ement de
l'ensemble  de cohomologie $H^2(k,L)$ est neutre ([CTGiPa04],
Prop. 5.4, noter que pour un tel  ${\overline H}$ on a ${\overline H}^{tor}=1$). 
Pour terminer la d\'emonstration, on suit alors la d\'emonstration
du \S 6 de [BoKu04], qui passe par la cohomologie de modules crois\'es,
et on utilise l'appendice A ci-dessous, qui permet de comparer
l'application $H^1(k,G) \to  \Ker [H^2(k,S) \to H^2(k,P)]$
avec la fl\`eche $H^1(k,G) \to H^1_{{ab}}(k,G)$ de Borovoi.
\qed

\medskip

{\it Remarques 8.4.1}

(1) L'\'enonc\'e (i) est d\'ej\`a \'etabli dans [CT05]  (Th\'eor\`eme 3 et remarque subs\'equente).

(2) De l'\'enonc\'e (ii) et du diagramme fondamental  (\S 3)
on peut d\'eduire le th\'eor\`eme 2.1  de [CTGiPa04], utilis\'e dans la
d\'emonstration ci-dessus.
Mais partant de ce diagramme et du th\'eor\`eme 2.1  de [CTGiPa04],
je n'ai  pas r\'eussi \`a donner une d\'emonstration
de l'\'enonc\'e (ii) \'evitant le recours \`a la cohomologie
des modules crois\'es, comme introduite par Borovoi dans [Bo 98]
et appliqu\'ee dans [BoKu04].

(3) Rappelons que les hypoth\`eses du th\'eor\`eme 8.4 sont satisfaites
dans chacun des cas suivants (voir [CTGiPa04] et [BoKu04]) :

-- Corps  $p$-adiques;  

-- Corps de nombres
totalement imaginaire ; 

-- Corps de type (ll) : $k$ est le corps des fractions d'un
anneau local int\`egre strictement hens\'elien de dimension 2
et de corps r\'esiduel (alg\'ebriquement clos) de caract\'eristique 0.
(Exemple : une extension finie du corps ${\bf C}((x,y))$
quotient de l'anneau des s\'eries formelles \`a deux variables.)

-- Corps de type (gl) : Un corps $k$ de fonctions de deux variables sur un
corps alg\'ebriquement clos de caract\'eristique 0
est de dimension cohomologique 2 ; 
l'hypoth\`ese indice/exposant est satisfaite
(th\'eor\`eme   de de Jong [dJ04]).
 On ne conna\^{i}t pas
dans ce dernier cas la dimension  cohomologique
de l'extension ab\'elienne maximale de $k$, mais l'hypoth\`ese
qu'elle est au plus $1$ n'est utilis\'ee dans le th\'eor\`eme
ci-dessus que lorsque $G$ poss\`ede un facteur de type $E_8$.

(4) 
Dans chacun des cas mentionn\'es dans la remarque pr\'ec\'edente, le groupe $H^1(k,S)$
est fini pour $S$ un $k$-tore flasque (voir [CTGiPa04], Thm. 3.4), donc dans chacun de ces cas le groupe  quotient $G(k)/R$ est  fini.

(5)
Une fois le th\'eor\`eme 8.4 \'etabli, il est facile d'en d\'eduire de fa\c con un peu plus fonctorielle,
sur les corps g\'eom\'etriques mentionn\'es ci-dessus, la g\'en\'eralisation
des r\'esultats 
de [CTGiPa04]  obtenue dans [BoKu04] (\S  5, contr\^ole de l'approximation
faible sur un corps de type (ll) ou (gl) ; \S 7, contr\^ole du d\'efaut 
du principe de Hasse sur un corps de type (ll)).

\bigskip

{\bf Proposition 8.5} {\it Soit $k$ un corps de caract\'eristique nulle,
de dimension cohomologique $\leq 2$, tel que sur toute extension finie $K$ de $k$,
indice et exposant co\"{\i}ncident pour les $K$-alg\`ebres simples centrales.
Soit $G$ un $k$-groupe r\'eductif
connexe. Si $G$ contient des facteurs de type
$E_{8}$, supposons que la dimension cohomologique de l'extension ab\'elienne
maximale de $k$ est au plus~1. Si un espace homog\`ene principal  $X$ de $G$
poss\`ede un z\'ero-cycle de degr\'e 1, alors il poss\`ede un point $k$-rationnel.}

\medskip

{\it D\'emonstration} Soit 
$$1 \to S \to H \to G \to 1$$
une r\'esolution flasque de $G$.
Les hypoth\`eses impliquent $H^1(k,H^{ss})=1$.
Soit $[X] \in H^1(k,G)$ la classe de l'espace homog\`ene principal $X$.
Si $X$ poss\`ede un z\'ero-cycle de degr\'e 1, il existe des extensions finies $k_{i}/k$,
de degr\'es premiers entre eux dans leur ensemble, telles que
$[X]$ ait une image triviale dans $H^1(k_{i},G)$. On conclut que l'image
de $[X]$ dans $H^2(k,S)$ devient triviale dans les groupes $H^2(k_{i},S)$.
Comme $S$ est commutatif, ceci implique (argument de norme) que l'image de $[X]$ dans $H^2(k,S)$
est nulle. D'apr\`es la proposition 8.3, $[X]$ est donc dans l'image de
$H^1(k,H^{ss})$ dans $H^1(k,G)$. Mais les hypoth\`eses sur $k$ assurent 
$H^1(k,H^{ss})=1$. Ainsi $[X]=1 \in H^1(k,G)$. \qed

\bigskip

{\it Remarque 8.5.1}
L'\'enonc\'e pr\'ec\'edent s'applique aux corps de nombres totalement imaginaires.
 Sur un corps de nombres $k$ qui admet un plongement r\'eel,
le m\^eme \'enonc\'e vaut (sans restriction sur les facteurs de type $E_{8}$).
On utilise le fait que $S(k)$ est dense dans le produit
des $S(k_v)$ pour $v$ parcourant les places r\'eelles, et que
l'application  $H^1(k,H) \to \prod_{v \hskip1mm r{\acute e}elle}H^1(k_v,H)$ est
une bijection.

\bigskip

{\bf \S 9.  Cohomologie galoisienne des groupes r\'eductifs connexes sur un corps de nombres }

\bigskip

Le th\'eor\`eme 8.4 permet de retrouver un grand nombre de r\'esultats classiques
sur les corps $p$-adiques et les corps de nombres totalement imaginaires. 
Dans ce paragraphe nous discutons aussi le cas d'un corps de nombres arbitraire.
Pour les raisons donn\'ees au d\'ebut du \S 8, on peut se limiter au cas des groupes
r\'eductifs connexes.

\bigskip

{\bf Th\'eor\`eme 9.1} {\it Soient $k$ un corps local de caract\'eristique nulle,
$G$ un $k$-groupe r\'eductif connexe
et $1 \to S \to H \to G \to 1$
une r\'esolution flasque de $G$.
Une telle r\'esolution induit :

(i) un isomorphisme de groupes ab\'eliens $G(k)/R \buildrel \simeq \over \rightarrow H^1(k,S)$;

(ii) si $k$ est local non archim\'edien une bijection
$$H^1(k,G) \buildrel \simeq \over \rightarrow \Ker [H^2(k,S) \to H^2(k,P)]$$
et des bijections 
$$H^1(k,G) \buildrel \simeq \over \rightarrow  \Hom(\Pic(G),\Q/\Z) \simeq   (\pi_{1}(G)_{\frak g})_{tors}   ;$$

(iii) si $k$ est r\'eel, une surjection
$$H^1(k,G) \to \Ker [H^2(k,S) \to H^2(k,P)]$$
et une application
$$H^1(k,G) \to \Hom(\Pic(G),\Q/\Z) \simeq   (\pi_{1}(G)_{\frak g})_{tors}  .$$

Dans (ii) et (iii) les applications $H^1(k,G) \to \Hom(\Pic(G),\Q/\Z)   $
sont ind\'ependantes du choix de la r\'esolution flasque de $G$.
}

\medskip

{\it D\'emonstration} L'\'enonc\'e (i) pour $k$ non archim\'edien est un  cas particulier
du th\'eor\`eme 8.4~(i). Si $k$ est r\'eel,
tout $k$-tore flasque est   quasi-trivial, et pout tout $k$-tore $T$ on a $T(k)/R=1$.
Ceci implique $G(k)/R=1$ et $H^1(k,S)=1$.

Pour $k$ non archim\'edien, la bijection $H^1(k,G) \buildrel \simeq \over \rightarrow \Ker [H^2(k,S) \to H^2(k,P)]$
est un  cas particulier du th\'eor\`eme 8.4 (ii). Combinant la proposition 3.3, les th\'eor\`emes de dualit\'e du corps de classes local pour la cohomologie des tores ([SeCG], Chap. II, \S 5.8)
et la finitude du groupe $\Pic(G)$,
on obtient l'isomorphisme
$$ \Ker [H^2(k,S) \to H^2(k,P)]  \simeq \Hom(\Pic(G),\Q/\Z).$$
L'isomorphisme $ \Hom(\Pic(G),\Q/\Z) \simeq   (\pi_{1}(G)_{\frak g})_{tors} $
est un fait g\'en\'eral, valable sur tout corps, il a fait l'objet de la proposition 6.3.

Supposons $k$ r\'eel. En utilisant le diagramme fondamental  (\S 3) et
un $k$-tore $T_{G} \subset G$ tel que $T_{G^{ss}}$ soit un ``tore fondamental''
(voir [Bo98] \S 5.3), on \'etablit que l'application
$H^1(k,G) \to \Ker [H^2(k,S) \to H^2(k,P)]$ est surjective.
 Par la dualit\'e locale sur les r\'eels ([MiADT], Chap. I, Thm. 2.13),
 le groupe fini $\Ker [H^2(k,S) \to H^2(k,P)]$ est dual du groupe
 $\Coker [{\hat H}^0(\frak{g},P^*) \to {\hat H}^0(\frak{g},P^*)] $, o\`u ici $\frak{g}=\Z/2$.
 De la proposition 3.3 on d\'eduit que le groupe $\Ker [H^2(k,S) \to H^2(k,P)]$
 est un sous-groupe de  $ \Hom(\Pic(G),\Z/2) \subset \Hom(\Pic(G),\Q/\Z).$

 La derni\`ere assertion de l'\'enonc\'e r\'esulte de la proposition 8.2 et de sa d\'emonstration.
 Les fl\`eches $H^1(k,G) \to \Hom(\Pic(G),\Q/\Z)$ sont en effet induites par
 l'accouplement naturel 
 $H^1(k,G) \times \Pic(G) \to \Br(k)$
 combin\'e avec le plongement $\Br(k) \hookrightarrow \Q/\Z$ donn\'e par
 la th\'eorie du corps de classes local.
\cqfd

\medskip

{\it Remarque 9.1.1} L'\'enonc\'e 9.1 (ii) est une variante d'un r\'esultat de Kottwitz
[Ko84, Prop. 6.4] [Ko86, Thm. 1.2] et de Borovoi [Bo98, Cor. 5.5 (i)]. L'\'enonc\'e 9.1 (iii) est une variante
d'un r\'esultat de Kottwitz [Ko86, Thm. 1.2] et de Borovoi [Bo98, Cor. 5.5 (ii)].
 
\bigskip

{\bf Proposition 9.2} {\it Soient $k$ un corps de nombres et
 $H$ un $k$-groupe r\'eductif quasi-trivial. Alors

(i) On a
$H^1(k_v,H)=1$ pour $v$ place non r\'eelle;

(ii) L'approximation faible vaut pour $H$;

(iii) On a 
$H^1(k,H) \buildrel \simeq \over \rightarrow \prod_{v \hskip1mm   r{\acute e}elle} H^1(k_v,H)$. }

\medskip

{\it D\'emonstration}
Soit $1\to H^{ss} \to H \to P \to 1$ la suite canonique attach\'ee
\`a $H$, avec $H^{ss}$ semi-simple simplement connexe et $P$
un $k$-tore quasi-trivial.

 \smallskip
 
L'assertion (i) r\'esulte imm\'ediatement de la suite exacte
d'ensembles point\'es
$$ H^1(k_v,H^{ss}) \to H^1(k_v,H) \to H^1(k_v,P),$$
o\`u $H^1(k_v,P)=0$ car $P$ est un tore quasi-trivial,
et $H^1(k_v,H^{ss})=1$ pour $v$ place non r\'eelle d'apr\`es
Kneser et Tits.

\smallskip

Soit $\Sigma$ un ensemble fini de places de $k$, dont on
peut supposer qu'il contient les places r\'eelles.
On a alors un diagramme commutatif de suites exactes d'ensembles point\'es
$$\diagram {H^{ss}(k) & \to & H(k) & \to & P(k) & \to & H^1(k,H^{ss})\cr
\downarrow{}{} &&\downarrow{}{} &&\downarrow{}{} && \downarrow{}{} \cr
\prod_{v \in
\Sigma}H^{ss}(k_v)& \to& \prod_{v \in \Sigma}H(k_v) &\to &\prod_{v \in
\Sigma}P(k_v) &\to&
\prod_{v \in \Sigma}H^1(k_v,H^{ss})}$$
D'apr\`es Kneser, Harder et Chernousov, on a une bijection d'ensembles finis
$$H^1(k,H^{ss}) \buildrel \simeq \over \rightarrow \prod_{v \hskip1mm   r{\acute e}elle} H^1(k_v,H^{ss}) =
\prod_{v \in \Sigma} H^1(k_v,H^{ss}).$$

\smallskip
D'apr\`es
Kneser et Platonov, 
tout
$k$-groupe semi-simple simplement connexe  satisfait
l'appro\-xi\-ma\-tion faible. Le tore quasi-trivial $P$ satisfait l'approxi\-mation
faible. Enfin chaque application $H(k_v) \to P(k_v)$ est ouverte,
puisque $H \to P$ est lisse. Une chasse au diagramme facile
montre alors que $H(k)$ est dense dans $\prod_{v \in \Sigma}H(k_v) $,
ce qui \'etablit le point (ii).

\smallskip

 Pour tout corps $F$ contenant $k$,
on a la suite  exacte d'ensembles point\'es
$$ P(F)  \to H^1(F,H^{ss}) \to H^1(F,H) \to H^1(F,P),$$
o\`u $P(F)$ agit sur $H^1(F,H^{ss})$. 
Comme $P$ est quasi-trivial,
on a $H^1(F,P)=0$. On a donc le diagramme de suites exactes
d'ensembles point\'es
$$  \diagram { P(k) & \to &H^1(k,H^{ss}) & \to & H^1(k,H) & \to & 0 \cr
\downarrow{}{}&&\downarrow{}{}&&\downarrow{}{}& &
\cr\prod_{v \hskip1mm   r{\acute e}elle}P(k_v)  &\to & \prod_{v \hskip1mm  
r{\acute e}elle}H^1(k_v,H^{ss})&
\to &\prod_{v \hskip1mm   r{\acute e}elle}H^1(k_v,H) &
\to &
0 }$$
Comme rappel\'e ci-dessus, la verticale m\'ediane est une bijection.
 En
utilisant la densit\'e de
$P(k)$  dans $\prod_{v \hskip1mm   r{\acute e}elle}P(k_v)$, on obtient le  point (iii).
\cqfd

L'\'enonc\'e suivant est une
variante  d'\'enonc\'es de P. Gille ([Gi97];  [Gi01]; appendice \`a [BoKu04])
et Borovoi et Kunyavski\u{\i} ([BoKu04, Thm. 8.4]).

\bigskip

{\bf Th\'eor\`eme  9.3} {\it Soit $k$ un corps de
nombres. 
Une r\'esolution flasque $1 \to S \to H \to G \to 1$ 
d'un $k$-groupe r\'eductif connexe $G$ induit une  suite exacte de
groupes 
finis
$$ H^{ss}(k)/R \to G(k)/R \to H^1(k,S) \to 1,$$
o\`u $H^{ss} \subset H$, le groupe  d\'eriv\'e de $H$,
est un groupe semi-simple
simplement connexe.
}

\medskip

{\it D\'emonstration}  On a la suite exacte
$$ G(k) \to H^1(k,S) \to H^1(k,H)$$
et les suites exactes locales
$$ G(k_v) \to H^1(k_v,S) \to H^1(k_v,H).$$
Pour $v$ place finie, $H^1(k_v,H)=1$ d'apr\`es la proposition pr\'ec\'edente.
Pour $v$ place r\'eelle, le $k_v$-tore flasque $S_{k_v}$ est quasi-trivial,
donc $H^1(k_v,S)=0$. Ainsi, pour toute place $v$, 
la fl\`eche
$H^1(k_v,S) \to H^1(k_v,H)$ a pour image l'\'el\'ement distingu\'e
de $H^1(k_v,H)$. Comme la fl\`eche $H^1(k,H) \to \prod_v H^1(k_v,H)$
a un noyau trivial (proposition pr\'ec\'edente), il s'en suit
que la fl\`eche $H^1(k,S) \to H^1(k,H)$ a une image r\'eduite 
\`a l'\'element distingu\'e de $H^1(k,H)$, et donc l'application
$G(k) \to H^1(k,S)$ est surjective.
D'apr\`es le th\'eor\`eme 8.1
on a donc une suite exacte de groupes
 $$ H(k)/R \to G(k)/R \to H^1(k,S) \to 1.$$

On trouve dans l'appendice de Gille \`a [BoKu04] (Theorem 1) une 
d\'emonstration du fait d\'elicat suivant  :
 sur un corps de nombres,
 pour toute suite
exacte de $k$-groupes r\'eductifs connexes
$$ 1 \to H^{ss} \to H \to T \to 1$$
avec $H^{ss}$ semi-simple simplement connexe et $T$ un $k$-tore,
on a une suite exacte induite
$$ H^{ss}(k)/R \to H(k)/R \to T(k)/R \to 1.$$
De la suite exacte
$$ 1 \to H^{ss} \to H \to P \to 1,$$ 
o\`u $P$ est un $k$-tore quasi-trivial, donc satisfait$P(k)/R=1$,
on d\'eduit  que l'application $H^{ss}(k)/R \to H(k)/R$ est surjective,
et donc la suite exacte de l'\'enonc\'e.

Il reste  \`a \'etablir la finitude des groupes intervenant dans cette suite.
La finitude de $H^1(k,S)$ pour $k$ un corps de nombres (et plus g\'en\'eralement un
corps de type fini sur le corps premier) et $S$ un $k$-tore flasque est connue
([CTSa77], Thm. 1 p. 192).
La finitude de $G(k)/R$ pour $G$ un $k$-groupe semi-simple simplement connexe
et $k$ un corps de nombres est une cons\'equence d'un th\'eor\`eme g\'en\'eral de Margulis.
En fait, il a \'et\'e \'etabli que pour un tel groupe,
$G(k)/R=1$ sauf peut-\^etre si $G$ contient un
facteur anisotrope de type $E_6$ (voir le livre de Platonov et
Rapinchuk  [PlRa] 
et l'article de Chernousov et
Timoshenko  [ChTi]). \qed

\bigskip

{\bf Th\'eor\`eme 9.4} {\it Soient $k$ un corps de nombres, $G$
un $k$-groupe r\'eductif connexe et
 $1 \to S \to H \to G \to 1$
une r\'esolution flasque de $G$. Alors

(i) Cette suite induit un isomorphisme entre le groupe
$A(G)$ qui mesure
le d\'efaut d'approxi\-mation faible pour $G$ et le
groupe  ab\'elien fini  $$ {\rm Coker} [H^1(k,S) \to \oplus_v
H^1(k_v,S)].$$

(ii)  Cette suite induit une bijection
de 
l'ensemble
$\X^1(k,G)$ avec le groupe ab\'elien fini $$\X^2(k,S)= {\rm Ker} \hskip1mm
[H^2(k,S) \to
\prod_v H^2(k_v,S)].$$

(iii) Cette suite induit une surjection
$$H^1(k,G) \to \Ker [H^2(k,S) \to H^2(k,P)].$$
Si $k$ est totalement imaginaire, cette application est une bijection.

(iv) (Kottwitz [Ko86], Borovoi [Bo98])
Cette suite induit une suite exacte d'ensembles point\'es
$$H^1(k,G) \to \oplus_{v} H^1(k_{v},G) \to \Hom(\Pic(G),\Q/\Z),$$
les fl\`eches $ H^1(k_{v},G) \to \Hom(\Pic(G),\Q/\Z)$
\'etant les fl\`eches compos\'ees
$$H^1(k_{v},G) \to \Hom(\Pic(G_{k_{v}}), \Q/\Z) \to \Hom(\Pic(G),\Q/\Z).$$

(v) (Sansuc [Sa81, Thm. 9.5])
 Soit $X$ une $k$-compactification lisse de $G$.
On a une suite exacte de groupes ab\'eliens finis
$$ 0 \to A(G) \to \Hom({\rm Br}(X)/{\rm Br}(k),\Q/\Z) \to \X^1(k,G) \to 0.$$
}

\medskip

{\it D\'emonstration} 

(i) Soit $K/k$ l'extension finie galoisienne d\'eployant
le $k$-tore $S$. Pour toute place finie $v$ non ramifi\'ee dans
$K$, le $k_v$-tore $S_{k_v}$ est d\'eploy\'e par une extension
cyclique, donc (Endo-Miyata, voir [CTSa77],  Prop. 2 p. 184) est un facteur direct d'un
$k_v$-tore quasi-trivial, donc $H^1(k_v,S)=0$.
Soit $\Sigma$ un ensemble fini de places de $k$  contenant les
places ramifi\'ees dans $K$ et les places r\'eelles.

 On consid\`ere le diagramme de suites exactes
$$ \diagram{H(k) &\to &G(k) &\to &H^1(k,S) &\to &H^1(k,H) \cr
\downarrow{}{}&&\downarrow{}{} &&\downarrow{}{} && \downarrow{}{}
\cr \prod_{v \in \Sigma} H(k_v)& \to & \prod_{v \in \Sigma} G(k_v) &\to &
\prod_{v \in
\Sigma} H^1(k_v,S)
&\to &\prod_{v \in \Sigma} H^1(k_v,H) .}$$
D'apr\`es la proposition 9.2 et la d\'emonstration du th\'eor\`eme 9.3, ce
diagramme se r\'eduit
\`a
$$ \diagram{H(k) &\to &G(k) &\to &H^1(k,S) &\to & 1 \cr
\downarrow{}{}&&\downarrow{}{} &&\downarrow{}{} && 
\cr \prod_{v \in \Sigma} H(k_v)& \to & \prod_{v \in \Sigma} G(k_v) &\to &
\prod_{v \in
\Sigma} H^1(k_v,S)
&\to &1}$$
et $H(k)$ est dense dans le produit $\prod_{v \in \Sigma} H(k_v)$.
Les groupes finis $H^1(k_v,S)$ sont nuls pour toute
place $v$ non ramifi\'ee dans  $K$. Ceci \'etablit l'\'enonc\'e (i).
Pour plus de d\'etails pour ce type d'argument qui remonte \`a Kneser, en particulier
pour le fait que l'adh\'erence de $G(k)$ dans  $ \prod_{v \in \Sigma} G(k_v) $
est un sous-groupe distingu\'e,
voir
[Sa81, \S 3].
A noter que la nullit\'e de $H^1(k_v,S)$ pour $v$ en dehors
des places non ramifi\'ees dans $K$ et l'argument ci-dessus assurent
 qu'en dehors de cet ensemble, l'approximation faible vaut pour $G$.

\medskip

(ii) 
La fl\`eche bord associ\'ee \`a la suite exacte
$1 \to S \to H \to G \to 1$
d\'efinit une application naturelle 
$\X^1(k,G) \to \X^2(k,S)$.
Si $\alpha \in \X^1(k,G) $ a une image triviale
dans $\X^2(k,S)$, il est l'image 
d'un \'el\'ement $\beta \in H^1(k,H)$ dont l'image
dans tout $H^1(k_v,H)$ provient de $H^1(k_v,S)$.
Aux places r\'eelles, le $k_v$-tore flasque $S_{k_v}$
est d\'eploy\'e par une extension cyclique. On
a donc $H^1(k_v,S)=0$ en une telle place.
Ainsi $\beta$ est dans le noyau de l'application
$H^1(k,H) \to \prod_{v \hskip1mm r{\acute e}elle} H^1(k_v,H)$,
et ce noyau est trivial (Prop. 9.2).
L'argument de torsion donn\'e par Sansuc 
([Sa81], preuve du th\'eor\`eme 4.3 pages
28 et 29) s'applique dans notre contexte et montre
que l'application est injective. 

D'apr\`es Kneser, Harder et Borovoi ([Bo98], Lemma 5.6.5) il existe
un $k$-tore maximal $T_{G^{sc}} \subset G^{sc} \subset H$ tel que $\X^2(k,T_{G^{sc}})=0$.
Soit $T _{H} \subset H$ le $k$-tore maximal qui est le centralisateur connexe
de $T_{G^{sc}}$ dans $H$. Cela d\'efinit un diagramme  du type pr\'ec\'edant la proposition A.1
ci-dessous. 
De la cohomologie de la suite exacte
$1 \to T_{G^{sc}} \to T_{H} \to P \to 1,$
o\`u $P$ est un tore quasi-trivial, et de $\X^2(k,T_{G^{sc}})=0$ on d\'eduit 
$\X^2(k,T_{H})=0$.
L'image de  tout \'el\'ement $\alpha \in \X^2(k,S)$ dans $H^2(k,T_{H})$ est donc nulle,
et tout tel \'el\'ement est dans l'image de $H^1(k,T_{G}) \to H^2(k,S)$, ce qui implique qu'il
 est dans  l'image de $H^1(k,G) \to H^2(k,S)$, et \'etablit le point (ii).

\medskip

(iii)  Soit $\alpha \in \ker [H^2(k,S) \to H^2(k,P)]$. Soit $\Sigma$ l'ensemble fini des places
de $k$ telles que $\alpha_{v }\neq 0 \in H^2(k_{v},S)$. D'apr\`es Kneser, Harder et Borovoi ([Bo98], Lemme 5.6.5), il existe un $k$-tore maximal $T_{G^{sc}} \subset G^{sc} \subset H$
tel que $H^2(k_{v},T_{G^{sc}})=0$ pour $v \in \Sigma$ et $\X^2(k,T_{G^{sc}})=0$.
Soit $T _{H} \subset H$ le $k$-tore maximal qui est le centralisateur connexe
de $T_{G^{sc}}$ dans $H$. Cela d\'efinit un diagramme  du type pr\'ec\'edant la proposition A.1
ci-dessous. De la cohomologie de la suite exacte
$1 \to T_{G^{sc}} \to T_{H} \to P \to 1,$
o\`u $P$ est un tore quasi-trivial, et de $\X^2(k,T_{G^{sc}})=0$ on d\'eduit comme ci-dessus
$\X^2(k,T_{H})=0$.
L'image de $\alpha  \in \ker [H^2(k,S) \to H^2(k,P)]$ dans $H^2(k,T_{H})$ 
appartient \`a $\X^2(k,T_{H})$, elle est donc nulle. Ainsi $\alpha \in H^2(k,S)$ est l'image d'un
\'el\'ement de $H^1(k,T_{G})$, a fortiori est-ce l'image d'un \'el\'ement de $H^1(k,G)$.
Que l'application soit bijective lorsque $k$ est totalement imaginaire est un cas particulier
du th\'eor\`eme 8.4.

\medskip

(iv) Pour tout $k$-tore $T$, un morceau de la suite de Tate-Nakayama ([MiADT], Chap. I, Thm. 4.20) est
$$ H^2(k,T) \to \oplus_{v} H^2(k_{v},T) \to \Hom((T^*)^{\frak{g}}, \Q/\Z) \to 0.$$
Du morphisme $S \to P$, de la proposition 3.3 et de   $\X^2(k,P)=0$
($P$ est quasi-trivial)
on d\'eduit donc la suite exacte
$$\Coker [H^2(k,S) \to H^2(k,P)] \to \oplus _{v}\Coker [H^2(k_{v},S) \to H^2(k_{v},P)]
\to \Hom(\Pic(G), \Q/\Z).$$
Soit $\xi_{v} \in \oplus_{v} H^1(k_{v},G)$ d'image triviale dans $\Hom(\Pic(G), \Q/\Z)$.
Soit $$\rho_{v } \in  \oplus _{v}\Coker [H^2(k_{v},S) \to H^2(k_{v},P)]$$ son image.
En combinant la suite exacte obtenue \`a l'instant, le th\'eor\`eme 9.1 et
le point (iii) ci-dessus, on trouve $\eta \in H^1(k,G)$ d'image $\xi_{v} \in H^1(k_{v},G)$
pour chaque place $v$ non r\'eelle et d'image $\rho_{v } \in \Coker [H^2(k_{v},S) \to H^2(k_{v},P)]$
pour $v$ place r\'eelle.
Pour l'argument de torsion d\'elicat permettant
de montrer que l'on peut trouver $\xi   \in H^1(k,G)$ d'image $\xi_{v}$ en toute place
de $k$, y compris aux places r\'eelles, je renvoie au th\'eor\`eme 5.11 de
 [Bo98].
 
\medskip

(v) Un autre morceau de la suite exacte de Tate-Nakayama ([MiADT], Chap. I, Thm. 4.20) pour le $k$-tore $S$
est
$$ H^1(k,S) \to \oplus_{v} H^1(k_{v},S) \to \Hom(H^1(k,S^*),\Q/\Z) \to H^2(k,S) \to 
\oplus_{v} H^2(k_{v},S).$$
D'apr\`es le th\'eor\`eme 7.1, on a 
 $H^1(k,S^*) \simeq H^1(k,{\rm
Pic}({\overline X}))
\simeq {\rm Br}(X)/{\rm Br}(k).$ 
L'\'enonc\'e (v) r\'esulte alors des \'enonc\'es (i) et (ii).
\cqfd

\medskip

{\it Remarques 9.4.1}

(1) Des isomorphismes comme aux \'enonc\'es 9.4 (i) et 9.4 (ii) sont \'enonc\'es
dans [BoKu04] (Thm. 8.11 et Thm. 8.16 (i)).

(2)  L'\'enonc\'e (iv) pour $G$ un $k$-tore $T$ est une cons\'equence imm\'ediate de
la dualit\'e de Tate-Nakayama et de l'isomorphisme bien connu $H^1(k,T^*) \simeq \Pic(T)$.
Lorsque $G$ est un $k$-groupe semi-simple, on  d\'eduit l'\'enonc\'e (iv) de la suite exacte de Poitou-Tate
pour la cohomologie des modules galoisiens finis et de la suite exacte $1 \to \mu \to G^{sc} \to G \to 1$. 
Dans le cas g\'en\'eral, cet \'enonc\'e est \'etabli,
sous une forme l\'eg\`erement diff\'erente,
par Kottwitz ([Ko86] (voir aussi [Bo98], \S 5). 
Le passage d'une formulation
\`a l'autre se fait en utilisant la proposition 6.3 et la remarque 6.3.1.
Cet \'enonc\'e a \'et\'e utilis\'e par Borovoi et Rudnick [BoRu95] et par Harari [Ha02].

(3) Dans (v), pour faire le lien avec l'obstruction de Brauer-Manin au principe de Hasse et
\`a l'approximation faible pour les espaces  homog\`enes principaux
de $G$, il
faudrait  suivre les  fl\`eches. Pour l'approxima\-tion faible cela ne doit
pas \^etre trop d\'elicat, dans la mesure o\`u l'on peut partir
de la r\'esolution flasque de $G$ donn\'ee par la restriction
d'un torseur universel trivial en $e_{G} \in G(k)$, et le lien
entre \'evaluation des torseurs et obstruction de Brauer-Manin est fait dans 
[CTSa87b].  Le cas de l'obstruction au principe de Hasse pour les compactifications
d'espaces  homog\`enes  principaux pourrait \^etre nettement plus
d\'elicat, comme on peut voir en suivant la d\'emonstration de
Sansuc [Sa81] ou la pr\'esentation de Skorobogatov ([Sk01], Chapitre 6.2).

\bigskip

\bigskip

{\bf Appendice A : Comparaison avec les travaux de Borovoi et de Borovoi-Kunyavski\u{\i}}

\bigskip

Soit
$$1  \to S \to H \to G  \to 1$$
une r\'esolution flasque du $k$-groupe r\'eductif $G$, 
avec $H$ extension d'un $k$-tore
quasi-trivial $P$ par le $k$-groupe semi-simple
simplement connexe $H^{ss}$.
Comme expliqu\'e au \S 3, une telle r\'esolution donne  naissance \`a un  diagramme commutatif de suites exactes
de $k$-groupes lin\'eaires
$$\diagram {&    & 1 && 1 && 1
\cr
  &    & \downarrow{}{} &&  \downarrow{}{} &&  \downarrow{}{}
&&    
\cr 
1 &\to &\mu &\to &H^{ss} &\to & G^{ss} &\to& 1
 \cr
 &    &  \downarrow{}{} &&  \downarrow{}{} &&  \downarrow{}{} &&     
   \cr
 1 & \to & S & \to & H & \to &G &\to  &  1 \cr
&    &  \downarrow{}{} &&  \downarrow{}{} &&  \downarrow{}{} &&   
    \cr
 1 &  \to & M    &\to & P &
\to & T & \to & 1
\cr
&    &  \downarrow{}{} &&  \downarrow{}{} &&  \downarrow{}{}
\cr
&    & 1 && 1 && 1.
}
$$
Dans ce diagramme, $P=H^{tor}$, $T=G^{tor}$, et $M$ est le $k$-groupe de type multiplicatif
noyau de la fl\`eche naturelle $P \to T$ induite par $H \to G$. 
Le $k$-groupe $H^{ss}$ s'identifie au rev\^etement simplement connexe $G^{sc}$ du groupe semi-simple $G^{ss}$, groupe d\'eriv\'e de $G$.

Fixons un tore maximal $T_{G}$ dans $G$.
En prenant des images r\'eciproques, et en  utilisant le
fait qu'un groupe
alg\'ebrique
extension d'un tore alg\'ebrique par un tore
alg\'ebrique est un tore alg\'ebrique (rappel {\bf 0.7}), on obtient 
 un diagramme de suites
exactes de $k$-groupes de type multiplicatif
$$\diagram {&    & 1 && 1 && 1
\cr
  &    & \downarrow{}{} &&  \downarrow{}{} &&  \downarrow{}{}
&&    
\cr 
1 &\to &\mu &\to &T_{G^{sc}}  &\to & T_{G^{ss}}  &\to& 1
 \cr
 &    &  \downarrow{}{} &&  \downarrow{}{} &&  \downarrow{}{} &&     
   \cr
 1 & \to & S & \to & T_H & \to & T_G &  \to &  1 \cr
&    &  \downarrow{}{} &&  \downarrow{}{} &&  \downarrow{}{} &&   
    \cr
 1 &  \to & M    &\to & P &
\to & T & \to & 1
\cr
&    &  \downarrow{}{} &&  \downarrow{}{} &&  \downarrow{}{}
\cr
&    & 1 && 1 && 1.
}
$$

Notons  que chacun des $k$-tores maximaux $T_{G} \subset G$,  $T_{H}\subset H$,
$T_{G^{ss}} \subset G^{ss}$,  $T_{H^{ss}}Ê\subset H^{ss}$ d\'etermine les autres.
Par exemple $T_{G}$ est le centralisateur de $T_{G^{ss}}$ dans $G$.

Du diagramme ci-dessus on tire en particulier les suites
exactes de $k$-groupes de type multiplicatif
$$ 1 \to \mu \to T_{G^{sc}} \to T_G \to T \to 1$$
et
$$ 1 \to \mu \to S \to P \to T \to 1.$$

On a les suites duales de groupes de caract\`eres.

\bigskip

\vfill\eject

{\bf Proposition A.1} {\it  Le complexe ci-dessus induit

(i) un diagramme commutatif  de suites exactes de $k$-groupes de type multiplicatif
$$\diagram{
1 & \to& \mu & \to &T_{G^{sc}}  & \to &T_{G}  &\to &T& \to& 1  \cr
&&\downarrow{=}  &&\downarrow{} &&\downarrow{} &&\downarrow{=} &&\cr
1  & \to &\mu & \to &T_{H} &\to &P\oplus T_{G}& \to &T& \to &1 \cr
&&\uparrow{=} &&\uparrow{} &&\uparrow{} &&\uparrow{=} &&\cr
1  & \to& \mu & \to & S& \to& P &\to& T& \to &1 ,
}$$
(ii) un diagramme commutatif de suites exactes de modules galoisiens de type fini
$$\diagram{
0 & \to& T^* & \to &T_{G}^* & \to & T_{G^{sc}}^*&\to &\mu^*& \to& 0  \cr
&&\uparrow{=}  &&\uparrow{} &&\uparrow{} &&\uparrow{=} &&\cr
0  & \to &T^* & \to &  P^*\oplus T_{G} ^*   &\to &      T_{H}^*   & \to &\mu^*& \to &0 \cr
&&\downarrow{=} &&\downarrow{} &&\downarrow{} &&\downarrow{=} &&\cr
0  & \to& T^* & \to &P^*  & \to& S^*  &\to& \mu^*& \to &0 ,}$$
(iii) un diagramme commutatif de suites exactes de modules galoisiens de type fini 
$$\diagram{
0 &\to &T_{G^{sc}*}  &\to& T_{G*}  &\to &R_{1}& \to& 0\cr
&&\downarrow&& \downarrow&& \downarrow{\simeq} &\cr
0 &\to &T_{H*} &\to &P_{*}\oplus T_{G*} &\to& R_{2}& \to& 0\cr
&&\uparrow&& \uparrow&& \uparrow{\simeq}&&\cr
0 &\to &S_{*}& \to &P_{*} & \to& R_{3}& \to& 0 .
}$$

\medskip

En particulier,

(iv) le complexe de $k$-tores  $[T_{G^{sc}} \to T_G]$  est naturellement quasi-isomorphe 
au
complexe de $k$-tores $[S
\to P]$,

(v) le complexe de modules galoisiens $[T_G^* \to T_{G^{sc}}^*]$
est naturellement quasi-isomorphe au complexe de modules galoisiens $[P^* \to S^*]$,

(vi)  le complexe de modules galoisiens  $[T_{G^{sc}*} \to T_{G*}]$
est naturellement quasi-isomorphe au complexe de modules galoisiens $[S_{*}\to P_{*}]$.
}
 
\medskip
Pour ces complexes de longueur 2, on place le terme de gauche en degr\'e
$-1$ et le terme de droite en degr\'e $0$. Le complexe de $k$-tores 
$[T_{G^{sc}} \to T_G]$ est le complexe
utilis\'e par Borovoi ([Bo96] et [Bo98]).
\medskip
\vfill\eject
{\it D\'emonstration}   
Il s'agit d'un exercice g\'en\'eral d'alg\`ebre
homologique. Supposons que l'on ait un diagramme de suites
exactes
$$\diagram {&    & 0 && 0 && 0\cr
  &    & \uparrow{}{} && \uparrow{}{} && \uparrow{}{}
&&    \cr 
0 &  \to & A_1    &\to & A_2 &
\to & A_3 & \to & 0 \cr
 &    & \uparrow{}{} && \uparrow{}{} && \uparrow{}{} &&       \cr
 0 & \to & B_1 & \to & B_2 & \to & B_3 &\to  &  0 \cr
&    & \uparrow{}{} && \uparrow{}{} && \uparrow{}{} &&       \cr
 0 &\to &C_1 &\to &C_2 &\to & C_3 &\to& 0
\cr
&    & \uparrow{}{} && \uparrow{}{} && \uparrow{}{}\cr
&    & 0 && 0 && 0
}
$$
 On a alors le diagramme commutatif
$$\diagram { 0 & \to & C_1 & \to & B_1 & \to &  A_2 & \to & A_3&  \to & 0\cr
&& \downarrow{}{=} && \downarrow{}{} && \downarrow{}{} && \downarrow{}{=}&\cr
 0 & \to &  C_1  & \to &  B_2  & \to  & A_2 \oplus B_3 &  \to  & A_3  & \to & 
0\cr
&& \uparrow{}{=} && \uparrow{}{} && \uparrow{}{} && \uparrow{=}{}&\cr
 0  & \to &  C_1  & \to  & C_2  & \to  & B_3  & \to  & A_3  & \to &  0}$$
Dans les suites horizontales, les fl\`eches sont les fl\`eches \'evidentes, \`a l'exception de
 la fl\`eche $(A_2 \oplus B_3) \to A_3$ dans la suite horizontale m\'ediane, qui 
est d\'efinie comme la diff\'erence ($(a,b) \mapsto f(a)-g(b)$) des deux fl\`eches 
$f : A_2 \to A_3$ et $g : B_3\to A_3$, et de la fl\`eche $B_{3} \to A_{3}$ dans la suite
horizontale inf\'erieure, qui est d\'efinie comme l'oppos\'ee de la fl\`eche donn\'ee $g : B_{3} \to A_{3}$.
On envoie la premi\`ere suite horizontale dans la seconde
par les applications \'evidentes et l'application
$A_2 \to A_2 \oplus B_3$ donn\'ee par $a \mapsto (a,0)$.
On envoie la troisi\`eme suite horizontale  dans la seconde
par les fl\`eches \'evidentes et  et l'application
$B_3 \to A_2 \oplus B_3$ donn\'ee par $b \mapsto (0,b)$.
On obtient ainsi un quasi-isomorphisme de complexes
$B_1 \to A_2$ vers $B_2 \to A_2 \oplus B_3$
et un quasi-isomorphisme de 
$ C_2 \to B_3$ vers $B_2 \to A_2 \oplus B_3$.
Ainsi les complexes $B_1 \to A_2$ et $ C_2 \to B_3$
sont quasi-isomorphes.

Ceci \'etablit l'\'enonc\'e (i). L'\'enonc\'e (ii) en r\'esulte imm\'ediatement, le passage au groupe des caract\`eres $M \mapsto \Hom_{\k-gp}({\overline M},\G_{m,\k})$ transformant suite exacte de $k$-groupes de type multiplicatif en suite exacte de modules galoisiens. On laisse au lecteur de soin d'\'etablir l'\'enonc\'e (iii) \`a partir de l'\'enonc\'e (ii). 
\cqfd

\bigskip

{\bf Groupe fondamental alg\'ebrique d'un groupe lin\'eaire et groupe de Brauer d'une compactification }

\medskip

Dans  [Bo96] et  [Bo98], Borovoi d\'efinit le groupe fondamental d'un groupe alg\'ebrique r\'eductif connexe
\`a partir du   diagramme de groupes de type multiplicatif
pr\'ec\'edant la proposition A.1 ci-dessus.
Sa d\'efinition est :

$$\pi_{1}^{Bor}(G) = {\rm Coker} [T_{G^{sc} *} \to T_{G *}].$$
Il d\'emontre que ce module galoisien ne d\'epend pas du choix du tore maximal $T_{G}$.

\bigskip

{\bf Proposition A.2} {\it Le  module galoisien $\pi_{1}^{Bor}(G)$ est naturellement  isomorphe au 
module galoisien $\pi_{1}(G)$ introduit au \S 6 ci-dessus.}

\medskip

{\it D\'emonstration} C'est une application de la proposition A.1 (iii).\cqfd

\bigskip

Soit $X$ une $k$-compactification lisse du groupe $G$.
Dans [BoKu00], sur un corps $k$ de caract\'eristique z\'ero, Borovoi et Kunyavski\u{\i} \'etablissent  des isomorphismes :
 $${\rm Br}(X)/{\rm Br}(k) \simeq  \X^1_{\omega}(k,[T_G^* \to T_{H^{ss}}^*])$$
 (dans le complexe de longueur 2, le terme de gauche est en degr\'e $-1$, celui de droite en
degr\'e $0$),
et
$${\rm Br}(X)/{\rm Br}(k) \simeq {\cyr X}^1_{\omega}(k,\Hom_{\Z}(\pi^{Bor}_{1}(G),\Q/\Z  )).$$
Pour ce faire,
ils utilisent les r\'esultats de [CTKu98] et recourent aux $z$-extensions. 

Le th\'eor\`eme 7.2 et la proposition  A.2 du pr\'esent article donnent imm\'ediatement une autre d\'emons\-tration 
de l'existence d'un 
isomorphisme entre les deux groupes de la seconde formule. Montrons comment l'existence d'un isomorphisme entre les deux groupes de la premi\`ere formule  r\'esulte aussi du pr\'esent article. D'apr\`es la proposition A.1, le complexe $[T_G^* \to
T_{H^{ss}}^*]$ et le complexe
$[P^* \to S^*]$ sont  quasi-isomorphes. Le groupe $\X^1_{\omega}(k,[T_G^* \to T_{H^{ss}}^*])$
est donc isomorphe au groupe $\X^1_{\omega}(k,[P^* \to S^*])$.
Il suffit donc de
montrer
$$ H^1(k,S^*) \simeq \X^1_{\omega}(k,[P^* \to S^*]).$$
On a la suite exacte courte de complexes :
$$ 0 \to [0 \to S^*] \to [P^* \to S^*] \to [P^* \to 0] \to 0.$$
Cette suite induit une suite exacte longue
$$H^1(k,P^*) \to H^1(k,S^*) \to H^1(k,[P^* \to S^*]) \to H^2(k,P^*)$$
et les suites similaires pour la cohomologie de tout sous-groupe ferm\'e
procyclique $\h$ de $\g={\rm Gal}(\k/k)$. Notons que $P^*$
\'etant de permutation, on a  $H^1(\g,P^*)=0=H^1(\h,P^*)$
et $\X^2_{\omega}(k,P^*)=0$. Par ailleurs, $S^*$
\'etant flasque, on a $H^1(\h,S^*)=0$.
Si l'on consid\`ere le diagramme commutatif de suites
exactes obtenu par la restriction de $\g$
\`a tous ses sous-groupes procycliques $\h$, le lemme des
5 donne imm\'ediatement
$$ H^1(k,S^*) \buildrel \simeq \over \rightarrow \X^1_{\omega}(k,[P^* \to S^*]). $$
D'apr\`es la proposition 7.1, on a $ {\rm Br}(X)/{\rm Br}(k) \simeq H^1(k,S^*) $.
Ceci ach\`eve la d\'emonstration.
\cqfd

\bigskip

{\bf Cohomologie ab\'elianis\'ee de Borovoi}

\medskip

Dans [Bo96] et [Bo98], pour $G$ un $k$-groupe r\'eductif, et $i=0,1$,
Borovoi d\'efinit des groupes de cohomologie galoisienne ab\'elianis\'es $H^{i}_{ab}(k,G)$.
Ce sont les groupes d'hypercohomologie 
$$H^{0}_{ab}(k,G)=H^{0}(k,[T_{G^{sc}} \to T_G])$$
et
$$H^{ 1}_{ab}(k,G)=H^{1}(k,[T_{G^{sc}} \to T_G])$$
du
complexe  de $k$-tores de longueur 2  $$[T_{G^{sc}} \to T_G],$$
o\`u les notations sont celles de la proposition A.1 ci-dessus, et  les tores sont plac\'es en degr\'e $-1$ et $0$.

Borovoi d\'efinit un homomorphisme
$$G(k) \to H^0_{ab}(k,G)$$ et une application
$$H^1(k,G) \to H^1_{ab}(k,G).$$

Montrons comment le point de vue adopt\'e ici permet de d\'efinir des fl\`eches analogues.
Soit $G$ un $k$-groupe r\'eductif connexe.
Soit
$1  \to S \to H  \to G  \to 1$
une r\'esolution flasque de $G$, 
avec $H$ extension d'un $k$-tore
quasi-trivial $P$ par un $k$-groupe semi-simple
simplement connexe. 
Comme indiqu\'e au \S 8, ces donn\'ees fournissent un homomorphisme
$$ G(k) \to {\rm Ker} [H^1(k,S) \to H^1(k,P)]$$
et une fl\`eche
$$ H^1(k,G) \to {\rm Ker} [H^2(k,S) \to H^2(k,P)].$$

D'apr\`es la proposition A.1 (iv), le complexe de  
 $k$-tores  $[T_{G^{sc}} \to T_G]$  est naturellement quasi-isomorphe 
au
complexe de $k$-tores $[S
\to P]$, on a donc
$$H^{i}(k,[S \to P]) \simeq H^{i}(k,[T_{G^{sc}} \to T_G]).$$

On a la suite exacte
$$ H^1(k,P) \to H^1(k, [S \to P]) \to H^2(k,S) \to H^2(k,P)$$
et comme $H^1(k,P)=0$ (Hilbert 90) on voit que la fl\`eche bord
$$H^1(k,G) \to \Ker [H^2(k,S) \to H^2(k,P)]$$ se transcrit en une fl\`eche
$$H^1(k,G) \to H^1(k, [S \to P]) $$
et donc en une fl\`eche
$$H^1(k,G) \to H^1(k,[T_{G^{sc}} \to T_G]))$$
soit encore
$H^1(k,G) \to  H^1_{ab}(k,G).$
On laisse au lecteur le soin d'\'etudier son lien avec l'application d\'efinie par Borovoi.

\medskip
Au niveau $H^0$, on a la suite exacte
$$ H^0(k,S) \to H^0(k,P) \to H^0(k, [S \to P]) \to H^1(k,S) \to H^1(k,P).$$
Borovoi d\'efinit une fl\`eche $G(k) \to H^0(k,[T_{G^{sc}} \to T_G]))$,
donc une fl\`eche $G(k) \to H^{0}(k, [S \to P])$.

On doit sans mal v\'erifier que cette fl\`eche et la 
fl\`eche \'evidente $G(k) \to H^1(k,S)$ sont compatibles.
Admettant ce point, la fl\`eche $G(k) \to H^0(k, [S \to P])$
contient plus d'information que
la fl\`eche bord $G(k) \to H^1(k,S)$ (laquelle passe au quotient par la $R$-\'equivalence).

Typiquement, si $G$ est semi-simple et $G^{sc}$ est son rev\^etement
simplement connexe, alors la fl\`eche de Borovoi contr\^ole
le quotient de $G(k)$ par l'image de $G^{sc}(k)$.

De fait, on a la suite exacte 
$$ G^{sc}(k) \to G(k) \to H^0_{ab}(k,G) \to H^1(k, G^{sc}) \to H^1(k, G) \to
H^1_{ab}(k,G).$$
(Borovoi, [Bo 96]  p. 405;  [Bo98]    (3.10.1) p. 24)
 Ainsi un \'el\'ement de $G(k)$ est dans le noyau de
la fl\`eche $G(k) \to H^0(k, [S \to P])$ si et seulement s'il
est dans l'image de $G^{sc}(k)$.

\bigskip
 
{\bf  Appendice B. Comparaison de deux complexes}

\bigskip

{\bf Lemme B.1} {\it Soient $k$ un corps, $X$ une $k$-vari\'et\'e lisse int\`egre,
$T$ un $k$-tore d\'eploy\'e, $\pi: Y \to X$ un $T$-torseur, $\lambda: T^* \to \pic(X)$
son type. Soient $U \subset X$ un ouvert non vide et $V =\pi^{-1}(U) \subset Y$.

(i) La fl\`eche $\pi^* : \Div_{X\setminus U} \to \Div_{Y\setminus V}$
est un isomorphisme.

(ii) On a une fl\`eche  naturelle $k[V]^{\times}  \to \Div(Y \setminus V)$
associant \`a une fonction son diviseur. La fibre g\'en\'erique $Y_{\eta}$ de
$\pi$ est un $T$-torseur trivial sur ${\rm Spec}(k(X))$, et le quotient du groupe des fonctions inversibles
sur $Y_{\eta}$ par le groupe multiplicatif $k(X)^{\times} $ est canoniquement isomorphe \`a $T^*$.
Ceci d\'efinit une fl\`eche $k[V]^{\times} \to T^*$.

(iii) Le diagramme suivant est anticommutatif
$$\diagram{k[V]^{\times}  & \to & T^* \cr
\downarrow & & \downarrow{}{\lambda} \cr
\Div(Y \setminus V) & & \cr
{\pi^*}\uparrow{\simeq}&&
 \cr
\Div(X \setminus U)& \to & \pic(X).
}$$}

{\it D\'emonstration} 
L'\'enonc\'e (i) est clair, puisque toutes les fibres 
g\'eom\'etriques de $\pi$ sont g\'eom\'etriquement int\`egres et non vides.
L'\'enonc\'e (ii) est un rappel de la version torseur du lemme de Rosenlicht
([CTSa87b], Prop. 1.4.2 p. 383), 
dans le cas simple d'un torseur trivial. 
Le tore $T$ est un produit $\prod_{i=1}^n\G_{m,k}$.
 La classe du $T$-torseur $Y \to X$
dans $H^1(X,T)=\oplus_{i=1}^n H^1(X,\G_m)$ s'\'ecrit $(L_1, \dots, L_n)$.
Soit, pour chaque $i$, $Y_i \to X$ un $\G_m$-torseur de classe $L_i \in 
\pic(X)$. On a un isomorphisme $Y \simeq Y_1\times_X\times \dots \times_XY_n$.
Pour \'etablir le r\'esultat, on peut supposer que $X$
 poss\`ede un $k$-point $p$ (si ce n'est pas le cas, on s'y ram\`ene en consid\'erant
 $X \times_{k}k(X)$ et en utilisant le fait que la fl\`eche naturelle
 $\pic(X) \to \pic(X\times_{k}k(X))$ est injective).
  On peut alors fixer un $k$-point $p_i \in Y_i(k)$ au-dessus de $p$.
Ceci permet de d\'efinir des immersions ferm\'ees $Y_i \hookrightarrow Y$.
On v\'erifie  que l'assertion du point (iii) est fonctorielle en
de tels morphismes. Ceci permet de ramener la d\'emonstration du point (iii)
au cas $T=\G_m={\rm Spec} k[t,1/t]$. 

Pour \'etablir ce r\'esultat, tr\`es certainement classique, on peut restreindre $U$, en particulier supposer
que le $\G_{m}$-torseur $Y \to X$ a une section  $U\to V =\pi^{-1}(U)$.
Identifions $V=\G_{m,U}$. 
On a $k[V]^{\times}/k[U]^{\times} \buildrel \simeq \over \rightarrow \Z$ la fl\`eche envoyant la classe de $t$ sur $1 \in \Z$.
L'image de $1 \in \Z$ dans $\pic (X)=H^1(X,\G_{m})$ est la classe du torseur $Y \to X$.
Le diviseur de la fonction rationnelle $t$ sur $Y$ a son support \'etranger \`a $V$,
il est donc l'image r\'eciproque d'un diviseur $\Delta$ sur $X$.
Soit ${\rm div}_{Y}(t) = \pi^*(\Delta)$.
Soit $U_{i}, i\in I$ un recouvrement ouvert de $X$, avec $U_{i_{0}}=U$
tel sur chaque ouvert $U_{i}$ on puisse \'ecrire $\Delta$ comme le diviseur d'une fonction
rationnelle $f_{i} \in k(X)$. Sur  $U_{i_{0}}=U$, choisissons $f_{i_{0}}=1$.
Soit $V_{i}=\pi^{-1}(U_{i})$. Sur $V_{i}$, le quotient $v_{i}=f_{i}/t  $ appartient
\`a $k[V_{i}]^{\times}$.
La classe de $-\Delta$ dans $\pic(X)=H^1(X,O_{X}^{\times}) $
est donn\'ee par le 1-cocycle $ u_{ij} =f_{j}/f_{i} $.
L'image r\'ecriproque de ce cocycle sur $Y$ est $tv_{j}/tv_{i}=v_{j}/v_{i}$.
On dispose des fl\`eches $\rho_{i } :V_{i} \to U_{i}\times \G_{m}$
donn\'ees par $y \mapsto (\pi(y),v_{i}(y))$.
Pour $i,j$ donn\'es, la fl\`eche
$U_{i}\times \G_{m} \to U_{j}\times \G_{m}$ donn\'ee par
$(x,z)  \mapsto (x,(f_{j}/f_{i})(x)z)$ est compatible avec
les fl\`eches $\rho_{i}$. On obtient ainsi un morphisme de $X$-sch\'emas
$Y \to Z$, o\`u $Z$ est le sch\'ema  obtenu par le recollement via les 
fl\`eches $f_{j}/f_{i}$, qui est donc le $\G_{m}$-torseur associ\'e \`a $-\Delta$.
Le $X$-sch\'ema $Y$ est un $\G_{m}$-torseur, le $X$-sch\'ema $Z$ est aussi
un $\G_{m}$-torseur sur $X$.
C'est le $\G_{m}$-torseur associ\'e au diviseur $\Delta$.
La restriction de $Y\to Z$ au-dessus de $U=U_{i_{0}}$ est l'identit\'e
et respecte l'action de $\G_{m}$.
Ceci implique que la fl\`eche $Y \to Z$ est un isomorphisme
de $\G_{m}$-torseurs. \qed

\bigskip

{\bf Proposition B.2} {\it 
 Soit $G$ un $k$-groupe r\'eductif connexe, soit $T=G^{tor}$
 et $\mu$ le noyau de $G^{sc} \to G^{ss}$.
 Soit $X$ une
$k$-compactification lisse de $G$. 

(i) On a   la suite
exacte de modules galoisiens
$$ 0 \to G^* \to {\rm Div}_{{\overline X} \setminus {\overline G}}({\overline X}) \to {\rm Pic}({\overline X})
\to {\rm Pic}({\overline G}) \to 0.$$

(ii) Le module $G^*$ s'identifie \`a $T^*$. Le module
${\rm Pic}({\overline G})$ s'identifie \`a ${\rm Pic}({\overline G}^{ss} )$,
c'est-\`a-dire \`a ${\mu}^*$.

(iii)
Soit $1 \to S \to H \to G \to 1$ une r\'esolution flasque de $G$, et soit $P=H^{tor}$.
Le complexe de modules galoisiens $[{\rm Div}_{{\overline X} \setminus {\overline G}}({\overline X}) \to {\rm
Pic}({\overline X})]$ est quasi-isomorphe au complexe de modules galoisiens $[P^* \to S^*]$.

(iv)  Avec les notations de l'appendice A, le complexe de modules galoisiens $[{\rm Div}_{{\overline X} \setminus {\overline G}}({\overline X}) \to {\rm
Pic}({\overline X})]$ est quasi-isomorphe au complexe de modules galoisiens
 $[T_G^* \to T_{G^{sc}}^*]$.}

\bigskip

{\it D\'emonstration} Comme $G$ est un ouvert de la $k$-vari\'et\'e projective, lisse et g\'eom\'etriquement int\`egre $X$, on a la suite exacte
$$ 0 \to {\overline k}[G]^{\times}/{\overline k}^{\times}  \to {\rm Div}_{{\overline X} \setminus {\overline G}}({\overline X}) \to {\rm Pic}({\overline X})
\to {\rm Pic}({\overline G}) \to 0.$$
L'\'enonc\'e (i) r\'esulte alors du lemme de Rosenlicht. L'\'enonc\'e (ii) 
est un cas particulier de l'\'enonc\'e {\bf 0.3}, appliqu\'e \`a la suite exacte de $k$-groupes
$ 1 \to G^{ss} \to G \to T \to 1,$ et de l'isomorphisme $\pic({\overline G}^{ss})={\mu}^*$.
Comme $S$ est un $k$-tore flasque, 
il existe un torseur $\pi : Y\to X$ de groupe structural $S$, de type 
\'etendant le torseur de groupe $S$ donn\'e par $H \to G$
([CTSa77, Prop. 9]).
Soit $\lambda : S^* \to \pic({\overline X})$ son type.
On a alors le diagramme commutatif de suites exactes suivant
$$\diagram { 
0  &\to &G^* & \to &  {\rm Div}_{{\overline X} \setminus {\overline G}}({\overline X})&  \to  &{\rm Pic}({\overline X}) & \to & {\rm Pic}({\overline G})&  \to & 0 \cr
&&\uparrow{=}{} &&\uparrow{}{}&&\uparrow{-\lambda}{}&&\uparrow{(-1)}{} & &\cr
0  &\to &G^* & \to &  H^*  &  \to  & S^* & \to & {\rm Pic}({\overline G})&  \to &0
}$$
La suite horizontale sup\'erieure est celle de l'\'enonc\'e (i). La suite horizontale inf\'erieure est donn\'ee par l'\'enonc\'e {\bf 0.3}, joint au fait que $\pic({\overline H})=0$.
La fl\`eche non \'evidente $H^* \to {\rm Div}_{{\overline X} \setminus {\overline G}}({\overline X})$
est la compos\'ee de la fl\`eche diviseur $H^* \to {\rm Div}_{{\overline Y} \setminus {\overline H}}({\overline Y})$ et de l'inverse
de l'isomorphisme $\pi^* : {\rm Div}_{{\overline X} \setminus {\overline G}}({\overline X}) \to {\rm Div}_{{\overline Y} \setminus {\overline H}}({\overline Y})$.
La commutativit\'e du carr\'e de gauche est alors \'evidente. Celle du carr\'e m\'edian est assur\'ee par
le lemme B.1. Celle du carr\'e de droite est une fonctorialit\'e \'evidente.
On a l'isomorphisme naturel $P^* \buildrel \simeq \over \rightarrow H^*$. On voit donc que le complexe 
$[P^* \to S^*]$ est quasi-isomorphe au complexe  $[{\rm Div}_{{\overline X} \setminus {\overline G}}({\overline X}) \to {\rm
Pic}({\overline X})]$, ce qui est l'assertion (iii).  

L'\'enonc\'e (iv) r\'esulte alors de la proposition A.1.
\qed

\bigskip

{\it Remarques B.2.1} 

(1) Soient $X$ une $k$-compactification lisse de $G$, $S^*=\pic({\overline X})$ et $\pi : Y  \to X$ le torseur universel sur
$X$ de fibre triviale au point $e \in G(k)$ (ce torseur est bien d\'efini \`a isomorphisme de  $S$-torseurs pr\`es).
On a montr\'e au \S 5 que l'on peut munir $H=Y\times_{X}G$ d'une structure de $k$-groupe quasi-trivial,
la projection $H \to G$ \'etant un  homomorphisme de noyau $S$.  Dans ce cas, les fl\`eches
$H^* \to {\rm Div}_{{\overline X} \setminus {\overline G}}({\overline X})$ et $\lambda : S^* \to \pic({\overline X})$ sont des isomorphismes :
dans le diagramme ci-dessus, toutes les fl\`eches sont des isomorphismes, le complexe
$[P^* \to S^*]$ est isomorphe au complexe  $[{\rm Div}_{{\overline X} \setminus {\overline G}}({\overline X}) \to {\rm
Pic}({\overline X})]$.

(2) Soit de fa\c con g\'en\'erale $X_{c}$ une $k$-compactification lisse d'une $k$-vari\'et\'e lisse
g\'eom\'e\-tri\-quement int\`egre $X$. Comme me l'a fait remarquer Skorobogatov, le diagramme
$$\diagram {
\k(X)^{\times}/\k^{\times} & \to & \Div({\overline X}) \cr
\uparrow & & \uparrow    \cr
 \k(X)^{\times}/\k^{\times} \oplus  {\rm Div}_{{\overline X_{c}} \setminus {\overline X}}   ({\overline X}_{c})
 & \to & \Div( {\overline X}_{c})    \cr
 \downarrow & & \downarrow \cr
 {\rm Div}_{{\overline X_{c}}  \setminus {\overline X}}     ({\overline X}_{c}) & \to & \pic({\overline X}_{c})
 }, $$
 o\`u les fl\`eches verticales de gauche sont les projections \'evidentes,
 la fl\`eche horizontale associe \`a $(f, \Delta)$ le diviseur ${\rm div}_{{\overline X}_{c}}(f) -\Delta$, et les autres fl\`eches
 sont \'evidentes, induit des quasi-isomorphismes entre les complexes horizontaux.
 
(3) La remarque pr\'ec\'edente et la proposition B.2 (iv) redonnent
 certains r\'esultats 
de Borovoi et van Hamel [BovH06].

\bigskip

{\bf Bibliographie}

\medskip

[Bo96] M. Borovoi, Abelianization of the first Galois cohomology
of reductive groups, IMRN~{\bf 8} (1996), 401--407.

[Bo98] M. Borovoi, Abelian Galois cohomology of reductive groups,
Mem. AMS, Vol.  132, Number 626 (1998).

[BoRu95] M. Borovoi et Z. Rudnick, Hardy-Littlewood varieties and semisimple groups,
Invent. Math. {\bf 119} (1995) 37--66.

[BoKu00] M. Borovoi et B. Kunyavski\u{\i},  Formulas for the
unramified Brauer group of a principal homogeneous space
of a linear algebraic group, J. of Algebra {\bf 225} (2000) 804--821.

[BoKu04] M. Borovoi et B. Kunyavski\u{\i}, Arithmetical birational
invariants of linear algebraic groups over two-dimensional
geometric fields, avec un appendice de P. Gille, J. of Algebra {\bf 276}
(2004) 292--339.

[BovH06] M. Borovoi et J. van Hamel, Extended Picard  complexes for algebraic groups
and homogeneous spaces,    C. R. Acad. Sci. Paris, S\'er. I  {\bf 342} (2006) 671--674.

[ChTi] V. I. Chernousov  et L. M. Timoshenko,  
 Sur le groupe des classes de $R$-\'equivalence des
groupes semi-simples sur des corps arithm\'etiques (en russe), Algebra i Analiz  {\bf 11}  (1999) 
 191--221;  trad. ang.,  St. Petersburg Math. J.  {\bf 11}  (2000),  
1097--1121.

[CT04] J.-L. Colliot-Th\'el\`ene, R\'esolutions flasques
des groupes r\'eductifs connexes, C. R. Acad. Sc. Paris,
 S\'erie I {\bf 339} (2004)
331--334.

[CT05] J.-L. Colliot-Th\'el\`ene, 
Un th\'eor\`eme de finitude pour le groupe de
Chow des z\'ero-cycles d'un groupe alg\'ebrique
lin\'eaire  sur un corps $p$-adique,  Invent. math. {\bf 159} (2005) 589--606.

[CTGiPa04] J.-L. Colliot-Th\'el\`ene, P. Gille et R. Parimala,
Arithmetic of linear algebraic groups over 2-dimensional
geometric fields, Duke Math. J. {\bf 121} (2004) 285--341.

[CTHaSk05] J.-L. Colliot-Th\'el\`ene, D. Harari et A. N. Skorobogatov,
Expositiones mathematicae.
Compactification \'equivariante d'un tore (d'apr\`es Brylinski et K\"unnemann),  
Expositiones mathematicae {\bf 23} (2005) 161-170.

[CTKu98] J.-L. Colliot-Th\'el\`ene et B. Kunyavski\u{\i},
Groupe de Brauer non ramifi\'e des espaces principaux homog\`enes
des groupes lin\'eaires, J. Ramanujan Math. Soc. {\bf 13} (1998)
37--49.

[CTSa77] J.-L. Colliot-Th\'el\`ene et J.-J. Sansuc, La $R$-\'equivalence sur
les tores, Ann. Sc. \'Ec. Norm. Sup. {\bf 10} (1977) 175--229.

[CTSa87a]  J.-L. Colliot-Th\'el\`ene et J.-J. Sansuc,
Principal homogeneous spaces under flasque tori: Applications,
J. Algebra {\bf 106} (1987) 148--205.

[CTSa87b] J.-L. Colliot-Th\'el\`ene et J.-J. Sansuc,  La descente sur les
vari\'et\'es rationnelles, II, Duke Math. J. {\bf 54} (1987) 375--492.

[dJ04]  A. J. de Jong, The period-index problem for the Brauer group of an algebraic surface,
Duke Math. J. {\bf 123} (2004) 71--94.

[Gi97] P. Gille, La $R$-\'equivalence sur les groupes alg\'ebriques
r\'eductifs d\'efinis sur un corps global,
Publications Math\'ematiques de l'I.H.\'E.S. {\bf 86} (1997) 
199--235.

[Gi01] P. Gille, Cohomologie galoisienne des groupes quasi-d\'eploy\'es sur des corps de dimension cohomologique $\leq 2$,  Compositio mathematica {\bf    125}
(2001) 283--325.

[Ha02] D. Harari, Groupes alg\'ebriques et points rationnels,
Math. Annalen {\bf 322} (2002) 811--826.

 [Ko84]  R. E. Kottwitz,
 Stable trace formula: cuspidal tempered terms,  Duke Math. J.  {\bf 51}  (1984),  
 611--650.
 
[Ko86] R. E. Kottwitz, Stable trace formula: elliptic singular terms, 
Math. Annalen {\bf 275} (1986) 365--399.

[Me93] A. S. Merkur'ev, Generic element in $SK_{1}$ for simple algebras,
K-Theory {\bf 7} (1993) 1--3.

[Me96] A. S. Merkur'ev,  $R$-equivalence and rationality problem for semi-simple adjoint classical algebraic groups.  Inst. Hautes \'Etudes Sci. Publ. Math. {\bf 84} (1996), 189--213 (1997).

[Me98] A. S. Merkur'ev, $K$-Theory and algebraic groups, in Proceedings
of the European Congress of Mathematicians (Budapest), Progress in
Mathematics {\bf 169}, p. 43--72, Birkh\"auser Verlag, 1998.

[MiADT] J. S. Milne, Arithmetic duality theorems, Perspective in Mathematics {\bf 1},
Academic Press (1986). Edition corrig\'ee (2004) disponible sur la page personnelle
de l'auteur.

[Mi-Sh82] J. S. Milne et K.-y. Shih, Conjugates of Shimura varieties,
in {\it Hodge Cycles, Motives and Shimura Varieties}, Springer L. N. M. {\bf 900}
(1982).

[Mi72] M. Miyanishi, On the algebraic fundamental group of an
algebraic group, J. Math. Kyoto Univ. {\bf 12-2} (1972) 351--367.

[PlRa] V.P. Platonov et A. S. Rapinchuk, Algebraic groups and number theory,
Academic Press, Boston, 1994.

[Sa81] J.-J. Sansuc,
Groupe de Brauer et arithm\'etique des groupes alg\'ebriques
lin\'eaires, J. f\"ur die reine und angew. Math. (Crelle) {\bf 327} (1981)
12--80.

[SGA3]  M. Demazure et A. Grothendieck,     
S\'eminaire de G\'eom\'etrie Alg\'ebrique du Bois Marie 1962/64, Sch\'emas en Groupes  I, II et III. Dirig\'e par M. Demazure et A. Grothendieck. Springer L.N.M. { \bf 151}, {\bf 152} et  {\bf 153} (1970).

[Se59]  J-P. Serre, 
Groupes alg\'ebriques et corps de classes, Actualit\'es scientifiques et
industrielles {\bf 1264}, Publications de l'Institut Math\'ematique
de Nancago VII, Hermann Paris, 1959.

[SeCG] J-P. Serre,  Cohomologie galoisienne,  cinqui\`eme \'edition, r\'evis\'ee et compl\'et\'ee,
Springer L. N. M. {\bf 5} (1973, 1994).

[Sk01] A. N. Skorobogatov, Torsors and rational points. Cambridge Tracts in Mathematics {\bf 144}. Cambridge University Press, Cambridge, 2001.

[Vo77] V. E. Voskresenski\u{\i}, Algebraicheskie tory, Nauka, Moscou, 1977.

[Vo98]  V. E. Voskresenski\u{\i}, Algebraic groups and their birational
invariants, Transl. Mathe\-ma\-tical Monographs {\bf 179}, Amer. Math.
Soc., 1998.

\vskip2cm
\vfill\eject

Jean-Louis Colliot-Th\'el\`ene

C.N.R.S., Math\'ematiques,

UMR 8628,

B\^atiment 425,

Universit\'e Paris-Sud,

F-91405 Orsay

FRANCE

colliot@math.u-psud.fr

\bye